\documentclass[10pt]{amsart}
\usepackage{amssymb,eucal}
\numberwithin{equation}{section}
\newtheorem{claim}{}[section]

\newtheorem{theorem}[claim]{Theorem}
\newtheorem{lemma}[claim]{Lemma}
\newtheorem{proposition}[claim]{Proposition}
\newtheorem{corollary}[claim]{Corollary}

\newcommand{\lup}{\begin{picture}(.4,.2)
\linethickness{.2mm}
\put(0.055,.13){$\scriptstyle{\to}$}
\put(.08,0){\line(0,1){.18}}
\end{picture}}

\newcommand{\rup}{\!\!\begin{picture}(.4,.2)
\linethickness{.2mm}
\put(-.01,.13){$\scriptstyle{\,\,\leftarrow}$}
\put(.33,0){\line(0,1){.18}}
\end{picture}}

\newcommand{\subrup}{\!\!\begin{picture}(.4,.2)
\linethickness{.2mm}
\put(.035,.15){$\scriptscriptstyle{\,\,\leftarrow}$}
\put(.33,0){\line(0,1){.18}}
\end{picture}}

\newcommand{\rdown}{\begin{picture}(.4,.2)
\linethickness{.2mm}
\put(0.055,-.080){$\scriptstyle{\to}$}
\put(.07,-.010){\line(0,1){.18}}
\end{picture}}

\newcommand{\subrdown}{\begin{picture}(.4,.2)
\linethickness{.2mm}
\put(0.055,-.050){$\scriptscriptstyle{\to}$}
\put(.07,-.010){\line(0,1){.18}}
\end{picture}\!\!}

\newcommand{\ldown}{\!\begin{picture}(.4,.2)
\linethickness{.2mm}
\put(-.01,-.08){$\scriptstyle{\,\,\leftarrow}$}
\put(.34,-.01){\line(0,1){.18}}
\end{picture}}

\newcommand{\subldown}{\!\!\begin{picture}(.4,.2)
\linethickness{.2mm}
\put(0,-.04){$\scriptscriptstyle{\,\,\leftarrow}$}
\put(.3,0){\line(0,1){.18}}
\end{picture}}

\numberwithin{equation}{section}
\newcounter{bean}

\catcode`\@=11                                                                  
                                                                                
\newcount\cnta\newcount\cntb\newcount\cntc

\def\showlasttree{%
     \g\cnta\count\l@stdiminfo
     \g\cntb\cnta
     \g\advance\cntb 5 
     \g\advance\cntb \count\l@sttreeheight
     \g\advance\cntb \count\l@sttreeheight
     \ifnum\count\l@sttreeheight=-1\relax
           \g\advance\cntb by 2
           \immediate\write16{Tree contour for dummy node:}
      \else\immediate\write16{Tree contour:}% 
        \fi
     \for\cntc:=\cnta\to\cntb\do\immediate\write16{\the\dimen\cntc}\od}

\def\lplain{lplain}                     % Set \LaTeXtrue if TreeTeX is          
\newif\ifLaTeX                          % used together with LaTeX,             
\ifx\fmtname\lplain\LaTeXtrue           % otherwise set \LaTeXfalse             
    \else\LaTeXfalse\fi                 % (LaTeX defines \fmtname=={lplain}).   

\immediate\write16{This is TreeTeX, Version 2.1, for use with \ifLaTeX LaTeX%     
                   \else plain TeX\fi.}                                         
                                                                                
\ifLaTeX \let\lineseg\line              % latex_picture is part of latex.tex,   
   \else \let\@line\line                % so you don't need it if you use       
\catcode`\@=11                                                                  
                                                                                
\def\@height{height}                                                            
\def\@depth{depth}                                                              
\def\@width{width}                                                              
                                                                                
\font\tenln=line10                                                              
\font\tencirc=circle10                                                          
\font\tenlnw=linew10                                                            
\font\tencircw=circlew10                                                        
                                                                                
\newcount\@tempcnta                                                             
\newcount\@tempcntb                                                             
\newdimen\@tempdima                                                             
\newdimen\@tempdimb                                                             
\newbox\@tempboxa                                                               
                                                                                
\def\@whilenoop#1{}                                                             
                                                                                
\def\@whiledim#1\do #2{\ifdim #1\relax#2\@iwhiledim{#1\relax#2}\fi}             
\def\@iwhiledim#1{\ifdim #1\let\@nextwhile=\@iwhiledim                          
        \else\let\@nextwhile=\@whilenoop\fi\@nextwhile{#1}}                     
                                                                                
\def\@ifnextchar#1#2#3{\let\@tempe #1\def\@tempa{#2}\def\@tempb{#3}\futurelet   
    \@tempc\@ifnch}                                                             
\def\@ifnch{\ifx \@tempc \@sptoken \let\@tempd\@xifnch                          
      \else \ifx \@tempc \@tempe\let\@tempd\@tempa\else\let\@tempd\@tempb\fi    
      \fi \@tempd}                                                              
                                                                                
\def\:{\let\@sptoken= } \:  % this makes \@sptoken a space token                
                                                                                
\def\:{\@xifnch} \expandafter\def\: {\futurelet\@tempc\@ifnch}                  
                                                                                
\def\@ifstar#1#2{\@ifnextchar *{\def\@tempa*{#1}\@tempa}{#2}}                   
                                                                                
\let\:=\>

\newdimen\@wholewidth                                                           
\newdimen\@halfwidth                                                            
\newdimen\unitlength \unitlength =1pt                                           
\newbox\@picbox                                                                 
\newdimen\@picht                                                                
                                                                                
\def\picture(#1,#2){\@ifnextchar({\@picture(#1,#2)}{\@picture(#1,#2)(0,0)}}     
                                                                                
\def\@picture(#1,#2)(#3,#4){\@picht #2\unitlength                               
\setbox\@picbox\hbox to #1\unitlength\bgroup                                    
\hskip -#3\unitlength \lower #4\unitlength \hbox\bgroup}                        
                                                                                
\def\endpicture{\egroup\hss\egroup\ht\@picbox\@picht                            
\dp\@picbox\z@\leavevmode\box\@picbox}                                          
                                                                                
\long\def\put(#1,#2)#3{\@killglue\raise#2\unitlength\hbox to \z@{\hskip         
#1\unitlength #3\hss}\ignorespaces}                                             
                                                                                
\long\def\multiput(#1,#2)(#3,#4)#5#6{\@killglue\@multicnt=#5\relax              
\@xdim=#1\unitlength                                                            
\@ydim=#2\unitlength                                                            
\@whilenum \@multicnt > 0\do                                                    
{\raise\@ydim\hbox to \z@{\hskip                                                
\@xdim #6\hss}\advance\@multicnt \m@ne\advance\@xdim                            
#3\unitlength\advance\@ydim #4\unitlength}\ignorespaces}                        
                                                                                
\def\@killglue{\unskip\@whiledim \lastskip >\z@\do{\unskip}}                    
                                                                                
\def\thinlines{\let\@linefnt\tenln \let\@circlefnt\tencirc                      
  \@wholewidth\fontdimen8\tenln \@halfwidth .5\@wholewidth}                     
\def\thicklines{\let\@linefnt\tenlnw \let\@circlefnt\tencircw                   
  \@wholewidth\fontdimen8\tenlnw \@halfwidth .5\@wholewidth}                    
                                                                                
\def\linethickness#1{\@wholewidth #1\relax \@halfwidth .5\@wholewidth}          
                                                                                
\def\shortstack{\@ifnextchar[{\@shortstack}{\@shortstack[c]}}                   
                                                                                
\def\@shortstack[#1]{\leavevmode                                                
\vbox\bgroup\baselineskip-1pt\lineskip 3pt\let\mb@l\hss                         
\let\mb@r\hss \expandafter\let\csname mb@#1\endcsname\relax                     
\let\\\@stackcr\@ishortstack}                                                   
                                                                                
\def\@ishortstack#1{\halign{\mb@l ##\unskip\mb@r\cr #1\crcr}\egroup}

\def\@stackcr{\@ifstar{\@ixstackcr}{\@ixstackcr}}                               
\def\@ixstackcr{\@ifnextchar[{\@istackcr}{\cr\ignorespaces}}                    
                                                                                
\def\@istackcr[#1]{\cr\noalign{\vskip #1}\ignorespaces}

\newif\if@negarg                                                                
                                                                                
\def\line(#1,#2)#3{\@xarg #1\relax \@yarg #2\relax                              
\@linelen=#3\unitlength                                                         
\ifnum\@xarg =0 \@vline                                                         
  \else \ifnum\@yarg =0 \@hline \else \@sline\fi                                
\fi}                                                                            
                                                                                
\def\@sline{\ifnum\@xarg< 0 \@negargtrue \@xarg -\@xarg \@yyarg -\@yarg         
  \else \@negargfalse \@yyarg \@yarg \fi                                        
\ifnum \@yyarg >0 \@tempcnta\@yyarg \else \@tempcnta -\@yyarg \fi               
\ifnum\@tempcnta>6 \@badlinearg\@tempcnta0 \fi                                  
\setbox\@linechar\hbox{\@linefnt\@getlinechar(\@xarg,\@yyarg)}%                 
\ifnum \@yarg >0 \let\@upordown\raise \@clnht\z@                                
   \else\let\@upordown\lower \@clnht \ht\@linechar\fi                           
\@clnwd=\wd\@linechar                                                           
\if@negarg \hskip -\wd\@linechar \def\@tempa{\hskip -2\wd\@linechar}\else       
     \let\@tempa\relax \fi                                                      
\@whiledim \@clnwd <\@linelen \do                                               
  {\@upordown\@clnht\copy\@linechar                                             
   \@tempa                                                                      
   \advance\@clnht \ht\@linechar                                                
   \advance\@clnwd \wd\@linechar}%                                              
\advance\@clnht -\ht\@linechar                                                  
\advance\@clnwd -\wd\@linechar                                                  
\@tempdima\@linelen\advance\@tempdima -\@clnwd                                  
\@tempdimb\@tempdima\advance\@tempdimb -\wd\@linechar                           
\if@negarg \hskip -\@tempdimb \else \hskip \@tempdimb \fi                       
\multiply\@tempdima \@m                                                         
\@tempcnta \@tempdima \@tempdima \wd\@linechar \divide\@tempcnta \@tempdima     
\@tempdima \ht\@linechar \multiply\@tempdima \@tempcnta                         
\divide\@tempdima \@m                                                           
\advance\@clnht \@tempdima                                                      
\ifdim \@linelen <\wd\@linechar                                                 
   \hskip \wd\@linechar                                                         
  \else\@upordown\@clnht\copy\@linechar\fi}                                     
                                                                                
\def\@hline{\ifnum \@xarg <0 \hskip -\@linelen \fi                              
\vrule \@height \@halfwidth \@depth \@halfwidth \@width \@linelen               
\ifnum \@xarg <0 \hskip -\@linelen \fi}                                         
                                                                                
\def\@getlinechar(#1,#2){\@tempcnta#1\relax\multiply\@tempcnta 8                
\advance\@tempcnta -9 \ifnum #2>0 \advance\@tempcnta #2\relax\else              
\advance\@tempcnta -#2\relax\advance\@tempcnta 64 \fi                           
\char\@tempcnta}                                                                
                                                                                
\def\vector(#1,#2)#3{\@xarg #1\relax \@yarg #2\relax                            
\@linelen=#3\unitlength                                                         
\ifnum\@xarg =0 \@vvector                                                       
  \else \ifnum\@yarg =0 \@hvector \else \@svector\fi                            
\fi}                                                                            
                                                                                
\def\@hvector{\@hline\hbox to 0pt{\@linefnt                                     
\ifnum \@xarg <0 \@getlarrow(1,0)\hss\else                                      
    \hss\@getrarrow(1,0)\fi}}                                                   
                                                                                
\def\@vvector{\ifnum \@yarg <0 \@downvector \else \@upvector \fi}               
                                                                                
\def\@svector{\@sline                                                           
\@tempcnta\@yarg \ifnum\@tempcnta <0 \@tempcnta=-\@tempcnta\fi                  
\ifnum\@tempcnta <5                                                             
  \hskip -\wd\@linechar                                                         
  \@upordown\@clnht \hbox{\@linefnt  \if@negarg                                 
  \@getlarrow(\@xarg,\@yyarg) \else \@getrarrow(\@xarg,\@yyarg) \fi}%           
\else\@badlinearg\fi}                                                           
                                                                                
\def\@getlarrow(#1,#2){\ifnum #2 =\z@ \@tempcnta='33\else                       
\@tempcnta=#1\relax\multiply\@tempcnta \sixt@@n \advance\@tempcnta              
-9 \@tempcntb=#2\relax\multiply\@tempcntb \tw@                                  
\ifnum \@tempcntb >0 \advance\@tempcnta \@tempcntb\relax                        
\else\advance\@tempcnta -\@tempcntb\advance\@tempcnta 64                        
\fi\fi\char\@tempcnta}                                                          
                                                                                
\def\@getrarrow(#1,#2){\@tempcntb=#2\relax                                      
\ifnum\@tempcntb < 0 \@tempcntb=-\@tempcntb\relax\fi                            
\ifcase \@tempcntb\relax \@tempcnta='55 \or                                     
\ifnum #1<3 \@tempcnta=#1\relax\multiply\@tempcnta                              
24 \advance\@tempcnta -6 \else \ifnum #1=3 \@tempcnta=49                        
\else\@tempcnta=58 \fi\fi\or                                                    
\ifnum #1<3 \@tempcnta=#1\relax\multiply\@tempcnta                              
24 \advance\@tempcnta -3 \else \@tempcnta=51\fi\or                              
\@tempcnta=#1\relax\multiply\@tempcnta                                          
\sixt@@n \advance\@tempcnta -\tw@ \else                                         
\@tempcnta=#1\relax\multiply\@tempcnta                                          
\sixt@@n \advance\@tempcnta 7 \fi\ifnum #2<0 \advance\@tempcnta 64 \fi          
\char\@tempcnta}

\def\@vline{\ifnum \@yarg <0 \@downline \else \@upline\fi}                      
                                                                                
\def\@upline{\hbox to \z@{\hskip -\@halfwidth \vrule \@width \@wholewidth       
   \@height \@linelen \@depth \z@\hss}}                                         
                                                                                
\def\@downline{\hbox to \z@{\hskip -\@halfwidth \vrule \@width \@wholewidth     
   \@height \z@ \@depth \@linelen \hss}}                                        
                                                                                
\def\@upvector{\@upline\setbox\@tempboxa\hbox{\@linefnt\char'66}\raise          
     \@linelen \hbox to\z@{\lower \ht\@tempboxa\box\@tempboxa\hss}}             
                                                                                
\def\@downvector{\@downline\lower \@linelen                                     
      \hbox to \z@{\@linefnt\char'77\hss}}

\def\dashbox#1(#2,#3){\leavevmode\hbox to \z@{\baselineskip \z@%                
\lineskip \z@%                                                                  
\@dashdim=#2\unitlength%                                                        
\@dashcnt=\@dashdim \advance\@dashcnt 200                                       
\@dashdim=#1\unitlength\divide\@dashcnt \@dashdim                               
\ifodd\@dashcnt\@dashdim=\z@%                                                   
\advance\@dashcnt \@ne \divide\@dashcnt \tw@                                    
\else \divide\@dashdim \tw@ \divide\@dashcnt \tw@                               
\advance\@dashcnt \m@ne                                                         
\setbox\@dashbox=\hbox{\vrule \@height \@halfwidth \@depth \@halfwidth          
\@width \@dashdim}\put(0,0){\copy\@dashbox}%                                    
\put(0,#3){\copy\@dashbox}%                                                     
\put(#2,0){\hskip-\@dashdim\copy\@dashbox}%                                     
\put(#2,#3){\hskip-\@dashdim\box\@dashbox}%                                     
\multiply\@dashdim 3                                                            
\fi                                                                             
\setbox\@dashbox=\hbox{\vrule \@height \@halfwidth \@depth \@halfwidth          
\@width #1\unitlength\hskip #1\unitlength}\@tempcnta=0                          
\put(0,0){\hskip\@dashdim \@whilenum \@tempcnta <\@dashcnt                      
\do{\copy\@dashbox\advance\@tempcnta \@ne }}\@tempcnta=0                        
\put(0,#3){\hskip\@dashdim \@whilenum \@tempcnta <\@dashcnt                     
\do{\copy\@dashbox\advance\@tempcnta \@ne }}%                                   
\@dashdim=#3\unitlength%                                                        
\@dashcnt=\@dashdim \advance\@dashcnt 200                                       
\@dashdim=#1\unitlength\divide\@dashcnt \@dashdim                               
\ifodd\@dashcnt \@dashdim=\z@%                                                  
\advance\@dashcnt \@ne \divide\@dashcnt \tw@                                    
\else                                                                           
\divide\@dashdim \tw@ \divide\@dashcnt \tw@                                     
\advance\@dashcnt \m@ne                                                         
\setbox\@dashbox\hbox{\hskip -\@halfwidth                                       
\vrule \@width \@wholewidth                                                     
\@height \@dashdim}\put(0,0){\copy\@dashbox}%                                   
\put(#2,0){\copy\@dashbox}%                                                     
\put(0,#3){\lower\@dashdim\copy\@dashbox}%                                      
\put(#2,#3){\lower\@dashdim\copy\@dashbox}%                                     
\multiply\@dashdim 3                                                            
\fi                                                                             
\setbox\@dashbox\hbox{\vrule \@width \@wholewidth                               
\@height #1\unitlength}\@tempcnta0                                              
\put(0,0){\hskip -\@halfwidth \vbox{\@whilenum \@tempcnta < \@dashcnt           
\do{\vskip #1\unitlength\copy\@dashbox\advance\@tempcnta \@ne }%                
\vskip\@dashdim}}\@tempcnta0                                                    
\put(#2,0){\hskip -\@halfwidth \vbox{\@whilenum \@tempcnta< \@dashcnt           
\relax\do{\vskip #1\unitlength\copy\@dashbox\advance\@tempcnta \@ne }%          
\vskip\@dashdim}}}\@makepicbox(#2,#3)}

\newif\if@ovt                                                                   
\newif\if@ovb                                                                   
\newif\if@ovl                                                                   
\newif\if@ovr                                                                   
\newdimen\@ovxx                                                                 
\newdimen\@ovyy                                                                 
\newdimen\@ovdx                                                                 
\newdimen\@ovdy                                                                 
\newdimen\@ovro                                                                 
\newdimen\@ovri                                                                 
                                                                                
\def\@getcirc#1{\@tempdima #1\relax \@tempcnta\@tempdima                        
  \@tempdima 4pt\relax \divide\@tempcnta\@tempdima                              
  \ifnum \@tempcnta > 10\relax \@tempcnta 10\relax\fi                           
  \ifnum \@tempcnta >\z@ \advance\@tempcnta\m@ne                                
    \else \@warning{Oval too small}\fi                                          
  \multiply\@tempcnta 4\relax                                                   
  \setbox \@tempboxa \hbox{\@circlefnt                                          
  \char \@tempcnta}\@tempdima \wd \@tempboxa}                                   
                                                                                
\def\@put#1#2#3{\raise #2\hbox to \z@{\hskip #1#3\hss}}                         
                                                                                
\def\oval(#1,#2){\@ifnextchar[{\@oval(#1,#2)}{\@oval(#1,#2)[]}}                 
                                                                                
\def\@oval(#1,#2)[#3]{\begingroup\boxmaxdepth \maxdimen                         
  \@ovttrue \@ovbtrue \@ovltrue \@ovrtrue                                       
  \@tfor\@tempa :=#3\do{\csname @ov\@tempa false\endcsname}\@ovxx               
  #1\unitlength \@ovyy #2\unitlength                                            
  \@tempdimb \ifdim \@ovyy >\@ovxx \@ovxx\else \@ovyy \fi                       
  \@getcirc \@tempdimb                                                          
  \@ovro \ht\@tempboxa \@ovri \dp\@tempboxa                                     
  \@ovdx\@ovxx \advance\@ovdx -\@tempdima \divide\@ovdx \tw@                    
  \@ovdy\@ovyy \advance\@ovdy -\@tempdima \divide\@ovdy \tw@                    
  \@circlefnt \setbox\@tempboxa                                                 
  \hbox{\if@ovr \@ovvert32\kern -\@tempdima \fi                                 
  \if@ovl \kern \@ovxx \@ovvert01\kern -\@tempdima \kern -\@ovxx \fi            
  \if@ovt \@ovhorz \kern -\@ovxx \fi                                            
  \if@ovb \raise \@ovyy \@ovhorz \fi}\advance\@ovdx\@ovro                       
  \advance\@ovdy\@ovro \ht\@tempboxa\z@ \dp\@tempboxa\z@                        
  \@put{-\@ovdx}{-\@ovdy}{\box\@tempboxa}%                                      
  \endgroup}                                                                    
                                                                                
\def\@ovvert#1#2{\vbox to \@ovyy{%                                              
    \if@ovb \@tempcntb \@tempcnta \advance \@tempcntb by #1\relax               
	\kern -\@ovro \hbox{\char \@tempcntb}\nointerlineskip                          
    \else \kern \@ovri \kern \@ovdy \fi                                         
    \leaders\vrule width \@wholewidth\vfil \nointerlineskip                     
    \if@ovt \@tempcntb \@tempcnta \advance \@tempcntb by #2\relax               
	\hbox{\char \@tempcntb}%                                                       
    \else \kern \@ovdy \kern \@ovro \fi}}                                       
                                                                                
\def\@ovhorz{\hbox to \@ovxx{\kern \@ovro                                       
    \if@ovr \else \kern \@ovdx \fi                                              
    \leaders \hrule height \@wholewidth \hfil                                   
    \if@ovl \else \kern \@ovdx \fi                                              
    \kern \@ovri}}                                                              
                                                                                
\def\circle{\@ifstar{\@dot}{\@circle}}                                          
\def\@circle#1{\begingroup \boxmaxdepth \maxdimen \@tempdimb #1\unitlength      
   \ifdim \@tempdimb >15.5pt\relax \@getcirc\@tempdimb                          
      \@ovro\ht\@tempboxa                                                       
     \setbox\@tempboxa\hbox{\@circlefnt                                         
      \advance\@tempcnta\tw@ \char \@tempcnta                                   
      \advance\@tempcnta\m@ne \char \@tempcnta \kern -2\@tempdima               
      \advance\@tempcnta\tw@                                                    
      \raise \@tempdima \hbox{\char\@tempcnta}\raise \@tempdima                 
        \box\@tempboxa}\ht\@tempboxa\z@ \dp\@tempboxa\z@                        
      \@put{-\@ovro}{-\@ovro}{\box\@tempboxa}%                                  
   \else  \@circ\@tempdimb{96}\fi\endgroup}                                     
                                                                                
\def\@dot#1{\@tempdimb #1\unitlength \@circ\@tempdimb{112}}                     
                                                                                
\def\@circ#1#2{\@tempdima #1\relax \advance\@tempdima .5pt\relax                
   \@tempcnta\@tempdima \@tempdima 1pt\relax                                    
   \divide\@tempcnta\@tempdima                                                  
   \ifnum\@tempcnta > 15\relax \@tempcnta 15\relax \fi                          
   \ifnum \@tempcnta >\z@ \advance\@tempcnta\m@ne\fi                            
   \advance\@tempcnta #2\relax                                                  
   \@circlefnt \char\@tempcnta}

\thinlines                                                                      
                                                                                
\newcount\@xarg                                                                 
\newcount\@yarg                                                                 
\newcount\@yyarg                                                                
\newcount\@multicnt                                                             
\newdimen\@xdim                                                                 
\newdimen\@ydim                                                                 
\newbox\@linechar                                                               
\newdimen\@linelen                                                              
\newdimen\@clnwd                                                                
\newdimen\@clnht                                                                
\newdimen\@dashdim                                                              
\newbox\@dashbox                                                                
\newcount\@dashcnt

         \let\lineseg\line              % has the command \line for geometric   
         \let\line\@line                % lines, and plain TeX has the same     
   \fi                                  % command for lines of text. Because    
\catcode`\@=11                                                                  
\let\g\global                                                                   
\def\gxdef{\global\xdef}                                                        
                                                                                
\def\newcount{\alloc@0\count\countdef\insc@unt}                                 
                                                                                
\def\for#1:=#2\to#3\do#4\od{%
   \def\f@rcount{#1}\def\upp@rlimit{#3}\def\b@dy{#4}\f@rcount=#2\relax\dof@r}

\def\dof@r{\ifnum\f@rcount>\upp@rlimit\relax\let\n@xt\relax
                 \else\b@dy\advance\f@rcount\@ne\let\n@xt\dof@r\fi
           \n@xt}
                                                                                
\newcount\@xcount
\newcount\t@mes                                                                    
                                                                                
\def\ex#1\times#2\xe{%                                                          
   \@xcount1 \t@mes#1\def\b@dy{#2}\do@x}

\def\do@x{\ifnum\@xcount>\t@mes\let\n@xt\relax
                \else\b@dy\advance\@xcount\@ne\let\n@xt\do@x\fi
          \n@xt}
                                                                              
\newskip\thickn@ss
\newskip\@nner
\newskip\@uter

\def\rect@ngle#1#2#3{\hbox to 0pt{%
     \thickn@ss#3%
     \g\@nner#2\g\advance\@nner-\thickn@ss
     \g\divide\@nner\tw@
     \g\@uter#2\g\advance\@uter\thickn@ss
     \g\divide\@uter\tw@
     \hskip 0pt minus .5fil%
     \vrule height\@uter depth\@nner width\thickn@ss
     \vrule height\@uter depth-\@nner width#1%
     \hskip 0pt minus 1fil%
     \vrule height-\@nner depth\@uter width#1%
     \vrule height\@nner depth\@uter width\thickn@ss
     \hskip 0pt minus .5fil%
     }% \hbox
     }% \def

\def\s@ries#1#2{%                                                               
     \g\t@mpcnta1                                                               
     \gdef\t@mp{#1}%                                                            
     \@ssign#2/\l@st  % \l@st is a sentinal element                                                        
     \expandafter\gdef\csname#1\endcsname##1{%                                  
                      \csname#1\romannumeral##1\endcsname}%                     
     }                                                                          
                                                                                
\def\@ssign#1/#2{%
      \expandafter\gdef\csname\t@mp\romannumeral\t@mpcnta\endcsname{#1}%    
      \g\advance\t@mpcnta\@ne                                            
      \ifx#2\l@st                                                               
          \g\let\n@xt\relax                                                     
     \else\g\let\n@xt\@ssign                                                    
       \fi                                                                      
      \n@xt}
                                                                                 
\newdimen\leftdist                                                              
\newdimen\rightdist                                                             
\newbox\TeXTree                                                                 
                                                                                
\newcount\sl@pe                                                                 
\newcount\l@vels                                                                
\newcount\s@ze                                                                  
                                                                                
\newbox\circleb@x                                                               
\newbox\squareb@x                                                               
\newbox\dotb@x                                                                  
\newbox\triangleb@x
\newbox\textb@x                                                             
\newbox\frameb@x
                                                                                
\newdimen\circlew@dth                                                           
\newdimen\squarew@dth                                                           
\newdimen\dotw@dth                                                              
\newdimen\trianglew@dth                                                         
\newdimen\textw@dth                                                             
\newdimen\framew@dth
                                                         
\newdimen\vd@st                                                                 
\newdimen\hd@st                                                                 
\newdimen\based@st                                                              
\newdimen\dummyhalfcenterdim@n                                                  
                                                                                
\newcount\t@mpcnta                                                              
\newcount\t@mpcntb                                                              
\newcount\t@mpcntc                                                              
\newcount\t@mpcntd                                                              
\newdimen\t@mpdima                                                              
\newdimen\t@mpdimb                                                              
\newdimen\t@mpdimc                                                              
\newbox\t@mpboxa                                                                
\newbox\t@mpboxb                                                                
                                                                                
\newbox\leftb@x                                                                 
\newbox\rightb@x                                                                
\newbox\centerb@x                                                               
\newbox\beneathb@x                                                              
\newtoks\typ@                                                                   
\newbox\centerb@@x                                                              
\newdimen\centerdim@n                                                           
\newdimen\halfcenterdim@n                                                       
                                                                                
\newdimen\mins@p                                                                
\newdimen\halfmins@p                                                            
\newdimen\tots@p                                                                
\newdimen\halftots@p                                                            
\newdimen\currs@p                                                               
\newdimen\adds@p                                                                
\newcount\l@ftht                                                                
\newcount\r@ghtht                                                               
\newcount\l@ftinfo                                                              
\newcount\r@ghtinfo                                                             
\newbox\l@ftbox                                                                 
\newbox\r@ghtbox                                                                
                                                                                
\newif\ifr@ghthigher  % true iff the right subtree is higher than the left one  
\newif\ifadds@p                                                                 
                                                                                
\newcount\@larg                                                                 
\newcount\@rarg                                                                 
                                                                                
\newif\ifl@fttop                                                                
\newif\ifl@ftonly                                                               
\newif\ifr@ghtonly                                                              
\newif\if@xt                                                                    
\newif\ifl@ftedge                                                               
\newif\ifr@ghtedge                                                              
\newif\ifext@nded                                                               
                                                                                
\newdimen\lm@ff                                                                 
\newdimen\rm@ff                                                                 
\newdimen\lb@ff                                                                 
\newdimen\rb@ff                                                                 
\newdimen\lt@p                                                                  
\newdimen\rt@p                                                                  
                                                                                
\newcount\l@sttreebox     % These four counter allocations have been copied
\newcount\l@sttreeheight  % to this position from the \Tree command
\newcount\l@stdiminfo     % (vs. 2.1). Previously each tree allocated its own
\newcount\l@sttreetype    % counters, using up counters for nothing.

\s@ries{xv@l}{0//1//1//1//1//2//1//3//2//3//4//5//1//6//%                       
              5//4//3//5//2//5//3//4//5//6}                                     
\s@ries{yv@l}{1//6//5//4//3//5//2//5//3//4//5//6//1//5//%                       
              4//3//2//3//1//2//1//1//1//1}

\def\hv@ldef{%                                                                  
     \for\t@mpcnta:=1\to24%                                                     
      \do\g\t@mpdima\vd@st\g\multiply\t@mpdima by\xv@l{\t@mpcnta}%              
         \g\divide\t@mpdima by\yv@l{\t@mpcnta}\g\multiply\t@mpdima by 2         
         \expandafter\gxdef\csname hv@l\romannumeral\t@mpcnta\endcsname{%       
                          \the\t@mpdima}%                                       
      \od}                                                                      
                                                                                
\def\hv@l#1{\csname hv@l\romannumeral#1\endcsname}

\def\p@s#1#2#3{%                                                                
     \g#1\csname#3info\endcsname                                                
     \gxdef\t@mp{\csname#3ht\endcsname}%                                        
     \ifnum\t@mp<0 \gxdef\t@mp{0}\fi                                            
     #2{#1}%                                                                    
     }                                                                          
                                                                                
\chardef\@lmoff0 \chardef\@rmoff1 \chardef\@ltop4 \chardef\@rtop5               
\chardef\@lboff2 \chardef\@rboff3 \chardef\@loff4 \chardef\@roff5               
                                                                                
\def\lmoff#1{\g\advance#1 by\@lmoff}                                            
\def\rmoff#1{\g\advance#1 by\@rmoff}                                            
\def\lboff#1{\g\advance#1 by\@lboff}                                            
\def\rboff#1{\g\advance#1 by\@rboff}                                            
\def\ltop#1{\g\advance#1 by\@ltop}                                              
\def\rtop#1{\g\advance#1 by\@rtop}                                              
\def\loff#1{\g\advance#1 by\@loff\g\advance#1 by\t@mp                           
     \g\advance#1 by\t@mp\relax}                                                
\def\roff#1{\g\advance#1 by\@roff\g\advance#1 by\t@mp                           
     \g\advance#1 by\t@mp\relax}                                                
                                                                                
\def\n@meinfo#1{%                                                               
     \n@me@nfo{#1}{lmoff}\n@me@nfo{#1}{rmoff}%                                  
     \n@me@nfo{#1}{lboff}\n@me@nfo{#1}{rboff}%                                  
     \n@me@nfo{#1}{ltop}\n@me@nfo{#1}{rtop}%                                    
     \n@me@nfo{#1}{loff}\n@me@nfo{#1}{roff}%                                    
     }                                                                          
                                                                                
\def\n@me@nfo#1#2{%                                                             
     \p@s\t@mpcnta{\csname#2\endcsname}{#1}%                                    
     \expandafter\gxdef\csname#1#2\endcsname{\dimen\the\t@mpcnta}}              
                                                                                
\def\n@metree#1#2#3#4#5{%                                                       
     \expandafter\gxdef\csname#5ht\endcsname{\count\the#1}%                     
     \expandafter\gxdef\csname#5info\endcsname{\count\the#2}%                   
     \expandafter\gxdef\csname#5box\endcsname{\the#3}%                          
     \expandafter\gxdef\csname#5type\endcsname{\toks\the#4}%                    
     \n@meinfo{#5}%                                                             
     }                                                                          
                                                                                
\chardef\@cntoff3 \chardef\@boxoff1 \chardef\@dimoff2 \chardef\@typeoff1        
                                                                                
\def\pr@vioustree{%                                                             
     \g\advance\l@sttreeheight by-\@cntoff                                      
     \g\advance\l@stdiminfo by-\@cntoff                                         
     \g\advance\l@sttreetype by-\@cntoff                                        
     \g\advance\l@sttreebox by-\@boxoff                                         
     \n@mel@st                                                                  
     }                                                                          
                                                                                
\def\@ddname#1#2{%                                                              
     \expandafter\gxdef\csname#2ht\endcsname{\csname#1ht\endcsname}%            
     \expandafter\gxdef\csname#2info\endcsname{\csname#1info\endcsname}%        
     \expandafter\gxdef\csname#2type\endcsname{\csname#1type\endcsname}%        
     \expandafter\gxdef\csname#2box\endcsname{\csname#1box\endcsname}%          
     \n@meinfo{#2}%                                                             
     }                                                                          
                                                                                
\def\n@xttree{%                                                                 
     \p@s\t@mpcnta\loff{l@st}\g\advance\t@mpcnta by\@dimoff                     
     \g\advance\l@sttreeheight by\@cntoff                                       
     \g\advance\l@stdiminfo by\@cntoff                                          
     \g\advance\l@sttreetype by\@cntoff                                         
     \g\advance\l@sttreebox by\@boxoff                                          
     \g\count\l@stdiminfo\t@mpcnta                                              
     }                                                                          
                                                                                
\def\@ppenddummy{% pushs a dummy onto the stack and names it `l@st'             
     \n@xttree \g\count\l@sttreeheight-\@ne\n@mel@st                            
     \l@sttype{circle}%                                                         
     \g\setbox\l@stbox\copy\voidb@x                                             
     \g\l@stlmoff=0pt\g\l@strmoff=0pt\g\l@stlboff=0pt\g\l@strboff=0pt%          
     \g\l@stltop=0pt\g\l@strtop=0pt\g\l@stloff=0pt\g\l@stroff=0pt%              
     }                                                                          
                                                                                
\def\g@tchildren{% enables us to talk about the left and the right child
     \ifl@fttop\@ddname{l@st}{l@ft}%                                            
               \pr@vioustree                                                    
               \@ddname{l@st}{r@ght}%                                           
          \else\@ddname{l@st}{r@ght}%                                           
               \pr@vioustree                                                    
               \@ddname{l@st}{l@ft}%                                            
            \fi                                                                 
     \ifnum\r@ghtht>\l@ftht\relax                                               
               \r@ghthighertrue                                                 
               \@ddname{r@ght}{m@x}%                                            
               \@ddname{l@ft}{m@n}%                                             
          \else\r@ghthigherfalse                                                
               \@ddname{l@ft}{m@x}%                                             
               \@ddname{r@ght}{m@n}%                                            
            \fi                                                                 
               }                                                                
                                                                                
\def\n@mel@st{%                                                                 
     \n@metree\l@sttreeheight\l@stdiminfo\l@sttreebox\l@sttreetype{l@st}}       
                                                                                
\def\beginTree{%                                                                
     \begingroup                                                                
     \unitlength 1pt%                                                           
     \divide\unitlength by 65536
     \l@sttreebox\count14                                                       
     \l@sttreeheight\count10                                                    
     \advance\l@sttreeheight by \@ne                                            
     \count\l@sttreeheight=-1                                                   
     \l@stdiminfo\l@sttreeheight                                                
     \advance\l@stdiminfo by \@ne                                               
     \count\l@stdiminfo\count11                                                 
     \advance\count\l@stdiminfo by -5                                           
     \l@sttreetype\l@stdiminfo                                                  
     \advance\l@sttreetype by\@ne                                               
     \count\l@sttreetype\count15                                                
     \n@mel@st\ignorespaces                                                     
     }

\def\endTree{%                                                                  
     \g\leftdist-\l@stlmoff\g\advance\leftdist by \l@stltop                     
     \g\rightdist\l@strmoff\g\advance\rightdist by\l@strtop                     
     \g\setbox\TeXTree\box\l@stbox\endgroup\ignorespaces}

\def\th@ck{\let\@linefnt\tenlnw                                                 
     \@wholewidth\fontdimen8\tenlnw\@halfwidth.5\@wholewidth}                   
                                                                                
\def\leftthick{\g\let\l@ftthick\th@ck}                                          
\def\rightthick{\g\let\r@ghtthick\th@ck}                                        
\def\lft#1{\g\setbox\leftb@x\hbox{#1\ }}                                        
\def\rght#1{\g\setbox\rightb@x\hbox{\ #1}}                                      
\def\cntr#1{\g\setbox\centerb@x\hbox{#1\strut}}                   
\def\bnth#1{\g\setbox\beneathb@x\hbox to0pt{\hss\strut#1\hss}}                        
\def\type#1{%                                                                   
     \g\setbox\centerb@@x\copy\csname#1b@x\endcsname                            
     \g\centerdim@n\csname#1w@dth\endcsname                                     
     \typ@{#1}%                                                                 
     \g\halfcenterdim@n=.5\centerdim@n}                                         
                                                                                
\def\ext@nded{\g\ext@ndedfalse} % This definition must precede                  

\def\node#1{%                                                                   
     \g\setbox\leftb@x\copy\voidb@x                                             
     \g\setbox\rightb@x\copy\voidb@x                                            
     \g\setbox\centerb@x\copy\voidb@x                                           
     \g\setbox\beneathb@x\copy\voidb@x                                          
     \type{circle}%                                                             
     \g\l@fttopfalse\g\l@ftonlyfalse\g\l@ftedgetrue                             
     \g\r@ghtonlyfalse\g\r@ghtedgetrue\g\@xtfalse\ext@nded\n@dummy              
     \g\let\l@ftthick\relax\g\let\r@ghtthick\relax                              
     #1% 
     \@pdcenter                                                                       
     \d@mmy
     \n@de
     \ignorespaces                                                              
     }

\def\@pdcenter{\csname\the\typ@ @cntr\endcsname}

\let\circle@cntr\relax
\let\square@cntr\relax
\let\triangle@cntr\relax
\let\dot@cntr\relax

\def\text@cntr{%
     \g\centerdim@n\wd\centerb@x
     \g\halfcenterdim@n.5\centerdim@n}                                                                          
     
\def\frame@cntr{% 
     \g\setbox\centerb@x\hbox{\ \unhcopy\centerb@x\ }
     \g\centerdim@n\wd\centerb@x
     \g\halfcenterdim@n.5\centerdim@n
     \g\setbox\centerb@@x\rect@ngle{\centerdim@n}{\squarew@dth}{.4pt}}
                                                                           
\def\leftonly{\g\l@ftonlytrue\g\r@ghtedgefalse\g\let\d@mmy\l@ftdummy}           
\def\rightonly{\g\r@ghtonlytrue\g\l@ftedgefalse\g\let\d@mmy\r@ghtdummy}         
\def\unary{\g\r@ghtedgefalse\g\let\d@mmy\@ndummy}                               
\def\external{\g\@xttrue\g\l@ftedgefalse\g\r@ghtedgefalse\g\let\d@mmy\@xtdummy} 
                                                                                
\def\lefttop{\g\l@fttoptrue}                                                    
                                                                                
\def\@xtdummy{%                                                                 
     \@ppenddummy                                                     
     \g\l@strtop-\halfmins@p                                         
     \@ppenddummy
     \g\l@stltop-\halfmins@p                                         
     }                                                                          
                                                                                
\def\n@dummy{\g\let\d@mmy\relax}                                                
                                                                                
\def\l@ftdummy{% cf. \g@tposition                                                               
     \@ppenddummy                                                      
     \g\l@stltop=\dummyhalfcenterdim@n                                          
     \g\l@strtop=\dummyhalfcenterdim@n                                          
     }                                                                          
                                                                                
\def\r@ghtdummy{% cf. \g@tposition                                                              
     \lefttop                                                                   
     \@ppenddummy                                                      
     \g\l@stltop=\dummyhalfcenterdim@n                                          
     \g\l@strtop=\dummyhalfcenterdim@n                                          
     }                                                                          
                                                                                
\def\@ndummy{%                                                                  
     \g\t@mpdima\l@strtop\relax                                                 
     \@ppenddummy                                                     
     \g\l@stltop-\mins@p\g\advance\l@stltop by-\t@mpdima                        
     \g\l@strtop=\t@mpdima                                                      
     }                                                                          
                                                                                
\def\n@de{%                                                                     
     \g@tposition       % naming children and calculating \sl@pe and \tots@p    
     \g@tlt@p\g@trt@p   % calculating \lt@p and \rt@p                           
     \g@tlm@ff\g@trm@ff % calculating \lm@ff and \rm@ff                         
     \g@tlb@ff\g@trb@ff % calculating \lb@ff and \rb@ff                         
     \@pdlroff          % updating loff and roff for all levels but the top one 
     \@pdloffl\@pdroffl % updating loff(1) and roff(1) of the parent tree             
     \@pddim            % updating ltop, rtop, lmoff, rmoff, lboff, and rboff   
     \@pdinfo\@pdht     % updating diminfo and treeheight                       
     \@pdbox            % updating treebox                                      
     \@pdtype           % updating type                                         
     \n@mel@st          % giving the name `l@st' to the new tree 
     \ignorespaces                                                              
     }                                                                          
                                                                                
\def\g@tposition{% naming children and calculating \sl@pe, \tots@p, and node offsets
     \g@tchildren\c@lcsep\c@lcslope\c@lcoffsets
     \ifext@nded\relax                                                          
           \else\ifl@ftonly\g\r@ghtrtop=-\tots@p                                
                           \g\advance\r@ghtrtop by\l@ftrtop                     
                        \fi                                                     
                \ifr@ghtonly\g\l@ftltop=-\tots@p                                
                            \g\advance\l@ftltop by\r@ghtltop                    
                         \fi                                                    
             \fi % cf. \l@ftdummy and \r@ghtdummy                                                               
     }                                                                          
                                                                                
\def\@pdinfo{% updating diminfo                                                 
     \g\l@stinfo=\m@xinfo\relax                                                 
     }                                                                          
                                                                                
\def\@pdht{% updating treeheight                                                
     \g\l@stht=\m@xht                                                           
     \g\advance\l@stht by\@ne                                                   
     }                                                                          
                                                                                
\def\@pdtype{% updating type                                                    
     \g\l@sttype\typ@                                                           
     }                                                                          
                                                                                
\def\g@tlt@p{% calculating \lt@p                                                
     \g\lt@p\wd\leftb@x\g\advance\lt@p by\halfcenterdim@n                       
     }                                                                          
                                                                                
\def\g@trt@p{% calculating \rt@p                                                
     \g\rt@p\wd\rightb@x\g\advance\rt@p by\halfcenterdim@n                      
     }                                                                          
                                                                                
\def\g@tlm@ff{% calculating \lm@ff                                              
     \g\lm@ff\l@ftlmoff                                                         
     \g\advance\lm@ff by-\l@ftltop                                              
     \g\advance\lm@ff by-\halftots@p                                            
     \g\advance\lm@ff by\lt@p\relax                                             
     \ifnum\l@ftht<\r@ghtht\relax                                               
           \g\t@mpdima\r@ghtlmoff                                               
           \g\advance\t@mpdima by-\r@ghtltop                                    
           \g\advance\t@mpdima by\halftots@p                                    
           \g\advance\t@mpdima by\lt@p\relax                                    
           \ifdim\t@mpdima<\lm@ff\relax                                         
                 \g\lm@ff\t@mpdima                                              
              \fi                                                               
        \fi                                                                     
     \ifdim0pt<\lm@ff\relax                                                     
           \g\lm@ff=0pt%                                                        
        \fi                                                                     
     }                                                                          
                                                                                
\def\g@trm@ff{% calculating \rm@ff                                              
     \g\rm@ff\r@ghtrmoff                                                        
     \g\advance\rm@ff by\r@ghtrtop                                              
     \g\advance\rm@ff by\halftots@p                                             
     \g\advance\rm@ff by-\rt@p\relax                                            
     \ifnum\r@ghtht<\l@ftht\relax                                               
           \g\t@mpdima\l@ftrmoff                                                
           \g\advance\t@mpdima by\l@ftrtop                                      
           \g\advance\t@mpdima by-\halftots@p                                   
           \g\advance\t@mpdima by-\rt@p\relax                                   
           \ifdim\t@mpdima>\rm@ff\relax                                         
                 \g\rm@ff\t@mpdima                                              
              \fi                                                               
        \fi                                                                     
     \ifdim0pt>\rm@ff\relax                                                     
           \g\rm@ff=0pt                                                         
        \fi                                                                     
     }                                                                          
                                                                                
\def\g@tlb@ff{% calculating \lb@ff                                              
     \if@xt\g\lb@ff0pt%
      \else\ifnum\l@ftht<\r@ghtht\relax                                               
                 \g\lb@ff\r@ghtlboff                                                  
                 \g\advance\lb@ff by-\r@ghtltop                                       
                 \g\advance\lb@ff by\halftots@p                                       
                 \g\advance\lb@ff by\lt@p\relax                                       
            \else\g\lb@ff\l@ftlboff                                                   
                 \g\advance\lb@ff by-\l@ftltop                                        
                 \g\advance\lb@ff by-\halftots@p                                      
                 \g\advance\lb@ff by\lt@p\relax                                       
              \fi
        \fi                                                                     
     }                                                                          
                                                                                
\def\g@trb@ff{% calculating \rb@ff                                              
     \if@xt\g\rb@ff0pt%
      \else\ifnum\r@ghtht<\l@ftht\relax                                               
                 \g\rb@ff\l@ftrboff                                                   
                 \g\advance\rb@ff by\l@ftrtop                                         
                 \advance\rb@ff by-\halftots@p                                        
                 \g\advance\rb@ff by-\rt@p\relax                                      
            \else\g\rb@ff\r@ghtrboff                                                  
                 \g\advance\rb@ff by\r@ghtrtop                                        
                 \g\advance\rb@ff by\halftots@p                                       
                 \g\advance\rb@ff by-\rt@p\relax                                      
              \fi
        \fi                                                                     
     }                                                                          
                                                                                
\def\@pdlroff{% updating loff and roff for all levels but the top one           
     \ifr@ghthigher\g\t@mpdima-\r@ghtltop\relax                                 
                   \g\t@mpdimb\l@ftlboff                                        
                   \g\advance\t@mpdimb by-\l@ftltop\relax                       
              \else\g\t@mpdima\l@ftrtop\relax                                   
                   \g\t@mpdimb\r@ghtlboff                                       
                   \g\advance\t@mpdimb by\r@ghtrtop\relax                       
                \fi                                                             
     \ifr@ghthigher\p@s\t@mpcnta\loff{m@n}% pointer to loff(1) of smaller tree  
                   \p@s\t@mpcntb\loff{m@x}% pointer to loff(1) of larger tree  
              \else\p@s\t@mpcnta\roff{m@n}% pointer to roff(1) of smaller tree  
                   \p@s\t@mpcntb\roff{m@x}% pointer to roff(1) of larger tree  
                \fi  % if the right tree is the higher one you have to shift    
     \ex\m@nht\times                                                            
        \g\advance\t@mpdima by\dimen\t@mpcntb                                   
        \g\dimen\t@mpcntb\dimen\t@mpcnta                                        
        \g\advance\t@mpcnta by-\@dimoff                                         
        \g\advance\t@mpcntb by-\@dimoff\relax                                   
     \xe                                                                        
     \ifnum\m@xht=\m@nht\relax                                                  
      \else\g\advance\dimen\t@mpcntb by\t@mpdima                                
           \ifnum\l@ftht<\r@ghtht\relax                                         
                 \g\advance\dimen\t@mpcntb by\tots@p                            
            \else\g\advance\dimen\t@mpcntb by-\tots@p                           
              \fi                                                               
           \g\advance\dimen\t@mpcntb by-\t@mpdimb                               
        \fi                                                                     
     }                                                                          
                                                                                
\def\@pdloffl{% updating loff(1) of parent tree                                 
     \p@s\t@mpcnta\loff{m@x}%                                                   
     \g\advance\t@mpcnta by \@dimoff\relax % pointer to loff(0) of parent tree  
     \g\dimen\t@mpcnta\lt@p                                                     
     \g\advance\dimen\t@mpcnta by-\halftots@p                                   
     \g\advance\dimen\t@mpcnta by-\l@ftltop\relax                               
     }                                                                          
                                                                                
\def\@pdroffl{% updating roff(1) of parent tree                                 
     \p@s\t@mpcnta\roff{m@x}%                                                   
     \g\advance\t@mpcnta by \@dimoff\relax % pointer to roff(0) of parent tree  
     \g\dimen\t@mpcnta-\rt@p                                                    
     \g\advance\dimen\t@mpcnta by\halftots@p                                    
     \g\advance\dimen\t@mpcnta by\r@ghtrtop\relax                               
     }                                                                          
                                                                                
\def\@pddim{% updating ltop, rtop, lmoff, rmoff, lboff, and rboff               
     \g\m@xlmoff=\lm@ff\g\m@xrmoff=\rm@ff                                       
     \g\m@xlboff=\lb@ff\g\m@xrboff=\rb@ff                                       
     \g\m@xltop=\lt@p\g\m@xrtop=\rt@p                                           
     }                                                                          
                                                                                
\def\@pdbox{% pushing the nodebox on the stack: updating treebox                
     \g\@xarg\xv@l\sl@pe\g\@yarg\yv@l\sl@pe                                     
     \ifnum\sl@pe=1 % vertical edge                                             
           \g\t@mpdima\vd@st                                                    
           \g\advance\t@mpdima by-\y@ff\typ@                                    
           \g\advance\t@mpdima by-\y@ff\l@fttype                                
           \g\@larg\t@mpdima % \@larg is a number register!                     
           \g\t@mpdima\vd@st                                                    
           \g\advance\t@mpdima by-\y@ff\typ@                                    
           \g\advance\t@mpdima by-\y@ff\r@ghttype                               
           \g\@rarg\t@mpdima % \@rarg is a number register!                     
      \else\g\t@mpdima\halftots@p                                               
           \g\advance\t@mpdima by-\x@ff\typ@                                    
           \g\advance\t@mpdima by-\x@ff\l@fttype                                
           \g\@larg\t@mpdima % \@larg is a number register!                     
           \g\t@mpdima\halftots@p                                               
           \g\advance\t@mpdima by-\x@ff\typ@                                    
           \g\advance\t@mpdima by-\x@ff\r@ghttype                               
           \g\@rarg\t@mpdima % \@rarg is a number register!                     
        \fi                                                                     
     \g\setbox\l@sttreebox\hbox{%
           \ifvoid\leftb@x\relax
             \else\hskip-\halfcenterdim@n\hskip-\wd\leftb@x
                  \unhcopy\leftb@x\hskip\halfcenterdim@n
               \fi
           \ifvoid\centerb@x\relax
             \else\g\t@mpdima-.5\wd\centerb@x\hskip\t@mpdima
                  \unhbox\centerb@x\hskip\t@mpdima
               \fi
           \ifvoid\rightb@x\relax
             \else\g\t@mpdima-\wd\rightb@x\hskip\halfcenterdim@n
                  \unhbox\rightb@x\hskip\t@mpdima\hskip-\halfcenterdim@n
               \fi
           \raise\based@st\copy\centerb@@x
           \if@xt\relax
                 \lower\s@ze pt\hbox to0pt{\hss\unhbox\beneathb@x\hss}%
            \else\hskip-\halftots@p
                 \lower\vd@st\box\l@ftbox
                 \ifl@ftedge\drawl@ftedge\else\hskip\halftots@p\fi
                 \ifr@ghtedge\drawr@ghtedge\else\hskip\halftots@p\fi
                 \lower\vd@st\box\r@ghtbox
                 \hskip-\halftots@p
              \fi
           }% of hbox
     }                                                                          
                                                                                
\def\drawl@ftedge{%                                                             
           \hskip\x@ff\l@fttype                                                 
           \g\t@mpdimc\y@ff\l@fttype\g\advance\t@mpdimc by\based@st
           \g\advance\t@mpdimc-\vd@st
           \raise\t@mpdimc                                                      
           \hbox{\l@ftthick\lineseg(\@xarg,\@yarg){\@larg}}%
           \hskip\x@ff\typ@                                                     
     }                                                                          
                                                                                
\def\drawr@ghtedge{%                                                            
           \hskip\x@ff\typ@                                                     
           \g\t@mpdimc\vd@st                                                      
           \g\advance\t@mpdimc by \based@st
           \g\advance\t@mpdimc by -\y@ff\typ@\relax
           \g\advance\t@mpdimc by- \vd@st
           \raise\t@mpdimc                                                      
           \hbox{\r@ghtthick\lineseg(\@xarg,-\@yarg){\@rarg}}% 
           \hskip\x@ff\r@ghttype                                                
     }                                                                          
                                                                                
\def\x@ff#1{%                                                                   
     \csname\the#1x@ff\endcsname\sl@pe                                          
     }                                                                          
                                                                                
\def\y@ff#1{%                                                                   
     \csname\the#1y@ff\endcsname\sl@pe                                          
     }                                                                          
                                                                                
\def\c@lcslope{%                                                                
     \g\sl@pe1                                                                  
       \loop                                                                    
      \ifdim\hv@l\sl@pe < \tots@p                                               
            \g\advance\sl@pe by1                                                
     \repeat
     \g\tots@p\hv@l\sl@pe                                                       
     \g\halftots@p.5\tots@p}

\def\c@lcsep{%                                                                  
     \g\tots@p\mins@p                                                           
     \g\advance\tots@p by\l@ftrtop                                              
     \g\advance\tots@p by\r@ghtltop\relax                                       
     \g\currs@p\mins@p                                                          
     \p@s\t@mpcnta\roff{l@ft}%                                                  
     \p@s\t@mpcntb\loff{r@ght}%                                                 
     \g\adds@pfalse                                                             
     \g\t@mpcntc\m@nht                               
     \ex\t@mpcntc\times                                                         
        \g\advance\currs@p by-\dimen\t@mpcnta                                   
        \g\advance\currs@p by \dimen\t@mpcntb                                   
        \ifdim\mins@p<\currs@p                                                  
         \else\g\adds@ptrue                                                     
           \fi                                                                  
        \ifdim\currs@p<\mins@p                                                  
              \g\advance\tots@p by\mins@p                                       
              \g\advance\tots@p by -\currs@p                                    
              \g\currs@p\mins@p                                                 
           \fi                                                                  
        \g\advance\t@mpcnta by -\@dimoff                                        
        \g\advance\t@mpcntb by -\@dimoff                                        
     \xe                                                                        
     \ifadds@p\g\advance\tots@p by\adds@p\fi}                                   
                                                                                
\def\tri@ngle{%                                                                 
     \vtop{\g\@xarg\xv@l\sl@pe \g\@yarg\yv@l\sl@pe                              
           \g\t@mpdimc\l@vels\vd@st 
           \g\advance\t@mpdimc by .5\squarew@dth
           \g\multiply\t@mpdimc\xv@l\sl@pe
           \g\divide\t@mpdimc\yv@l\sl@pe
           \g\@larg\t@mpdimc                                                    
           \offinterlineskip                                                    
           \vskip0pt% Force the reference point to the top                                                         
           \hbox to0pt{\hss\lineseg(\@xarg,\@yarg){\@larg}%                     
                       \hskip\t@mpdimc\rlap{\lineseg(-\@xarg,\@yarg){\@larg}}%  
                       \hss}%                                                   
           \setbox\t@mpboxa                                                          
           \hbox to0pt{\hss\vrule height.2pt depth.2pt width2\t@mpdimc\hss}%
           \t@mpdimc-.5\squarew@dth\advance\t@mpdimc\based@st
           \ht\t@mpboxa0pt\dp\t@mpboxa\t@mpdimc\copy\t@mpboxa    
          }%                                                                    
     }                                                                          
                                                                                
\def\lvls#1{\g\l@vels#1}                                                        
\def\slnt#1{\g\sl@pe#1}                                                         
                                                                                
\def\treesymbol#1{%                                                             
     \g\setbox\leftb@x\copy\voidb@x                                             
     \g\setbox\rightb@x\copy\voidb@x                                            
     \g\setbox\centerb@x\copy\voidb@x                                           
     \g\setbox\beneathb@x\copy\voidb@x                                          
     \lvls{1}\slnt{3}%                                                          
     #1%                                                                        
     \g\centerdim@n\trianglew@dth                                               
     \g\halfcenterdim@n.5\trianglew@dth                                         
     \n@xttree                                                                  
     \g\count\l@sttreeheight\l@vels% \g\advance\count\l@sttreeheight by\tw@  
     \g\toks\l@sttreetype{triangle}%                                            
     \n@mel@st                                                                  
     \g\hd@st\hv@l\sl@pe \g\divide\hd@st by\tw@                                 
     \g\l@stltop=\halfcenterdim@n\g\advance\l@stltop by\wd\leftb@x              
     \g\l@strtop=\halfcenterdim@n\g\advance\l@strtop by\wd\rightb@x             
     \g\l@stlboff=-\hd@st \g\multiply\l@stlboff by\l@vels   
     \g\advance\l@stlboff by\wd\leftb@x
     \g\l@strboff=\hd@st \g\multiply\l@strboff by\l@vels 
     \g\advance\l@strboff by-\wd\rightb@x
     \g\l@stlmoff=\l@stlboff\relax
     \ifdim\l@stlmoff>0pt\relax\g\l@stlmoff=0pt\fi   
     \g\l@strmoff=\l@strboff
     \ifdim\l@strmoff<0pt\relax\g\l@strmoff=0pt\fi 
     \g\t@mpcnta\l@stinfo\g\advance\t@mpcnta by6% preliminary                   
     \ex\l@vels\times                                                           
        \g\dimen\t@mpcnta-\hd@st\g\advance\t@mpcnta by\@ne                      
        \g\dimen\t@mpcnta\hd@st\g\advance\t@mpcnta by\@ne                       
     \xe                                                                        
     \g\advance\t@mpcnta by-\tw@                                                
     \g\advance\dimen\t@mpcnta by\wd\leftb@x                                    
     \g\advance\t@mpcnta by\@ne                                                 
     \g\advance\dimen\t@mpcnta by-\wd\rightb@x                                  
     \g\setbox\l@stbox\vtop % to\l@vels\vd@st
            {\offinterlineskip                   
             \g\setbox\t@mpboxa                                                 
             \hbox{\hskip-\halfcenterdim@n\hskip-\wd\leftb@x\unhbox\leftb@x
                   \hskip\halfcenterdim@n
                   \raise\based@st\tri@ngle
                   \hskip\halfcenterdim@n\t@mpdima-\wd\rightb@x
                   \unhbox\rightb@x\hskip\t@mpdima\hskip-\halfcenterdim@n}                                     
             \g\ht\t@mpboxa=0pt\box\t@mpboxa                                    
             \setbox\centerb@x\hbox to0pt{\hss\unhbox\centerb@x\hss}%
             \ht\centerb@x0pt\dp\centerb@x0pt\box\centerb@x
             \vskip\s@ze pt
             \ht\beneathb@x0pt\box\beneathb@x                                                 
             \vskip-\dp\beneathb@x\vskip-\ht\beneathb@x}%
     \ignorespaces                                                              
     }

\def\norm@ff{% everything set up for 10pt node size                                                                 
\s@ries{circley@ff}{0.50000pt//0.49320pt//0.49029pt//0.48507pt//%               
                    0.47434pt//0.46424pt//0.44721pt//0.42875pt//%               
                    0.41603pt//0.40000pt//0.39043pt//0.38411pt//%               
                    0.35355pt//0.32009pt//0.31235pt//0.30000pt//%               
                    0.27735pt//0.25725pt//0.22361pt//0.18570pt//%               
                    0.15811pt//0.12127pt//0.09806pt//0.08220pt}%                
     }                                                                          
                                                                                
\def\dotx@ff#1{0pt}
\def\doty@ff#1{0pt} 

\def\trianglex@ff#1{0pt}
\def\triangley@ff#1{0pt} 

\def\c@lcoffsets{%
     \ifnum\sl@pe=\@ne\relax
           \xdef\circlex@ff##1{0pt}%
      \else\g\t@mpcnta26 % number of slopes + 2
           \g\advance\t@mpcnta-\sl@pe
           \xdef\circlex@ff##1{\circley@ff\t@mpcnta}%
        \fi
     \ifnum\sl@pe<13\relax % incoming edge meets upper border of a square node
           \g\t@mpdima.5\squarew@dth
           \xdef\squarey@ff##1{\the\t@mpdima}%
           \g\multiply\t@mpdima\xv@l\sl@pe
           \g\divide\t@mpdima\yv@l\sl@pe
           \xdef\squarex@ff##1{\the\t@mpdima}%
      \else\g\t@mpdima.5\squarew@dth
           \xdef\squarex@ff##1{\the\t@mpdima}%
           \g\multiply\t@mpdima\yv@l\sl@pe
           \g\divide\t@mpdima\xv@l\sl@pe
           \xdef\squarey@ff##1{\the\t@mpdima}%
        \fi
     \g\t@mpdima.5\squarew@dth
     \xdef\texty@ff##1{\the\t@mpdima}%
     \g\multiply\t@mpdima\xv@l\sl@pe
     \g\divide\t@mpdima\yv@l\sl@pe
     \xdef\textx@ff##1{\the\t@mpdima}% 
     \let\framex@ff\textx@ff
     \let\framey@ff\texty@ff
    }
                                                                             
\def\upds@ze#1{%                                                                
     \for\t@mpcntc:=1\to24                                                      
      \do\g\t@mpdimc=\csname#1\romannumeral\t@mpcntc\endcsname\relax            
         \g\multiply\t@mpdimc by\s@ze                                           
         \expandafter\gxdef\csname#1\romannumeral\t@mpcntc\endcsname            
                            {\the\t@mpdimc}%                                    
      \od}                                                                      
                                                                                
\def\nodes@ze{%                                                                 
     \begingroup                                                                
     \unitlength 1pt%                                                           
     \divide\unitlength by 65536                                                
     \g\based@st\s@ze pt\g\divide\based@st by 10 % \based@st is 10 % of         
     \g\dummyhalfcenterdim@n=\s@ze pt\g\divide\dummyhalfcenterdim@n by\tw@      
     \g\circlew@dth=\s@ze pt%                                                   
     \g\t@mpcntc\s@ze\g\multiply\t@mpcntc by 65536                              
     \g\setbox\circleb@x\hbox to0pt{\circle{\t@mpcntc}\hss}%                     
     \upds@ze{circley@ff}%                                  
     \g\squarew@dth.9pt\g\multiply\squarew@dth by\s@ze                          
     \g\setbox\squareb@x\rect@ngle{\squarew@dth}{\squarew@dth}{.4pt}%      
     \g\dotw@dth=\s@ze pt\g\divide\dotw@dth by 5                                
     \ifdim\dotw@dth < 1pt\relax                                                
           \g\dotw@dth1pt\relax                                                 
        \fi                                                                     
     \g\t@mpcntc\dotw@dth                                               
     \g\setbox\dotb@x\hbox to 0pt{\circle*{\t@mpcntc}\hss}%
     \g\trianglew@dth=\s@ze pt\g\multiply\trianglew@dth by \tw@                 
     \g\divide\trianglew@dth by 3 
     \g\textw@dth=0pt%
     \g\setbox\textb@x\copy\voidb@x
     \g\framew@dth0pt%
     \g\setbox\frameb@x\copy\voidb@x
     \hv@ldef                                                                   
     \endgroup                                                                  
     }                                                                          
                                                                                
\def\treefonts#1{#1}                                                            
\def\vdist#1{\g\vd@st=#1\relax}                                                 
\def\minsep#1{\g\mins@p=#1\relax\g\halfmins@p=.5\mins@p}                        
\def\addsep#1{\g\adds@p=#1\relax}                                               
\def\extended{\def\ext@nded{\g\ext@ndedtrue}}                                   
\def\noextended{\def\ext@nded{\g\ext@ndedfalse}}                                
\def\nodesize#1{\g\t@mpdima=#1\relax\g\s@ze=\t@mpdima                           
     \g\divide\s@ze by 65536\relax} % conversion from dimension to number       
\def\Treestyle#1{\norm@ff#1\nodes@ze\ignorespaces}                                           
                                                                                
\newcount\b@nno % number of nodes                                               
\newcount\b@nlv % number of complete levels                                     
\newcount\b@ndl % number of nodes on incomplete level                           
                                                                                
\def\ld(#1,#2,#3){% #1, #2, and #3 must be counter registers.                   
     #2=0 #3=1                                                                  
       \loop\ifnum #1>\@ne\relax                                                
            \divide #1 by\tw@ % this is integer division                        
            \advance #2 by\@ne                                                  
            \multiply #3 by\tw@                                                 
     \repeat}                                                                   
                                                                                
\def\b@nary#1{% draws a complete binary tree with #1 internal nodes,            
     \b@nno=#1\relax\advance\b@nno by \@ne                                      
     \ld(\b@nno,\b@nlv,\b@ndl)%                                                 
     \b@ndl=-\b@ndl\advance\b@ndl by #1\advance\b@ndl by\@ne                    
     \b@n}                                                                      
                                                                                
\def\b@n{%                                                                      
     \ifnum\b@nlv>\@ne                                                          
           \advance\b@nlv by-\@ne                                               
           \b@n                                                                 
           \b@n                                                                 
           \advance\b@nlv by\@ne                                                
           \node{}                                                              
      \else\ifnum\b@ndl>\@ne                                                    
                 \advance\b@ndl by-\tw@                                         
                 \node{\le@f\external}\node{\le@f\external}\node{}%             
                 \node{\le@f\external}\node{\le@f\external}\node{}%             
                 \node{}%                                                       
            \else\ifnum\b@ndl=\@ne                                              
                       \advance\b@ndl by-\@ne                                   
                       \node{\le@f\external}\node{\le@f\external}\node{}%       
                       \node{\le@f\external}%                                   
                       \node{}%                                                 
                  \else\node{\le@f\external}\node{\le@f\external}\node{}%       
                    \fi                                                         
              \fi                                                               
        \fi}                                                                    
                                                                                
\def\circleleaves{\def\le@f{\type{circle}}}                                     
\def\squareleaves{\def\le@f{\type{square}}}                                     
                                                                                
\newcount\no@                                                                   
\def\no#1{\no@=#1\relax}                                                        
                                                                                
\def\binary#1{%                                                                 
     \no{1}\circleleaves                                                        
     #1%                                                                        
     \b@nary{\no@}}                                                             
                                                                                
\newcount\f@bht                                                                 
                                                                                
\def\f@b{% draws a Fibonacci tree of depth #1                                   
     \ifnum\f@bht>1                                                             
           \advance\f@bht by-\@ne\f@b\advance\f@bht by\@ne                      
           \advance\f@bht by-\tw@\f@b\advance\f@bht by\tw@                      
           \ifunn@des\node{\unary}                                              
                  \fi                                                           
           \node{\lefttop}                                                      
      \else\ifnum\f@bht=1                                                       
                 \node{\external\le@f}                                          
                 \node{\external\le@f}                                          
                 \node{}                                                        
            \else\node{\external\le@f}                                          
              \fi                                                               
        \fi}                                                                    
                                                                                
\newif\ifunn@des                                                                
                                                                                
\let\unarynodes\unn@destrue                                                     
\def\hght#1{\f@bht=#1\relax}                                                    
                                                                                
\def\fibonacci#1{%                                                              
     \hght{0}\unn@desfalse\circleleaves                                         
     #1%                                                                        
     \f@b}

\Treestyle{%  
\ifLaTeX\treefonts{\normalsize\rm}%                                        
\else\treefonts{
\fontsize{10}{12pt}\rmfamily}                                                    
\fi                                                                   
     \vdist{60pt}%                                                              
     \minsep{20pt}%                                                             
     \addsep{0pt}%                                                              
     \nodesize{20pt}%                                                           
     }

\Treestyle{\vdist{20pt}\minsep{16pt}} 

\def\proclaim #1. #2\par{\medbreak
\noindent{\bf#1.\enspace}{\sl#2}\par\medbreak}
\begin{document}

\title[Zimmermann Cancellation in Free Algebras]{Zimmermann Type Cancellation \\in
the Free Fa\`a di Bruno Algebra}

\begin{abstract}
The $N$-variable Hopf algebra introduced by Brouder, Fabretti, and
Krattenaler (BFK) in the context of non-commutative Lagrange inversion can
be identified with the inverse of the incidence algebra of $N$-colored
interval partitions. The (BFK) antipode and its reflection determine the
(generally distinct) left and right inverses of power series with
non-commuting coefficients and $N$ non-commuting variables. As in the case
of the Fa\`{a} di Bruno Hopf algebra, there is an analogue of the Zimmermann 
cancellation formula. The summands of
the (BFK) antipode can indexed by the \emph{depth first} ordering of vertices on
contracted planar trees, and the same applies to the interval partition antipode. Both
can also be indexed by the \emph{breadth first} ordering of vertices in the non-order
contractible planar trees in which precisely one non-degenerate vertex occurs on each
level.
\end{abstract}

\date{April 17, 2005}

\author{Michael Anshelevich}
\address[Michael Anshelevich]{Department of Mathematics\\
University of California
Riverside, CA 92521}
\email{manshel@math.ucr.edu}
\author{Edward G. Effros}
\address[Edward G. Effros]{Department of Mathematics\\
UCLA, Los Angeles, CA 90095-1555}
\email{ege@math.ucla.edu}
\author{Mihai Popa}
\address[Mihai Popa]{Department of Mathematics\\
UCLA, Los Angeles, CA 90095-1555}
\email{mvpopa@math.ucla.edu}
\thanks{Anshelevich and Effros were partially supported by the National Science
Foundation (DMS-0400860,}

\maketitle

%[Lagrange Hopf algebra]
%\large
\section{Introduction}

Non-commutative power series play an important role in a number of areas,
including combinatorics, free probability, and quantum field theory. A
striking aspect of these applications is that one can effectively manipulate series
in which neither the coefficients nor the variables commute. Such
calculations are often simplified through the use of combinatorial indices
such as trees and graphs. In turn, these somewhat \emph{ad hoc} techniques
can frequently be systematized by using Hopf algebras. This approach was
pioneered by Rota and his colleagues in their studies of combinatorics \cite
{J}. More recently Kreimer \cite{K}, and Kreimer and Connes \cite{C} have
used Hopf-theoretic methods to rationalize various Feynman diagram methods
used in perturbative quantum field theory.

An important example of this theory was described by Haiman and Schmitt 
\cite{H}, who showed that calculating the antipode for the reduced Fa\`{a}
di Bruno Hopf algebra is equivalent to finding an explicit Lagrange
inversion formula for factorial power series with commuting coefficients.
For this purpose they realized the Hopf algebra as the incidence algebra of the
colored partitions of finite colored sets, and they proved a summation formula for the antipode
in which the terms are indexed by colored trees.  
They showed that many of the terms
in this sum cancel, and that it suffices to use the reduced trees,  
in which each non-leaf vertex is non-degenerate, i.e, has more than one offspring 
(see \cite{H}, Theorem 7). In their argument they passed from general trees to reduced trees by 
contracting the appropriate edges. As emphasized in \cite{F} (see the
discussion of Zimmermann's formula in \S12.2), Haiman and Schmitt's procedure may be regarded as an
elementary illustration  of the cancellations that
play such an important role in perturbative quantum field theory.

In a recent paper, Brouder, Fabretti, and Krattenaler \cite{B} described a
Hopf algebra, called the \emph{left Lagrange algebra} $\mathcal{
\ L}=\mathcal{L}^{N}$ below, which is related to Lagrange
inversion for power series with non-commuting coefficients and $N$ non-commuting variables. As
they pointed out, the situation is more delicate, since the ``composition''
of such polynomials is generally not associative. Furthermore, they made the
important observation that the antipode $S_{\mathcal{L}}$ is not involutory,
i.e., $S_{\mathcal{L}}^{2}\neq id$ (correcting the statement \cite{S}, Prop.
4.4). As a result, each left Lagrange algebra has a correspoding \emph{%
inverse} Hopf algebra $\mathcal{H}=\mathcal{H}^{N}$ associated with the
antipode $S_{\mathcal{H}}=S_{\mathcal{L}}^{-1}$ (see \S\S5-6). As we will
see, $\mathcal{H}^{N}$ is just the incidence Hopf algebra of the $N$-colored 
interval (i.e., ordered) partitions.

As one might expect, one has a formula for the antipode $S_{\mathcal{H}}$ in which the terms
are indexed by \emph{planar trees}, but one must in addition
keep track of the \emph{breadth first} ordering of the vertices (see (\ref{omegaantipode})).
In contrast to the commutative situation considered by Haiman and Schmitt, one cannot 
contract the singular edges in these trees since that procedure can disrupt the ordering. 

The primary result in this paper is that despite the new complications, a 
non-commutative analogue of 
Haiman and Schmitt's reduced tree formula is valid provided one instead uses the 
\emph{depth first} ordering on
the vertices (Corollary \ref{redantipodeth}). This is initially proved by induction. We also show
that the reduced tree antipode formulae for the
the three Hopf algebras $\mathcal{H}, \mathcal{L}$ and the related \emph{right Lagrange
Hopr algebra} $\mathcal{R}$ can be derived from each other.
 
In \S\S8-9 we show that the reduced tree formulae can be derived by suitable
``ordered'' cancellations. We begin by using order-preserving contractions and
expansions to cancel out all but the ``order reduced simple layered trees''
$\mathbf{OST}$ in the breadth first formula (Theorem \ref{ordantipodeth}). In \S9 we prove that
$\mathbf{OST}$ is precisely the class of layered trees on which  left-to-right breadth first and
right-to-left depth first orderings \emph{coincide} on the  non-degenerate vertices.  We then show
that contraction of
\emph{all} of the singular edges provides a one-to-one correspondence between the trees in
$\mathbf{OST}$ and the reduced planar trees $\mathbf{RT}$. This leads to a constructive proof of the
reduced tree formula for $S_{\mathcal{H}}$, which suggests that similar arguments might be used
for studying antipodes in the non-commutative versions of the Connes-Kreimer algebras. 
   
In \S 10 we prove that despite the fact that the substitution operation for
non-commutative power series is not multiplicatively associative, one can
still use the antipodes of the Lagrange Hopf algebra $\mathcal{L}$ (which \emph{is} associative) and
its reflection $\mathcal{R}$ to find the left and right substitutional inverses of power
series in which neither the variables nor the constants commute. As we
illustrate, these inverses are generally distinct.

In order to make the material more accessible to functional analysts and
mathematical physicists, we have included a summary of the
relevant constructions from algebraic combinatorics. Those familiar with this 
material might prefer to skip the initial sections. 

\section{Ordered sets and their colorings}

A \emph{partially ordered set} $(P,\leq )$ is a set $P$ together with a
relation $\leq $ such that $x\leq y$ and $y\leq x$ if and only if $x=y,$ and 
$x\leq y\leq z$ implies $x\leq z$. We say $P$ is a \emph{linearly ordered set%
} or simply an \emph{ordered set} if $x\leq y$ or $y\leq x$ for all $x,y\in
X.$ If a partially ordered set $P$ has a minimum element, we denote it by $%
0_{P},$ and similarly we write $1_{P}$ for a maximum element. Given $x,y\in
P $ with $x\leq y,$ we let $[x,y]$ denote the \emph{segment} $\left\{ z\in
P:x\leq z\leq y\right\} $. We denote a finite ordered set $S$ by $
(x_{1},\ldots ,x_{p})$ where $x_{1}<\cdots <x_{p}.$ In particular if $p\in 
\Bbb{N},$ we let $[p]=(1,\ldots ,p)$. Given partially ordered sets $P$ and $%
Q,$ we let $P\times Q$ have the \emph{product} partial ordering $%
(x_{1},x_{2})\leq (y_{1},y_{2})$ if $x_{1}\leq y_{1}$ and $x_{2}\leq y_{2}.$

We fix a finite (unordered) set of ``colors'' $\Gamma$, which in most cases
will be $\langle N \rangle=\{1,2,3,\ldots ,N\}$. A \emph{colored} (or $%
\Gamma $-\emph{colored}) \emph{ordered set} $(S,c)$ is an ordered set $S$
together with a function $c =c_{S}:S\rightarrow \Gamma$ (we place no
restrictions on $c$). We refer to $c(x)$ as the \emph{color} of a point $%
x\in S$, and we say that $c $ is a \emph{coloring} of $S.$ If $
S=(x_{1},\ldots ,x_{p})$ and we let $v(i)=c(x_{i})$, $c$ is determined by
the word $v=v(1)\cdots v(p)$ in the free monoid $\Gamma^{*}$ generated by $%
\Gamma$. We will often identify $c:S\to \Gamma$ with the word $v$. The
identity $1$ of this monoid is the empty word $\emptyset $. We write $\left|
v\right| $ for the length of an element $v$ in $\Gamma^{*}$. Two colored
ordered sets $(S,v)$ and $(T,w) $ are \emph{isomorphic} or \emph{have the
same coloring} if there is an order isomorphism of $\theta:S\to T$ that
preserves the coloring, i.e., $c_{T}(\theta(s))=c_{S}(s)$. Since these are
linearly ordered sets, the mapping $\theta$ will necessarily be unique.

Given a finite colored ordered set $S$, any subset $R\subseteq S$ is itself
linearly ordered in the relative order, and we let $R$ have the restricted
coloring $c|R.$ Given disjoint $\Gamma$-colored sets $S$ and $T,$ with
colorings $v$ and $w,$ we let $S\sqcup T$ denote the union with the left to
right ordering, and the coloring $c_{S\sqcup T}=vw$.

\section{Planar trees and their colorings.}

The \emph{planar forests} (or simply ``forests'') are defined recursively.
For transparency we use terms associated with botanical and genealogical
trees. To construct a plane forest $F$ we first choose an ordered set $L_{0}$
of \emph{vertices} $(x_{1},\ldots ,x_{r})$ called the $\emph{roots}$ or the
zeroth \emph{level} of $F.$ For each root $x_{i}$ in $L_{0}$ we then choose
a possibly empty ordered set of vertices $(x_{i1},\ldots ,x_{ir_{i}})$
called the \emph{offspring} or \emph{children} of $x_{i}.$ The entire collection
$L_{1}$ of these offspring is called the first level, which we linearly order
first by their parents and then among children of a given parent by the given
order. Having chosen the $(n-1)$-st level, we choose an ordered set of vertices
for each vertex in that set. These new vertices consititute the $n$-th level,
and we order them in the same manner. We consider only finite forests. Owing
to the conventions we have adopted, we regard the $(n-1)$-st level as being
``higher'' than the $n$-th level. A \emph{tree} is a forest with only one
root.

We identify a forest with a graph in the plane in the usual manner. The
levels are placed in horizontal rows, parents are joined by edges to their
offspring, and the left to right order reflects the recursively defined
order on the parents, and the given order on the offspring in each family. A
typical planar forest is illustrated below: 

\vspace{-.4in}
\setlength{\unitlength}{1cm} 
 
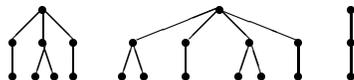
\[
\begin{picture}(8,1.5)
\put(.1,.3){\line(1,1){.4}} 
\put(.5,.7){\line(0,-1){.5}}
\put(.5,.25){\circle*{.1}}
\put(.5,.25){\line(1,-3){.15}}
\put(.4,-.2){\line(1,3){.15}}
\put(.5,.7){\line(1,-1){.4}}

\put(1.55,-.2){\line(1,3){.15}}
\put(1.7,.25){\line(1,-3){.15}}
\put(2.4,.2){\line(0,-1){.4}}
\put(1.55,-.2){\circle*{.1} }
\put(1.85,-.2){\circle*{.1} }
\put(3.9,.2){\line(0,-1){.4}}
\put(2.4,-.2){\circle*{.1} }
\put(3.25,.3){\line(1,-3){.17}}
\put(3.1,-.2){\line(1,3){.15}}
\put(3.1,-.2){\circle*{.1} }
\put(3.4,-.2){\circle*{.1} }
\put(3.9,-.2){\circle*{.1} }

\put(1.75,.28){\line(5,2){1.1}} 
\put(2.45,.3){\line(1,1){.4}}
\put(2.85,.7){\line(1,-1){.4}}
\put(2.9,.67){\line(5,-2){1}}

\put(.1,-.2){\circle*{.1} }
\put(.1,-.2){\line(0,1){.45} }
\put(.4,-.2){\circle*{.1} }
\put(.9,-.2){\circle*{.1} }
\put(.9,-.2){\line(0,1){.45}}
\put(.65,-.2){\circle*{.1} }
\put(.1,.25){\circle*{.1} }
\put(.9,.25){\circle*{.1} }
\put(1.7,.25){\circle*{.1} }
\put(2.4,.25){\circle*{.1} }
\put(3.25,.25){\circle*{.1} }
\put(3.9,.25){\circle*{.1} }

\put(2.85,.68){\circle*{.1} }
\put(.5,.68){\circle*{.1} }

\put(4.6,-.2){\line(0,1){.45} }
\put(4.6,.68){\circle*{.1} }
\put(4.6,-.2){\circle*{.1} }
\put(4.6,.25){\circle*{.1} }
\put(4.6,.3){\line(0,1){.4}}

\end{picture}
\]

\hspace{.2in}

We use the notations $\mathbf{V}(F)$ and $\mathbf{E}(F)$ for the vertices and edges of a forest $F$. A
vertex $x\in \mathbf{V}(F)$ is a 
\emph{leaf} if it has no offspring, it is \emph{unary} if it has precisely
one offspring, and it is \emph{non-degenerate} if it has more than one
offspring. The offspring of a vertex are said to be
\emph{siblings}. An edge is \emph{singular} if it descends from a
unary vertex.

A tree $T$ is said to be
\emph{reduced} if it has no unary vertices, or equivalently,
all of its non-leaf vertices are non-degenerate. Given an arbitrary tree $T$, we 
let $\rho(T)$ be the corresponding reduced tree, in which each singular 
edge from a unary vertex $x$ is contracted to the vertex $x$.

A \emph{layer} (respectively, \emph{branch}) is a forest (respectively,
tree) in which all vertices are roots or leaves. We define the $n$-th \emph{%
layer of a forest} $F$ to be the forest obtained by considering only the
vertices in the $(n-1)$-st and $n$-th levels together with the edges joining
them in $F$. We say that a forest $F$ is \emph{uniformly layered} or simply \emph{layered}
if leaves occur only at the same lowest level. Equivalently, any vertex at a higher level
must have at least one offspring. 

We say that a layered tree is \emph{proper} if each non-leaf level has at least
one non-degenerate vertex (or equivalently, the $n$-th level is larger than the 
$(n-1)$-st level for each $n$). Unless otherwise indicated, 
all layered forests in this paper are assumed to be
proper.

A non-degenerate vertex is \emph{simple} if it is the only non-degenerate
vertex on that level. If that is the case, we say that the corresponding
level is simple. A forest is \emph{simple} if all of its non-degenerate
vertices are simple.

\vspace{0.2cm} The ascending \emph{breadth first} ordering $\ll $
on the vertices of a \emph{layered} forest $F$ is defined as follows. We
write $x\ll y$ if $x$ is on a lower level than $y$, or $x$ is in the same
level as $x$ and lies to the left of $y$ as illustrated in the forest below,
in which $w_{1}\ll w_{2}\ll \dots $. In the literature the term breadth first
usually refers to the corresponding reverse or descending ordering $\gg$. 

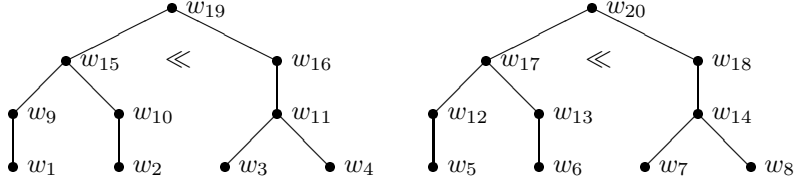
\[
\begin{Tree}
\node{\external\type{dot}{\rght{$w_{1}$}}}
\node{\unary\type{dot}\rght{$w_{9}$}}
\node{\external\type{dot}{\rght{$w_{2}$}}}
\node{\unary\type{dot}\rght{$w_{10}$}}
\node{\type{dot}\rght{$w_{15}$\,\,\,\,\,\,\,\,\,\,$\ll$}}
\node{\external\type{dot}\rght{$w_{3}$}}
\node{\external\type{dot}\rght{$w_{4}$}}
\node{\type{dot}\rght{$w_{11}$}}
\node{\unary\type{dot}\rght{$w_{16}$}}
\node{\type{dot}\rght{$w_{19}$}}
\end{Tree}
\hskip\leftdist\box\TeXTree\hskip\rightdist\qquad 
\begin{Tree}
\node{\external\type{dot}{\rght{$w_{5}$}}}
\node{\unary\type{dot}\rght{$w_{12}$}}
\node{\external\type{dot}{\rght{$w_{6}$}}}
\node{\unary\type{dot}\rght{$w_{13}$}}
\node{\type{dot}\rght{$w_{17}$\,\,\,\,\,\,\,\,\,\,$\ll$}}
\node{\external\type{dot}\rght{$w_{7}$}}
\node{\external\type{dot}\rght{$w_{8}$}}
\node{\type{dot}\rght{$w_{14}$}}
\node{\unary\type{dot}\rght{$w_{18}$}}
\node{\type{dot}\rght{$w_{20}$}}
\end{Tree}
\hskip\leftdist\box\TeXTree\hskip\rightdist\qquad 
\]

We will also use the \emph{left} and \emph{right depth first} total orderings of an arbitrary
tree $T$ and their inverse orderings. Given distinct
vertices $x$ and $y$ in $T$, we write $x\lup y$ (respectively, $x
\rup y$) if either $y$ is an ancestor of $x$, or letting $z$ be the
first common ancestor of $x$ and $y$, the branch headed towards $x$ lies to
the left (respectively, right) of the branch headed towards $y$. We define 
the descending orderings $\rdown$ and $\ldown $ to be the inverses of $\rup$ and $\lup$,
respectively. We have, for example, that $x_{1}\lup x_{2} \lup x_{3} \lup \ldots
$ and
$y_{1}\rup y_{2}\rup y_{3}\rup \ldots $ in the following trees, 
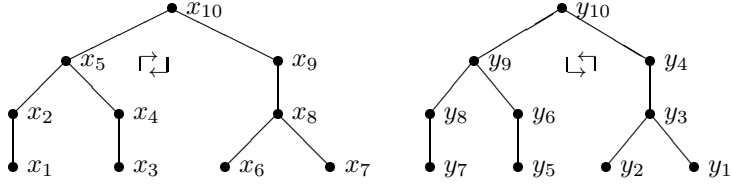
\[
\begin{Tree}
\node{\external\type{dot}{\rght{$x_{1}$}}}
\node{\unary\type{dot}\rght{$x_{2}$}}
\node{\external\type{dot}{\rght{$x_{3}$}}}
\node{\unary\type{dot}\rght{$x_{4}$}}
\node{\type{dot}\rght{$x_{5}$}}
\node{\external\type{dot}\rght{$x_{6}$}}
\node{\external\type{dot}\rght{$x_{7}$}}
\node{\type{dot}\rght{$x_{8}$}}
\node{\unary\type{dot}\rght{$x_{9}$}}
\node{\type{dot}\rght{$x_{10}$}}
\end{Tree}
\hskip\leftdist\raisebox{-.3in}{$\lup \!\!\!\!
\ldown$}\box\TeXTree\hskip\rightdist\qquad  
\begin{Tree}
\node{\external\type{dot}{\rght{$y_{7}$}}}
\node{\unary\type{dot}\rght{$y_{8}$}}
\node{\external\type{dot}{\rght{$y_{5}$}}}
\node{\unary\type{dot}\rght{$y_{6}$}}
\node{\type{dot}\rght{$y_{9}$}}
\node{\external\type{dot}\rght{$y_{2}$}}
\node{\external\type{dot}\rght{$y_{1}$}}
\node{\type{dot}\rght{$y_{3}$}}
\node{\unary\type{dot}\rght{$y_{4}$}}
\node{\type{dot}\rght{$y_{10}$}}
\end{Tree}
\hskip\leftdist\box\TeXTree\raisebox{-.3in}{$\rdown \!\!\!
\rup$}\hskip\rightdist\qquad
\]
whereas $x_{10}\ldown x_{9}\ldown x_{8}\ldown\ldots $ and $y_{10}\rdown y_{9}
\rdown y_{8}\rdown \ldots$ in those same diagrams.
The linear orderings $\ldown$
and $\rdown$
 are known as the (descending) right
and left
\emph{depth first orderings}. To numerically label the vertices
according to $x\ldown y$ one begins at the root, and
then successively chooses right branches going to successive generations,
backtracking to enumerate the next right-most uncounted vertex when necessary.

It is important to remember that if $x\lup y$ or $x\rup y$, one cannot
conclude that the level of $x$ is lower than or equal to that of $y$.
Nevertheless the root is the greatest element in both the $\lup$ and
$\rup$ orderings. If two vertices $x$ and $y$ are on the same level,
then $x\lup y$ if and only if $x\ll y$. Finally $x\lup y$ and $ x\rup
y$ will both hold if and only if $x$ is a descendant of $y$. We will write
$x\Rsh
\!\!\!\!\Rsh y$ if $x\lup y$  but it is not true that $%
x\rup y$, i.e., $x$ is not a descendant of $y$.

A \emph{colored forest} $(F,c)$ (or simply $F$) consists of a forest $F$
together with a coloring of the vertices $c:F\rightarrow \Gamma $ such that
if $x$ is a unary vertex with offspring $y$, then $c(y)=c(x).$ The following
is a 2-colored layered tree. \Treestyle{\nodesize{13pt}\vdist{25pt} %
\minsep{24pt}} 
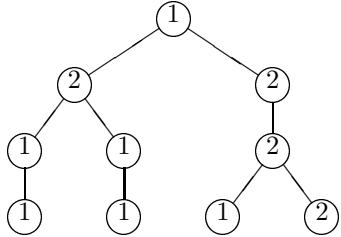
\[
\begin{Tree}
\node{\external\cntr{1}}
\node{\unary\cntr{1}}
\node{\external\cntr{1}}
\node{\unary\cntr{1}}
\node{\cntr{2}}
\node{\external\cntr{1}}
\node{\external\cntr{2}}
\node{\cntr{2}}
\node{\unary\cntr{2}}
\node{\cntr{1}}
\end{Tree}
\hskip\leftdist\box\TeXTree\hskip\rightdist\qquad 
\]

We say that colored layered forests $F$ and $G$ are \emph{compatible} if the
set of roots of $F$ and the leaves of $G$ are colored and order isomorphic.
If that is the case, we define the \emph{right join} $F\vartriangleright G$ be the
forest that results if one identifies these two sets ($G$ lies ``over'' $F$
in the resulting forest). It is apparent that every layered forest $F$ has a
unique decomposition $F=L_{1} \vartriangleright L_{2}\vartriangleright
\ldots \vartriangleright L_{r},$ where $L_{k}$ is the $(r-k+1)$ level of $F$.

If the vertices of compatible layered forests $F$ and $G$ are given by $x_{1}\ll
\ldots \ll x_{m}$ and $y_{1}\ll \ldots \ll y_{n,}$ then the vertices of $
F\vartriangleright G$ are ordered by 
\begin{equation}
x_{1}\ll \ldots \ll x_{m}\ll y_{1}\ll \ldots \ll y_{n}.  \label{breadthchain}
\end{equation}

Given $u,v\in \langle N\rangle^{*}$, we let 
$\mathbf{F}^{v}_{u}, \mathbf{RF}^{v}_{u}, \mathbf{SF}^{v}_{u}$, and 
$\mathbf{LF}^{v}_{u}$ denote the equivalence classes of
general, reduced, simple and (proper) layered forests, respectively, with roots colored by $v$ and
leaves colored by $u$. If $v=i\in\langle N \rangle$, we let
$\mathbf{T}^{i}_{u},\mathbf{RT}^{i}_{u}, \mathbf{ST}
^{i}_{u}$, and $\mathbf{LT}^{i}_{u}$ denote the corresponding classes of trees. 

Given a colored tree $T$, we define $T^{*}$ to be the ``reflected'' tree in
which we reverse the order of the vertices at each level, keeping the same
parental relation. Each ``reflected'' vertex $x^{*}$
is given the color of $x$. It is evident that if $x$ and $y$ are vertices in $T$,
then $x\lup y$ in $T$ if and only if $y^{*}\,\rup x^{*}$ in $T^{*}$, and that
$x\ldown y$ if and only if $y^{*}\rdown x^{*}$.

\section{Colored ordered partitions and their segments}

An \emph{ordered} (or \emph{interval}) \emph{partition} $\sigma
=(B_{1},\ldots ,B_{q})$ of an ordered set $S$ is a collection of subsets for
which $\bigcup B_{k}=S$ and $B_{1}<\ldots <B_{q}$ in the given ``left to
right'' ordering. We may use the layer (see \S3) \setlength{\unitlength}{1cm}
\vspace{-0.7cm} 
\begin{equation}
\begin{picture}(7,1.5) \put(.1,.3){\line(1,1){.4}}
\put(.5,.7){\line(1,-1){.4}} \put(4.92,.3){\line(5,2){1}}
\put(5.55,.3){\line(1,1){.4}}\put(5.95,.7){\line(1,-1){.4}}\put(6,.67){%
\line(5,-2){1}}\put(.1,.25){\circle*{.1} } \put(.9,.25){\circle*{.1} }
\put(1.85,.25){\circle*{.1} } \put(1.85,.25){\line(6,5){.5}}
\put(2.4,.25){\circle*{.1} } \put(2.4,.25){\line(0,1){.4}}
\put(2.9,.25){\circle*{.1} } \put(2.4,.68){\line(6,-5){.5}}
\put(3.9,.25){\circle*{.1} } \put(3.9,.68){\circle*{.1} }
\put(3.9,.25){\line(0,1){.4}} \put(4.9,.25){\circle*{.1} }
\put(5.5,.25){\circle*{.1} } \put(6.35,.25){\circle*{.1} }
\put(7.05,.25){\circle*{.1} } \put(2.4,.68){\circle*{.1} }
\put(.5,.68){\circle*{.1} } \put(5.95,.68){\circle*{.1} }
\put(7.5,0.1){$\longleftarrow S$} \put(7.3,0.6){$\longleftarrow S/\sigma$}
\end{picture}  \label{tree}
\end{equation}
\noindent or the parenthetical expression 
\[
(12)(345)(6)(789\overline{{10}}) 
\]
to denote the partition $\sigma =(B_{1},\ldots ,B_{4})$ of $[10]$ with 
\[
B_{1}=(1,2),B_{2}=(3,4,5), B_{3}=(6),B_{4}=(7,8,9,10) 
\]
(we regard a block as an ordered set). We identify a partition with the
corresponding equivalence relation on $S$ and the ordered quotient set $%
S/\sigma$ with $(B_{1},\ldots,B_{q})$. We say that a block is a \emph{%
singleton} if it has only one element, and that its element is \emph{unary}.

There is an alternative approach to partitions that is useful. We define an
(abstract) partition $\sigma $ of an ordered set $S$ to be an increasing map 
$f_{\sigma}:S\to S_{\sigma}$ of $S$ onto an ordered set $T=S_{\sigma}.$ We
then have the ordered partition $(B_{t})_{t\in T},$ where $
B_{t}=f_{\sigma}^{-1}(t).$ Conversely given a partition $\sigma=(B_{1},%
\ldots ,B_{q})$ in our initial sense, we have a corresponding increasing
surjection $f=f_{\sigma}:S\rightarrow S/\sigma=(B_{1},\ldots ,B_{q}),$ where $%
f_{\sigma}(x)=B_{j}$ if $x\in B_{j}.$ (\ref{tree}) may be regarded as the mapping
diagram of $f_{\sigma}$ in that example.

Given a second ordered partition $\pi =(C_{1},\ldots ,C_{r})$ of $S,$ we
write $\sigma \preceq \pi $ if every set $B_{j}$ is contained in some set $%
C_{k},$ i.e., $\sigma$ is a \emph{refinement} of $\pi.$ Equivalently, $
f_{\pi}=g\circ f_{\sigma}$ for some increasing function $g:S_{\sigma}%
\rightarrow S_{\pi}.$ We write $\sigma\prec \pi$ if $\sigma \preceq \pi $
and $\sigma \neq\pi $. Given partitions $\sigma$ of $S$ and $\pi$ of $T,$  
\[
\sigma\sqcup \pi=\left\{ B\sqcup C:B\in \sigma,C\in \pi\right\} 
\]
is an ordered partition of $S\sqcup T$ (we use the left to right ordering).

A \emph{colored partition} $(\mathcal{\sigma },c_{\mathcal{\sigma }})$ of a
colored ordered set $(S,c_{S})$ is an ordered partition $\mathcal{\sigma }
=(B_{1},\ldots ,B_{q})$ of $S$ together with a coloring $c _{\sigma }$ of $%
(B_{1},\ldots ,B_{q}),$ such that if $B_{j}=(x),$ then $c _{\sigma
}(B_{j})=c _{S}(x),$ i.e., each singleton has the same color as its unique
element. Equivalently, we have a colored ordered set $S_{\sigma }$ and an
order preserving surjection $f=f_{\sigma }:S\rightarrow S_{\sigma }$ with
the property that if $x$ is unary, then $c (f(x))=c (x).$ Given colored
ordered sets $S$ and $T$, we say that a mapping $f:S\rightarrow T$ is \emph{%
proper}, and write $f:S\twoheadrightarrow T$ if it satisfies these
conditions, i.e., it is an order-preserving surjection with the singleton color
condition. We say that $\sigma $ is a \emph{partition of }$S$ \emph{with the
coloring} $w=c _{\sigma }$ and that $(S/\sigma,w)$ is a colored ordered set.
We have that $(\mathcal{\sigma }\sqcup \mathcal{\mathcal{\tau }},vw)$ is a
colored ordered partition of the colored ordered set $S\sqcup T$.

We use the colored planar forest {\setlength{\unitlength}{1cm}{\ \vspace{%
-0.5cm} 
\begin{equation}
\begin{picture}(7,1.7) \put(.1,.25){\line(1,1){.4}}
\put(.65,.65){\line(1,-1){.4}} \put(2,.3){\line(2,1){.8}}
\put(2.95,.27){\line(0,1){.37}} \put(3.07,.69){\line(2,-1){.83}}
\put(4.95,.27){\line(0,1){.36}} \put(6.9,.25){\line(1,1){.4}}
\put(7.5,.67){\line(1,-1){.4}} \put(6,.3){\line(5,2){1.2}}
\put(7.55,.8){\line(5,-2){1.3}} \put(.55,.8){\circle{.4}}
\put(2.95,.8){\circle{.4}} \put(4.95,.8){\circle{.4}}
\put(7.38,.8){\circle{.4}} \put(.45,.7){2} \put(2.85,.7){1} \put(4.85,.7){2}
\put(7.3,.7){3} \put(.1,.1){\circle{.4}} \put(1.1,.1){\circle{.4}}
\put(2,.1){\circle{.4}} \put(2.95,.1){\circle{.4}}
\put(3.95,.1){\circle{.4}} \put(4.95,.1){\circle{.4}}
\put(5.95,.1){\circle{.4}} \put(6.95,.1){\circle{.4}}
\put(7.95,.1){\circle{.4}} \put(8.9,.1){\circle{.4}} \put(0,0){3}
\put(1,0){1} \put(1.9,0){2} \put(2.85,0){3} \put(3.85,0){3} \put(4.85,0){2}
\put(5.85,0){1} \put(6.85,0){2} \put(7.85,0){2} \put(8.85,0){3}
\put(9.5,0){$S$} \put(9.2,.7){$S/\sigma$} \end{picture}  \label{onelayer}
\end{equation}
} \vspace{0.3cm} {\normalsize \noindent }}or the parenthetical expression {%
{\normalsize 
\[
{(}31{)}_{2}{(}233{)}_{1}{(}2{)}_{2}(1223)_{3} 
\]
}}to indicate the colored ordered partition $\sigma =((12),(345),(6),(789%
\overline{{10}}{))},$ $w=c_{\sigma }=2123$ of the colored set 
\[
([10],3123321223). 
\]
In this example we have the colored order isomorphism 
\[
([10]/\sigma ,w)\cong ([4],2123). 
\]

If $S$ is a colored ordered set, we define $\mathcal{Y}(S)$ (respectively $%
\mathcal{Y}_{q}(S)$) to be the collection of all colored ordered partitions $%
\sigma $ of the colored set $S$ (respectively, with $q$ blocks). We
partially order $\mathcal{Y}(S)$ by $(\sigma ,v)\preceq (\pi ,w)$ if (1) $%
\sigma \preceq \mathcal{\pi },$ and (2) for any $B\in \sigma \cap \pi ,$ $%
c_{\sigma }(B)=c_{\pi }(B)$ (i.e., the common block $B$ has the same color
in either partition). For simplicity we simply write $\sigma \preceq 
\mathcal{\pi }$. Owing to the second condition, if $\sigma \preceq \mathcal{%
\pi }$ and $\pi \preceq \sigma $, then $\sigma =\pi $ as colored sets. It is
evident that $\sigma \preceq \pi $ if and only if $f_{\pi }=g\circ f_{\sigma
},$ for a (necessarily unique) proper function $g:S_{\sigma
}\twoheadrightarrow S_{\pi }$. If $S=(x_{1},\ldots ,x_{p})$ has the coloring 
$v=v(1)\cdots v(p)\in \Gamma ^{*},$ then $\mathcal{Y}(S)$ has the minimum
element $\mathbf{0}_{v}=(((x_{1})\ldots (x_{p})),v)$ and the maximal
elements $\mathbf{1}_{j}=((x_{1}\ldots x_{p}),j)$ for each $j\in \Gamma $.

We turn next to segments of colored partitions $P=[\sigma ,\tau ]$, where $%
\sigma \preceq \tau \in \mathcal{Y}(S).$ $P$ has the relative partial
ordering $\preceq $, and an element $\lambda \in P$ may be regarded as a
colored ordered set, as a colored partition of the elements in $S/\sigma $,
or as a colored partition of $S$. We say that segments $P\subseteq \mathcal{Y%
}(S)$ and $Q\subseteq \mathcal{Y}(T)$ are \emph{isomorphic} and write $%
P\simeq Q,$ if there exists an order isomorphism $\theta :P\rightarrow Q$
such that for each $\lambda \in P,$ $\lambda $ and $\theta (\lambda )$ are
isomorphic colored ordered sets. In particular, for any partitions $\sigma
\in \mathcal{Y}(S)$ and $\tau \in \mathcal{Y}(T)$, the segments $P=[\sigma
,\sigma ]$ and $Q=[\tau ,\tau ]$ are isomorphic if and only if $\sigma $ and 
$\tau $ are isomorphic colored ordered sets. Given a coloring $v$ of $[p]$
and $j\in [p]$, we let $Y_{v}^{j}=[\mathbf{0}_{v},\mathbf{1}_{j}]$. If $%
\sigma \preceq \pi $ and $\sigma ^{\prime }\preceq \pi ^{\prime },$ then we
have a unique order isomorphism 
\begin{equation}
\theta :[\sigma ,\pi ]\times [\sigma ^{\prime },\pi ^{\prime }]\simeq
[\sigma \sqcup \sigma ^{\prime },\pi \sqcup \pi ^{\prime }], \label{cart}
\end{equation}
where for each $(\lambda ,\lambda ^{\prime })$, $\theta ((\lambda ,\lambda ^{\prime }))$
and $ (\lambda ,\lambda ^{\prime })$ are isomorphic colored ordered sets.

\begin{lemma}
Let us suppose that $S$ is a colored ordered set and that $(\sigma
,v)\preceq (\pi ,w)$ in $\mathcal{Y}(S).$ Then letting $\pi =(C_{1},\ldots
,C_{q})$, and $v_{k}=v|C_{k}$ we have a natural isomorphism 
\begin{equation}
\theta :[\sigma ,\pi ]\cong Y_{v_{1}}^{w(1)}\times \ldots \times %
Y_{v_{q}}^{w(q)}  \label{intervals}
\end{equation}
where for each $\lambda \in [\sigma ,\pi ],$ $\lambda $ and $\theta (\lambda
)$ $=(\lambda _{1},\ldots ,\lambda _{q})$ are isomorphic colored ordered
sets.
\end{lemma}

\proof
Consider the mapping $g:S_{\sigma }\twoheadrightarrow S_{\pi } $ described
above. Let us identify $S_{\pi }$ with $[q]$. The intermediate colored
partitions correspond to factorizations $S_{\sigma }\twoheadrightarrow
T\twoheadrightarrow [q].$ To construct such a diagram it suffices to
choose for each $j\in [q]$ a factorization $g^{-1}(j)\twoheadrightarrow
T_{j}\twoheadrightarrow \{j\},$ i.e., an element $\lambda _{j}$ of $%
Y_{v_{j}}^{j},$ where $v_{j}$ is the coloring of the interval $g^{-1}(j).$
It is evident that we have a one-to-one order preserving correspondence $%
\lambda \leftrightarrow (\lambda _{1},\ldots ,\lambda _{q})$ with the
desired coloring property. \endproof

If $S=((x),j)$, then $\pi_{j}=((x),j)$ is the only element in $\mathcal{Y}%
(S) $, and the only segment is isomorphic to $Y_{j}^{j}$. In order to obtain
the incidence Hopf algebras, it is necessary to impose a more inclusive
equivalence relation, which identifies all such intervals with a
multiplicative identity. Given a colored partition $(\sigma ,v)$ of $S,$ we
let $S_{ns}$ be the union of the non-singleton blocks in $\sigma ,$ $\sigma
_{ns}$ be the collection of non-singleton blocks, and $v_{ns}$ be the
restriction of $v$ to $\sigma _{ns} $. We say that order intervals $%
P=[\sigma ,\tau ]\subseteq \mathcal{Y}(S)$ and $Q=[\sigma ^{\prime },\tau
^{\prime }]\subseteq \mathcal{Y}(T)$ are \emph{similar,} and write $P\sim Q$%
, if there is an order-preserving bijection $\theta :P\rightarrow Q$%
\thinspace such that for each $\lambda \in P,$ $\lambda _{ns}$ and $\theta
(\lambda )_{ns}$ are isomorphic colored sets. If the non-singleton sets $%
S_{ns}$ and $T_{ns}$ are empty, the latter is a vacuous restriction. 

For any
colored partition $\sigma $ and segment $P,$ 
\[
\lbrack \sigma ,\sigma ]\times P\sim P. 
\]
To prove this, consider the mapping 
\[
\theta :[\sigma ,\sigma ]\times P\rightarrow P:\sigma \sqcup \lambda \mapsto
\lambda 
\]
This is clearly a bijection and order preserving. Since $(\sigma \sqcup
\lambda )_{ns}=(\lambda )_{ns},$ $\theta \ $satisfies the coloring
condition. Similarly we have that $P\times [\sigma ,\sigma ]\sim P$ for any $%
\sigma $ and $P.$ Finally it is easy to verify that if $P\sim P^{\prime }$
and $Q\sim Q^{\prime },$ then 
\begin{equation}\label{symproduct}
P\times Q\sim P^{\prime }\times Q^{\prime }.
\end{equation}

We let $\mathcal{P}=\mathcal{P}^{N}$ denote the collection of all similarity
classes $P_{\sim}$ of segments
$P=[\sigma ,\pi ]\subseteq \mathcal{Y}(S)$ for arbitrary finite $N$-colored sets $S.$ We define an
associative product on $\mathcal{P}$ by 
\begin{equation}
[\sigma ,\tau ]_{\sim }[\sigma ^{\prime },\tau ^{\prime }]_{\sim
}=[\sigma \times \sigma ^{\prime },\tau \times \tau ^{\prime }]_{\sim }. 
\label{monoid}
\end{equation}
It follows from (\ref{symproduct}) that this is well-defined. The corresponding multiplicative
identity is given by
$1=[Y_{i}^{i}]$, where 
$i\in \langle N \rangle$ is arbitrary.

The intervals $Y_{u}^{i}$ with $\left| u\right| >1$ are all non-similar.
When confusion is unlikely, we will dispense with the similarity class
notation $[\,\,]_{\sim }$. It is evident that $\mathcal{P}$ is just the free
monoid on the symbols $Y_{u}^{i}$ with $\left| u\right| >1.$

Given a segment $P=[\sigma ,\tau ]$ in $Y_{v}^{i}$ we define a \emph{chain} $%
\gamma $ in $P$ to be a sequence $\sigma =\sigma _{0}\prec \sigma _{1}\prec
\ldots \prec \sigma _{r}=\tau .$ Given such a chain, the blocks of $\sigma
_{k-1}$ form a colored partition of the blocks in $\sigma _{k}.$ Thus the
interval $[\sigma _{k-1},\sigma _{k}]$ determines a one-layer colored forest 
$L_{k}$ (see (\ref{tree})). The colored layers $L_{1},\ldots ,L_{r}$ are
consecutively compatible, and letting $v$ be the colors of $L_{1}$ and $w$ be the 
colors of the roots of $L_{r}$, we may
associate the layered forest 
\[
F(\gamma )=L_{1}\vartriangleright L_{2}\vartriangleright \ldots
\vartriangleright L_{r} \in \mathbf{LF}_{v}^{w}(N)
\]
with $\gamma .$ Conversely for each $F\in \mathbf{LF}_{v}^{w}$, we may use the
layers of $F$ to determine a unique chain $\gamma (F)$ in $[\sigma,\tau]$ with
$F(\gamma (F))=F.$

We define the \emph{length} of a chain $\gamma $ in a segment $P$ to be the number 
of intervals in $\gamma$, or equivalently the number
of layers in $F(\gamma).$ The \emph{grading} $\rho (P)$ of $P$ is the
maximal length of a chain in $P.$ If $P=Y_{u}^{i},$ any chain $\gamma in P$ 
can be extended to a chain $\gamma ^{\prime }$ such that for each $
k,\sigma _{k-1}^{\prime }$ is obtained from $\sigma _{k}^{\prime }$ by
splitting precisely one of the blocks of $\sigma _{k}$ into two blocks. The
corresponding colored trees are simple, layered, and each non-degenerate vertex has two
offspring. Examining the tree of such a chain it follows that the maximal
chains in $Y_{u}^{i}$ all have $\left| u\right| -1$ elements, and thus $\rho
(Y_{u}^{i})=\left| u\right| -1.$ Given an arbitrary segment $P=[\sigma ,\tau
]\cong Y_{v_{1}}^{w(1)}\cdots Y_{v_{q}}^{w(q)},$ it is evident that 
\[
\rho (P)=\rho (Y_{v_{1}}^{w(1)})+\ldots +\rho (Y_{v_{q}}^{w(q)}). 
\]

\section{Hopf algebras}

We briefly recall some elementary notions from the theory of Hopf algebras.
More complete discussions can be found in \cite{Ka}, \cite{Sw}, \cite{F},
and \cite{M}.

Given a vector space $V,$ we let $L(V)$ denote the algebra of all linear
mappings $T:V\rightarrow V$. The identity mapping $I:V\rightarrow V$ is a
multiplicative identity for $L(V)$. If $(A,1)$ is a unital algebra, the 
\emph{tensor product algebra }$A\otimes A$ is given the associative
multiplication 
\[
(x_{1}\otimes y_{1})(x_{2}\otimes y_{2})=x_{1}x_{2}\otimes y_{1}y_{2} 
\]
and the multiplicative unit $1\otimes 1.$ A \emph{bialgebra} $(H,m,\eta
,\Delta ,\varepsilon )$ consists of a vector space $A$ with an associative
product $m:H\otimes H\rightarrow H,$ a homomorphism $\eta :$ $\Bbb{C}
\rightarrow H:\alpha \rightarrow \alpha 1$, where $1$ is a multiplicative
unit for $H,$ a coassociative coproduct $\Delta :H\rightarrow H\otimes H,$
and a counit $\varepsilon :H\rightarrow \Bbb{C}$ with the linking property
that $\Delta :H\rightarrow H\otimes H$ is a unital homomorphism. We employ
Sweedler's notation 
\[
\Delta a=\sum_{(a)}a_{(1)}\otimes a_{(2)}. 
\]

An \emph{antipode} for a bialgebra $H$ is a mapping $S:H\rightarrow H$ such
that for any $a\in H$ 
\[
\sum_{(a)} S(a_{(1)})a_{(2)}=\sum_{(a)}a_{(1)}S(a_{(2)})=\varepsilon (a)1. 
\]
or equivalently, $m(S\otimes I)\Delta =m(I\otimes S)\Delta =\eta \circ
\varepsilon .$ We say that $S$ is a \emph{left antipode} if one just has the
first and third terms are equal, and a \emph{right antipode}, if one has the
second equality. If $H$ has an antipode, then any left (respectively right)
antipode automatically coincides with $S,$ and in particular, $S$ is unique.
An antipode $S$ is automatically a unital antihomorphism, i.e., we have 
\begin{eqnarray*}
S(gh) &=&S(h)S(g) \\
S(1) &=&1
\end{eqnarray*}
(see \cite{Sw}, Prop. 4.0.1). A Hopf algebra $(H,m,\eta ,\Delta ,\varepsilon
,S)$ is a bialgebra $(H,m,\eta ,\Delta ,\varepsilon )$ together with an
antipode $S.$

Given $\varphi ,\psi \in L(H),$ we define the \emph{convolution }$\varphi
*\psi \in L(H)$ by 
\[
\varphi *\psi (a)=\sum_{(a)}\varphi (a_{(1)})\psi (a_{(2)}) 
\]
or equivalently, $\varphi *\psi =m\circ (\varphi \otimes \psi )\circ \Delta
. $ This determines an associative product on $L(H)$ with the multiplicative
unit $u=\eta \circ \varepsilon .$ It is evident form the definition that $S$
is an antipode if and only if $S*I=I*S=u,$ i.e., it is the convolution
inverse of $I.$

Given any vector space $V,$ we define the \emph{flip} $\tau :V\otimes
V\rightarrow V\otimes V$ by $\tau (v\otimes w)=w\otimes v.$ Given a
bialgebra $H=(H,m,1,\Delta ,\varepsilon ),$ we let $m^{op }=m\circ \tau
$ and $\Delta ^{cop }=\tau \circ \Delta $ and We
define the 
\emph{opposite }and \emph{co-opposite } bialgebras to be $%
H^{op}=(H,m^{op},1,\Delta ,\varepsilon )$ and $H^{cop}=(H,m,1,\Delta ^{op
},\varepsilon ).$ It is shown in \cite{Ka} Cor. III.5.5 that these are
indeed bialgebras, and the antipode $S:S:H^{op}\rightarrow H^{cop}$ is a
Hopf algebra homomorphism. The latter corresponds to the fact that the
antipode $S:H\to H$ is both an anti-endormorphism and an
anti-coendomorphism. If the antipode $S:H\rightarrow H$ is a linear
isomorphism, then $S^{-1}$ is the antipode for both $H^{op}$and $H^{cop},$
and $S:H^{op}\rightarrow H^{cop}$ is a Hopf algebraic isomorphism. We refer
to either Hopf algebra as the \emph{inverse} Hopf algebra of $H$.

The following will play an important role in what follows.

\begin{proposition}\label{Magidth} If $H$ is a bialgebra for which both $H$ and $H^{cop}$
have antipodes
$S$ and $S^{cop}$, respectively, then $S$ is invertible,
and $S^{cop}=S^{-1}$.\end{proposition}
\proof See \cite{Mag}, Exercise 1.3.3 (the solution is provided). \endproof

A bialgebra $(H,m,\eta ,\Delta ,\varepsilon ,S)$ is said to be \emph{%
filtered }if one has an increasing sequence of subspaces $H^{n}$ $(n\geq 0)$
with $\bigcup H^{n}=H,$ for which 
\begin{eqnarray*}
H^{m}H^{n} &\subseteq &H^{m+n} \\
\Delta H^{n} &\subseteq &\sum_{p+q=n}H^{p}\otimes H^{q}.
\end{eqnarray*}
and it is said to be \emph{connected} if in addition $H^{0}=\Bbb{C}1.$ We
let $H_{+}$ denote the kernel of $\varepsilon$, and $H^{n}_{+}=H^{n}\cap H_{+}$.

A connected filtered bialgebra is automatically a Hopf algebra, i.e., it
has an antipode. The following is proved in \cite{F} (see also \cite{S}, p. 238).

\begin{theorem}\label{geometricth}
If $H$ is a connected filtered bialgebra, then it has an antipode given by
the ``geometric series'' 
\begin{equation}
Sa=(u-(u-I))^{-1}(a)=\sum_{k=0}^{\infty }(u-I)^{*k}a.  \label{geometrics}
\end{equation}
The sum is finite for each $a$ since if $a\in H^{n},$ then $%
(u-I)^{*(n+1)}(a)=0.$ 
\end{theorem}

\begin{corollary} The antipode of a connected filtered Hopf algebra is invertible.
\end{corollary}
\proof
We have that $H^{cop}$ is a filtered bialgebra since 
\[
\Delta^{cop} H^{n} \subseteq \sum_{p+q=n}H^{q}\otimes H^{p}.
\]
Thus from Theorem \ref{geometricth}, $H^{cop}$ has an antipode and from Proposition
\ref{Magidth}, $S$ is invertible.
\endproof

The proof of Theorem \ref{geometricth} in \cite{F} is based upon the following lemma.
We have included a proof since we will need the calculation in the discussion
that follows.

\begin{lemma}
Suppose that $H$ is a connected Hopf algebra. Then for any $a\in H_{+}^{n},$ 
\[
\Delta a=a\otimes 1+1\otimes a+y
\]
where $y\in H_{+}^{n-1}\otimes H_{+}^{n-1}$
\end{lemma}

\begin{proof}
Let us define $y$
by this relation. Since $\Delta a\in \sum H^{k}\otimes H^{n-k},$ 
\[
y=\Delta a-a\otimes 1-1\otimes a\in \sum_{k=0}^{n} H^{k}\otimes H^{n-k},
\]
and we have  $y=\sum b_{k}\otimes c_{k},$ where $b_{k}\in H_{k}$ and $%
c_{k}\in H_{n-k}.$ Applying the right coidentity relation to $a\in H_{+},$ 
\[
a=\varepsilon \otimes id(\Delta a)=\varepsilon (a)1+a+\sum \varepsilon
(b_{k})c_{k}=a+\sum \varepsilon (b_{k})c_{k},
\]
and thus $\sum \varepsilon (b_{k})c_{k}=0.$ It follows that
\[
y=\sum (b_{k}-\varepsilon (b_{k})1)\otimes c_{k}=\sum b_{k}^{\prime }\otimes
c_{k}
\]
where $b_{k}^{\prime }=b_{k}-\varepsilon (b_{i})1\in H_{+}^{k}.$ Similarly,
from the right coidentity relation, 
\[
y=\sum b_{k}^{\prime }\otimes c_{k}^{\prime }
\]
with $c_{k}^{\prime }=c_{k}-\varepsilon (c_{k})1\in H_{+}^{n-k}.$ Since $%
H_{+}^{0}=\{0\},$ we obtain 
\begin{equation}
\Delta y\in \sum_{k=1}^{n-1}H_{+}^{k}\otimes H_{+}^{n-k}\subseteq
H^{n-1}_{+}\otimes H^{n-1}_{+}.\label{recursion}
\end{equation}
\end{proof}

There is a simple recursive characterization of the antipode in a connected
filtered Hopf algebra. If $a\in H_{+}^{n}$, then from (\ref{recursion}), 
\begin{equation}
\Delta a=a\otimes 1+1\otimes a+\sum_{k=1}^{n-1}a_{k}\otimes b_{n-k}
\label{deltared}
\end{equation}
where $a_{k}\in H^{k}$, $b_{n-k}\in H^{n-k},$ and thus 
\[
0=\varepsilon (a)1=S(a)+a+\sum_{k=1}^{n-1} S(a_{k})b_{n-k}. 
\]
$S$ is thus recursively determined by $S(1)=1$ and if $a\in H_{n}$ with $%
\varepsilon (a)=0,$ then 
\begin{equation}
S(a)=-a-\sum_{k=1}^{n-1}S(a_{k})b_{n-k}, \label{inducantip}
\end{equation}
where $a_{k}\in H^{k}$ and $b_{n-k}\in H^{n-k}$.

Restricting our attention to the co-opposite algebra, the antipode $S^{-1}$
is characterized by the relations 
\[
\sum_{(a)}S^{-1}(a_{(2)})a_{(1)}=\sum_{(a)} a_{(2)}S^{-1}(a_{(1)})=\varepsilon (a)1. 
\]
In particular assuming that $H$ is a connected filtered Hopf algebra, it is
recursively determined by $S^{-1}(1)=1$ and if $a\in H^{n}_{+}$ then 
\begin{equation}
S^{-1}(a)=-a-\sum_{k=1}^{n-1}b_{k}S^{-1}(a_{n-k}),  \label{inverserec}
\end{equation}
where $a_{k}\in H^{k}$, $b_{n-k}\in H^{n-k}.$

A bialgebra $H$ is \emph{graded} if there are subspaces $H_{n}$ of $H$
with
$ H_{n}\cap H_{m}=\left\{ 0\right\} $ and $\sum H_{n}=H$ such that 
\begin{eqnarray*}
H_{m}H_{n} &\subseteq &H_{m+n} \\
\Delta H_{n} &\subseteq &\sum_{p+q=n}H_{p}\otimes H_{q}.
\end{eqnarray*}
It is immediate that the subspaces
$ H^{n}=\sum_{i=0}^{n}H_{i}$
determine a filtration of $H$.

Finally let us suppose that $H$ is a Hopf algebra and that $\theta
:H\rightarrow H$ is a linear isomorphism such that $\theta (1)=1$ and $%
\varepsilon \circ \theta =\varepsilon .$ It is evident that $H^{\theta
}=(H,m^{\theta },1,\Delta ^{\theta },\varepsilon ,S^{\theta })$ is again a
Hopf algebra, where $m^{\theta }=\theta \circ m\circ (\theta ^{-1}\otimes
\theta ^{-1}),$ $\Delta ^{\theta }=(\theta \otimes \theta )\circ \Delta
\circ \theta^{-1} ,$ and $S=\theta \circ S\circ \theta ^{-1}.$ We refer to $%
H^{\theta }$ as the $\theta $\emph{-transformed Hopf algebra}. In all of the
cases considered below, $\theta$ is involutory, i.e., $\theta^{2}=id$
and thus $\theta^{-1}=\theta$.

\section{The interval and Lagrange Hopf incidence algebras}

We begin by constructing the relevant incidence Hopf algebra (see \cite{S} for
more details). As in \S4, we let
$\mathcal{P}=\mathcal{P}^{N}$ be the collection of all similarity classes of segments
$[\sigma,\tau]$ where $\sigma \preceq \tau$ are 
colored ordered partitions of $N$-colored sets. To be more explicit, we
may assume that each $[\sigma,\tau]$ is a subset of some $\mathcal{Y }(S_{u})$, 
where $S_{u}$ is the ordered set $[r]=(1,\ldots ,r)$ with the $N$-coloring $u.$ We
define the \emph{interval partition incidence bialgebra} $\mathcal{H}=\mathcal{H}^{N}$ to be the
vector space with basis (labelled by) the elements of
$\mathcal{P}$. The monoid operation on $\mathcal{P}$ determined by (\ref{monoid}) determines a
multiplication on $\mathcal{H}$. It is evident that $%
\mathcal{H}$ is the free unital algebra on the segments $Y_{u}^{i}$ $%
(1\leq i\leq N,\left| u\right| \geq 1).$ Following Joni and Rota \cite{J},
p. 98, we define a coproduct
\[
\Delta :\mathcal{H}\rightarrow \mathcal{H}\otimes \mathcal{H}%
\]
by letting  
\[
\Delta (P)=\sum_{0_{P}\preceq\pi\preceq 1_{P}}[0_{P},\pi ]\otimes [\pi ,1_{P}],
\]
for $P \in \mathcal{P}$, and then extending linearly to $\mathcal{H}.$ Given $P_{1}$ and
$P_{2}\in
\mathcal{P},$ the equality $%
\Delta (P_{1}P_{2})=\Delta (P_{1})\Delta (P_{2})$ is a simple consequence of 
(\ref{monoid}), hence $\Delta$ is a homomorphism. We also have a
coidentity homomorphism
$\varepsilon :\mathcal{H}\rightarrow \Bbb{C}$ determined by
$\varepsilon (1)=1$ and $\varepsilon (Y_{v}^{i})=0,$ for $|u|\geq 2$ and we see that
$\mathcal{H}$ is a bialgebra. 
As an algebra, $\mathcal{H}$ is
the free associative unital algebra generated by the segments $Y_{u}^{i}$ $(1%
\leq i\leq N,\left| u\right| \geq 2)$ with unit $1.$ In order to simplify
the notation below we let $Y_{j}^{i}=\delta _{ij}1$ if $i,j\in \langle
N\rangle .$ 

Given $(\pi ,w)\in Y_{u}^{i}$, $\pi =(C_{1},\ldots ,C_{q})$ is an ordered
partition of the colored set $([p],u)$ where $p=\left| u\right| ,$ and $%
w=w(1)\cdots w(q)$ is a coloring of $\pi .$ We let $u_{i}=u|C_{i}.$ From
Lemma 1, $[\mathbf{0}_{u},\pi ]\thicksim Y_{u_{1}}^{v(1)}\ldots
Y_{u_{q}}^{v(q)}$ and $[\pi ,\mathbf{1}_{i}]$ $=[\mathbf{0}_{v},\mathbf{1}%
^{i}]\thicksim Y_{v}^{i}$ from which we conclude that 
\begin{equation}
\Delta (Y_{u}^{i})=\sum_{q=1}^{p}\,\,\,\sum_{\pi \in \mathcal{Y}%
_{q}(p)}\,\,\,\sum_{v\in \langle N\rangle ^{q}}Y_{u_{1}}^{v(1)}\ldots
Y_{u_{q}}^{v(q)}\otimes Y_{v}^{i}  \label{comult}
\end{equation}
Owing to our convention that $Y_{j}^{i}=\delta _{ij}1,$ for $i,j\in \langle
N\rangle$, many of these terms are zero. For example if we let $q=1,$ then $%
\mathcal{Y}_{1}(p)=\left\{ ((1,\ldots ,p),j)\right\} $ and we have only the
summand 
\[
\,\sum_{v\in \langle N\rangle }Y_{u}^{v}\otimes Y_{v}^{i}=\sum_{v}\delta
_{vu}1\otimes Y_{v}^{i}=1\otimes Y_{u}^{i}. 
\]
On the other hand if $q=p,$ then $\mathcal{Y}_{q}(p)=\left\{ (((1),\ldots
,(p)),u(1))\ldots u(p))\right\} $ and 
\[
\sum_{v\in \langle N\rangle ^{p}}Y_{u(1)}^{v(1)}\ldots
Y_{u(p)}^{v(p)}\otimes Y_{v}^{i}=\sum_{v\in \langle N\rangle ^{p}}\delta
_{u(1)}^{v(1)}\ldots \delta _{u(p)}^{v(p)}1\otimes Y_{v}^{i}=1\otimes
Y_{u}^{i}. 
\]
It follows that 
\begin{equation}
\Delta (Y_{u}^{i})=Y_{u}^{i}\otimes 1+1\otimes
Y_{u}^{i}+\sum_{q=2}^{p-1}\,\,\,\sum_{\pi \in \mathcal{Y}_{q}(p)}\,\,\,%
\sum_{v\in \langle N\rangle ^{q}}Y_{u_{1}}^{v(1)}\ldots
Y_{u_{q}}^{v(q)}\otimes Y_{v}^{i}\label{expanddelta} 
\end{equation}

We define $\mathcal{H}_{n}$ to be the linear subspace spanned by the
segments $P$ with $\rho (P)=n$ (see \S4). Given a generator $%
Y_{u}^{i}\in \mathcal{H}_{n}$, (i.e., with $|u|=n+1$), a partition $\pi =(C_{1},\ldots ,C_{q})$ of
$[n+1],$ and $v\in \langle N\rangle ^{q},$ 
\begin{eqnarray*}
\rho (Y_{u_{1}}^{v(1)}\ldots Y_{u_{q}}^{v(q)})+\rho (Y_{v}^{i}) &=&(\left|
u_{1}\right| -1)+\ldots +\left| u_{q}\right| -1)+(\left| v\right| -1) \\
&=&\left| u\right| -q+(q-1)=n,
\end{eqnarray*}
and thus 
\[
\Delta (Y_{u}^{i})\in \sum_{p+q=n}\mathcal{H}_{p}\otimes \mathcal{H}_{q}. 
\]
For any $y\in \mathcal{H}_{p}$, $y^{\prime }\in \mathcal{H}_{p^{\prime }}$, $%
z\in \mathcal{H}_{q},z^{\prime }\in \mathcal{H}_{q^{\prime }},$ 
\[
(x\otimes y)(x^{\prime }\otimes y^{\prime })=xx^{\prime}\otimes yy^{\prime}\in
\mathcal{H}_{p+p^{\prime }}\otimes \mathcal{H}_{q+q^{\prime }}.
\] Since $\Delta $ is a multiplicative
homomorphism, we conclude that if $P$ is an arbitrary interval, i.e., a
product of terms of the form $Y_{u}^{i}$, and $\rho (P)=m,$ then 
\[
\Delta (P)\in \sum_{p^{\prime \prime }+q^{\prime \prime }=m}\mathcal{H}%
_{p^{\prime \prime }}\otimes \mathcal{H}_{q^{\prime \prime }} 
\]
and thus $\mathcal{H}$ is a graded and connected bialgebra. From the previous
section, $\mathcal{H}$ has an invertible antipode $S_{\mathcal H}$, and in particular it
is a Hopf algebra.

The following antipode formula of Schmitt may be regarded as a transcription
of (\ref{geometrics}) (see \cite{F}, \S 11.1 and \cite{S}, Th. 4.1). 
\begin{equation}
S_{\mathcal{H}}(Y_{u}^{i})=\sum_{k\geq 0}\,\,\,\sum_{\mathbf{0}_{u}=\sigma
_{0}\prec
\sigma _{1}\prec \ldots \prec \sigma
_{k}=\mathbf{1}_{i}}(-1)^{k}\prod_{h=1}^{k}[\sigma _{h -1},\sigma _{h}].
\label{antipode2}
\end{equation}
As we have seen in \S 4, there is a one-to-one correspondence between the
chains $\gamma=(\mathbf{0}_{u}=\sigma _{0}\prec \sigma _{1}\prec \ldots \prec
\sigma _{k}=\mathbf{1}_{i})$ in $Y_{u}^{i}$ and the colored layered trees 
\[
T(\gamma )=F_{1}\rhd \ldots \rhd F_{k}\in \mathbf{LT}_{u}^{i}(N),
\]
where the one layered forest $F_{h}$ corresponds to the interval $[\sigma
_{h-1},\sigma _{h}].$ The roots $x_{1}^{h}\ll \ldots \ll x_{q_{h}}^{h}$ in $%
F_{h}$ determine the factors in the
decomposition 
\[
[\sigma _{h-1},\sigma _{h}]=Y(x_{1}^{h})\ldots
Y(x^{h}_{q_{h}}).
\]
 Furthermore since we have identified the factors
$Y_{j}^{j}$ with the multiplicative identity in $\mathcal{H}$, we may employ
just the non-degenerate vertices. Relabelling the non-degenerate vertices in
$T(\gamma )$ by $x_{1}\ll \cdots
\ll x_{k}$ and letting 
\begin{equation}
\Omega(T(\gamma ))=\prod^{\ll }Y(x)=Y(x_{1})\ldots Y(x_{k}),
\end{equation}
we conclude that 
\[
\prod_{h=1}^{k}[\sigma _{h -1},\sigma _{h}]=\Omega (T(\gamma )), 
\]
and we may rewrite (\ref{antipode2}) in the form: 
\begin{equation}
S_{\mathcal{H}}(Y_{u}^{i})=\,\,\sum_{T\in \mathbf{LT}_{u}^{i}}(-1)^{\ell
(T)}\Omega (T).  \label{omegaantipode}
\end{equation}

We call the inverse Hopf algebra 
\[
\mathcal{L}=\mathcal{L}^{N}=(\mathcal{H}^{N},m,\varepsilon ,\Delta
^{op},\eta ,S_{\mathcal{L}}=S_{}^{-1}) 
\]
the \emph{left Lagrange} Hopf algebra. We may also regard $\mathcal{L}$ as a
transformed algebra of $\mathcal{H}$, as we next show.

Given a word $u=u_{1}\cdots u_{n}\in \langle N\rangle^{*},$ we let $u^{*}$
be its \emph{reflection} $u_{n}\ldots u_{1}$. Since $\mathcal{H}$ is
freely generated by the $Y_{u}^{j},$ the inclusion mapping 
\[
\mathbf{s}:Y_{u}^{i}\mapsto Y_{u^{*}}^{i}\in \mathcal{H}^{op} 
\]
extends to an involutory algebra isomorphism $\mathbf{s}:\mathcal{H}
\rightarrow \mathcal{H}^{op}.$ This may be regarded as an antiisomoprhism $%
\mathbf{s}:\mathcal{H}\rightarrow \mathcal{\mathcal{H}}$ satisfying 
\begin{equation}
\mathbf{s}(Y_{u_{1}}^{i_{1}}\cdots
Y_{u_{n}}^{i_{n}})=Y_{u_{n}^{*}}^{i_{n}}\cdots Y_{u_{1}^{*}}^{i_{1}}.
\end{equation}

\begin{lemma}
We have that $m^{op}=\mathbf{s}\circ m\circ (\mathbf{s\otimes s})$, $\Delta =(%
\mathbf{s}\otimes \mathbf{s})\circ \Delta \circ \mathbf{s}$, and
$S_{\mathcal L}=\mathbf{s} S_{\mathcal{H}}\mathbf{s}$.
\end{lemma}

\proof The first result is immediate since $\mathbf{s}$ is an antihomomorphism.

If $u=u(1)\ldots u(p)$ then $u^{*}=u(p)\ldots u(1)$ i.e., $%
u^{*}(k)=u(p+1-k)$ and 
\begin{eqnarray*}
(\mathbf{s}\otimes \mathbf{s})\circ \Delta \circ \mathbf{s}(Y_{u}^{i}) &=&(\mathbf{s%
}\otimes \mathbf{s})\circ \Delta (Y_{u^{*}}^{i}) \\
&=&(\mathbf{s}\otimes \mathbf{s})\sum Y_{u^{*}|B_{1}}^{v(1)}\cdots
Y_{u^{*}|B_{q}}^{v(q)}\otimes Y_{v}^{i} \\
&=&\sum_{v} Y_{(u^{*}|B_{q})^{*}}^{v(q)}\cdots
Y_{(u^{*}|B_{1})^{*}}^{v(1)}\otimes Y_{v^{*}}^{i}
\end{eqnarray*}
Given that $B_{1}=(1,\ldots ,j_{p(1)}),$ $\cdots ,B_{q}=(j_{p(q-1)}+1,\ldots
,j_{p(q)}=p)$ we have that 
\begin{eqnarray*}
u^{*}|B_{1} &=&u(p)u(p-1)\cdots u(p-j_{p(1)}+1), \\
u^{*}|B_{2} &=&u(p-j_{p(1)})\cdots u(p-j_{p(2)}+1 ) \\
&&\ldots \\
u^{*}|B_{q} &=&u(p-j_{p(q-1)})\cdots u(1)
\end{eqnarray*}
and thus 
\begin{eqnarray*}
(u^{*}|B_{q})^{*} &=&u(1)\cdots u(p-j_{p(q-1)})=u|C_{1} \\
&&\ldots \\
(u^{*}|B_{1})^{*} &=&u(p-j_{p(1)}+1)\cdots u(p)=u|C_{q}.
\end{eqnarray*}
where $C_{k}=B_{q-k+1}.$ If we let $w(k)=p+1-k,$ we conclude that 
\[
(\mathbf{s}\otimes \mathbf{s})\circ \Delta \circ \mathbf{s}(Y_{u}^{i})=\sum_{v}
Y_{u|C_{1}}^{v(1)}\cdots Y_{u|C_{q}}^{v(q)}\otimes Y_{v}^{i}=\Delta
(Y_{u}^{i}). 
\]
Since $\mathbf{s}$ is an antihomomorphism and $\Delta $ is a homomorphism, $(%
\mathbf{s}\otimes \mathbf{s})\circ \Delta \circ \mathbf{s}$ is a homomorphism and
thus the relation holds for all elements of $\mathcal{H}.$

Finally from (\ref{comult}),

\begin{eqnarray*}
\sum_{(u)} \mathbf{s}S_{\mathcal H}\mathbf{s}((Y_{u}^{i})_{(2)})(Y_{u}^{i})_{(1)} &=&\sum
_{v}\mathbf{s} (S_{\mathcal{H}}(Y_{(v_{n}\cdots v_{1})}^{i}))Y_{u_{1}}^{v(1)}\ldots
Y_{u_{q}}^{v(q)} \\ &=&\mathbf{s}(\sum_{v} Y_{u_{q}^{*}}^{v(q)}\ldots
Y_{u_{1}^{*}}^{v(1)}S_{\mathcal{H}}(Y_{v^{*}}^{i})) \\
&=&\mathbf{s}(\sum_{(u)} (Y_{u^{*}}^{i})_{(1)}S_{\mathcal{H}}((Y_{u^{*}}^{i})_{(2)})) \\
&=&\mathbf{s}(\varepsilon (Y_{u^{*}}^{i})1)=\delta _{u^{*}}^{i}=\delta
_{u}^{i},
\end{eqnarray*}

\endproof

\begin{corollary}
The left Lagrange algebra $\mathcal{L}$ is isomorphic to the 
$\mathbf{s}$-transformed algebra $\mathcal{H}^{\mathbf{s}}.$
\end{corollary}

\proof As we remarked in \S 5, the mapping 
\[
S_{\mathcal H}:(\mathcal{H},m^{op},1, \Delta ,S^{-1}_{\mathcal{H}},\varepsilon
)\rightarrow (\mathcal{H} ,m,1,\Delta ^{op},S^{-1}_{\mathcal{H}},\varepsilon
)=\mathcal{L} 
\]
is a Hopf algebraic isomorphism. \endproof

The left Lagrange algebra can also be regarded as the incidence algebra of
the \emph{reversed }interval partitions. For this purpose we define a new
ordering of colored partitions by letting $(\sigma ,v)\trianglelefteq (\tau
,w)$ if $(\tau ,w)\preccurlyeq (\sigma ,v).$ This does not effect the
algebraic structure of $\mathcal{H}.$ To be more specific, let us
notationally identify the segments $[\tau ,\sigma ]_{\trianglelefteq }$ with
the segments $[\sigma ,\tau ]_{\preccurlyeq }.$ The underlying algebra of
the Hopf algebra $\mathcal{H}_{\trianglelefteq }$ associated with this
reordered system is then identified with the free algebra on the symbols $%
Y_{u}^{k}$ with $\left| u\right| \geq 2.$ On the other hand the
comultiplication is determined by

\[
\Delta _{\trianglelefteq }([\tau ,\sigma ]_{\trianglelefteq })=\sum_{\tau
\trianglelefteq \mu \trianglelefteq \sigma }[\tau ,\mu ]_{\trianglelefteq
}\otimes [\mu ,\sigma ]_{\trianglelefteq }=\sum_{\sigma \preccurlyeq \mu
\preccurlyeq \tau }[\mu ,\tau ]\otimes [\sigma ,\mu ]=\Delta ^{op
}([\sigma,\mu]). 
\]

We have a related involutory anti-isomorphism $\mathbf{t}:\mathcal{H}%
\rightarrow \mathcal{H}$ determined by the identity mapping 
\[
\mathbf{t}:Y_{u}^{i}\mapsto Y_{u}^{i}\in \mathcal{L}^{op}, 
\]
or equivalently, 
\[
\mathbf{t}(Y_{u_{1}}^{i(1)}\cdots Y_{u_{n}}^{i(n)})=Y_{u_{n}}^{i(n)}\cdots
Y_{u_{1}}^{i(1)}. 
\]
We have that $\mathbf{t}=\alpha \mathbf{s}=\mathbf{s}\alpha ,$ where $\alpha $:$%
\mathcal{H\rightarrow \mathcal{H}}$ is the involutory \emph{automorphism}
determined by $\alpha (Y_{u}^{i})=Y_{u^{*}}^{i}.$

We define the \emph{right Lagrange algebra }$\mathcal{R}=\mathcal{R}^{N}$ to
be the $\mathbf{t} $-transformed Hopf algebra, i.e., 
\[
\mathcal{R}^{N}=(\mathcal{H}^{N},m^{\mathbf{t}},\varepsilon ,\Delta ^{\mathbf{t}%
},\eta ,\mathbf{t}S\mathbf{t}) 
\]
We note that $m^{\mathbf{t}}$ coincides with $m^{op}$ since $\mathbf{t}$ is an
antiisomorphism and 
\[
\mathcal{R}^{N}=((\mathcal{H}^{N})^{\mathbf{s}})^{\alpha }=(\mathcal{L}%
^{N})^{\alpha }. 
\]
Owing to the latter relation, we also refer to $\mathcal{R}^{N}$ as the 
\emph{reflection} of $\mathcal{L}^{N}$.

\section{The reduced tree formulae for the antipodes}

Given a reduced colored tree $T$, each $x\in \mathbf{V}(T)$ determines a
generator $Y(x)=Y_{v}^{j},$ where $j$ is the color of $x,$ and $u=u(1)\ldots
u(k)$ are the colors (in order) of its offspring. For any reduced tree $T$ we
let $\mathbf{v}(T)$ denote the number of non-leaf vertices in $\mathbf{V}(T)$.

We define 
\begin{equation}
\Lambda_{\subldown}\!(T)=\prod_{x\in \mathbf{V}(T)}^{\subldown}Y(x)=Y(x_{1})\cdots Y(x_{r}), 
\label{lambda}
\end{equation}
where $x_{1}\ldown\ldots \ldown x_{r}$ are the non-leaf vertices
of $T$ (and thus $x_{1}$ is the root).

We may regard an ordered set of $n$ reduced trees $(T_{1},\ldots ,T_{n})$ as
a forest. Letting $T_{x}$ be a branch (see \S3) with root $x$ and
offspring indexed by the colored roots $x_{1},\ldots ,x_{n}$ of $%
T_{1},\ldots ,T_{n},$ 
\[
c_{x}(T_{1},\ldots ,T_{n})=(T_{1},\ldots
,T_{n}) \vartriangleright T_{x}
\]
is the tree obtained by introducing a new colored root $x$ with color $i$,
and edges joining each of the roots $x_{j}$ (with color $i_{j})$ of $T_{j}$
to $x.$ It is evident that with the exception of the unique one layer tree 
$T_{u}^{i}\in \mathbf{R}_{u}^{i}$, every tree $T\in \mathbf{R}_{u}^{i}$ has a
unique representation of the form $T=c_{x}(T_{1},\ldots ,T_{n})$ with $%
T_{j}\in \mathbf{R}_{u}^{i_{j}}$.

\begin{lemma}
Suppose that we are given an ordered $n$-tuple of reduced trees $%
(T_{1},\ldots ,T_{n})$ $(n\geq 2),$ and that the root $x_{j}$ of $T_{j}$ has
color $i_{j}.$ Then 
\begin{equation}
\Lambda_{\subldown} (c_{x}(T_{1},\ldots ,T_{n}))=Y_{i_{1}\ldots
i_{n}}^{i}\Lambda_{\subldown} (T_{n})\ldots \Lambda_{\subldown} (T_{1})
\end{equation}
and $\mathbf{v}(c_{x}(T_{1},\ldots ,T_{n}))=\sum \mathbf{v}(T_{j})+1.$
\end{lemma}

\proof Let us suppose that the non-leaf vertices of $T_{k}$ are given by $%
x_{k,1}\ldown \ldots 
\ldown x_{k,p_{k}},$ and thus $x_{k,1}$ is the root of $T_{k}$. The new root
$x$ is not a leaf in the tree $T=c_{x}(T_{1},\cdots ,T_{n})$ and we have the
ordering 
\[
x\ldown x_{n,1}\ldown x_{n,2} \ldown \ldots \ldown x_{n,p_{n}}\ldown
x_{n-1,1}\ldown
x_{n-1,2}\ldown\ldots \ldown x_{1,p_{1}}. 
\]
It follows that 
\begin{eqnarray*}
\Lambda_{\subldown} (T) &=&Y(x)Y(x_{n,1})Y(x_{n,2})\cdots Y(x_{n-1,1})\cdots
Y(x_{1,p_{1}}) \\
&=&Y_{i_{1}\cdots i_{n}}^{i}\Lambda_{\subldown} (T_{n})\cdots \Lambda_{\subldown}
(T_{1}).
\end{eqnarray*}
The second relation is immediate.\endproof

\begin{theorem}
The left Lagrange antipode is given by 
\begin{equation}
S_{\mathcal{L}}(Y_{u}^{i})=\sum_{T\in \mathbf{RT}_{u}^{i}}(-1)^{\mathbf{v}%
(T)}\Lambda_{\subldown} (T).  \label{depth}
\end{equation}
\end{theorem}

\proof Since $S_{\mathcal{L}}$ satisfies $S_{\mathcal{L}}(1)=1,$ and it is
an antihomomorphism, this relation indeed determines $S_{\mathcal{L}}$ on $%
\mathcal{L}.$ We use the recursive characterization (\ref{inverserec}) for $%
S_{\mathcal{L}}=S_{\mathcal{H}}^{-1}$ . If $u=jk, $ $\mathbf{RT}_{u}^{i}$
contains only the branch $T=T_{jk}^{i}$ and $\Lambda
(T_{jk}^{i})=Y_{jk}^{i}. $ We have
from (\ref{inverserec}),
\[
S_{\mathcal{H}}^{-1}(Y_{jk}^{i})=-Y_{jk}^{i}=(-1)^{\mathbf{v}(T)}\Lambda_{\subldown}
(T), 
\]
which coincides with the right side of (\ref{depth}).

Let us suppose that the formula is true for $\rho (u) = p-1.$ From
(\ref{inverserec}), if $u=u(1)\ldots u(p+1)$ (and thus $\rho (u) \leq p$)
\[
S_{\mathcal{H}}^{-1}(Y_{u}^{i})=-Y_{u}^{i}-\sum {}^{'} (Y_{u}^{i})_{(2)}
S_{\mathcal{H}}^{-1}((Y_{u}^{i})_{(1)}) 
\]
where the prime indicates that we are considering sums of terms
$b_{k}S_{\mathcal{H}}^{-1}(a_{k})$  with $a(k),b(k)\in \mathcal{H}_{+}$. It follows from 
(\ref {expanddelta}) that
\begin{eqnarray*}
S_{\mathcal{H}}^{-1}(Y_{u}^{i}) &=&-Y_{u}^{i}-\sum_{q=2,\pi ,w}^{p-1}Y_{w}^{i}S_{%
\mathcal{H}}^{-1}(Y_{u_{1}}^{w(1)}\cdots Y_{u_{q}}^{w(q)}) \\
&=&-Y_{u}^{i}-\sum_{q=2,\pi
,w}^{p-1}Y_{w}^{i}S_{\mathcal{H}}^{-1}(Y_{u_{q}}^{w(q)})%
\cdots S_{\mathcal{H}}^{-1}(Y_{u_{1}}^{w(1)}) \\
&=&-Y_{u}^{i}-\sum\limits_{q=2,\pi ,w}^{p-1}\sum_{T_{k}\in \mathbf{RT}%
_{u_{k}}^{w(k)}}(-1)^{\sum \mathbf{v}(T_{k})}Y_{w}^{i}\Lambda_{\subldown}
(T_{q})\cdots
\Lambda_{\subldown} (T_{1}) \\
&=&-Y_{u}^{i}+\sum\limits_{q=2,\pi ,w}^{p-1}\sum_{T_{k}\in \mathbf{RT}%
_{u_{k}}^{w(k)}}(-1)^{\mathbf{v}(c_{x}(T_{1},\ldots ,T_{q}))}\Lambda_{\subldown}
(c_{x}(T_{1},\ldots ,T_{q})) \\
&=&-Y_{u}^{i}+\sum_{T\in \mathbf{RT}_{u}^{i}\backslash \{T_{u}^{i}\}}(-1)^{%
\mathbf{v}(T)}\Lambda_{\subldown} (T)
\end{eqnarray*}
where the one layer tree $T_{u}^{i}\in \mathbf{RT}_{u}^{i}$ is not assembled from
non-trivial reduced subtrees, and on the other hand, a reduced tree cannot
have the form $c_{x}(T_{1})$, i.e., $q>1$. Since
$\Lambda_{\subldown} (T_{u}^{i})=Y_{u}^{i}$ and $\mathbf{v}(T_{u}^{i})=1,$  
\[
S_{\mathcal{L}}(Y_{u}^{i})=S_{\mathcal{H}}^{-1}(Y_{u}^{i})=\sum_{T\in
\mathbf{RT}_{u}^{i}}(-1)^{\mathbf{v} (T)}\Lambda _{\subldown}(T). 
\]
\endproof

If $y_{1}\rup y_{2}\ldots \rup y_{r}$ are the non-degenerate vertices in $T$
(and thus $y_{r}$ is the root), we define 
\[
\Lambda _{\subrup }(T)=Y(y_{1})\cdots Y(y_{r}) 
\]
Given $T\in \mathbf{RT}_{u}^{i}$, then $T^{*}\in \mathbf{RT}_{u^{*}}^{i}$, and
this determines a one-to-one correspondence between these two families of
colored ordered trees.

\begin{corollary}\label{redantipodeth} The antipode of the Hopf interval algebra is determined by 
\[
S_{\mathcal{H}}(Y_{u}^{i})=\sum_{T\in \mathbf{RT}_{u}^{i}}(-1)^{\mathbf{v}
(T)}\Lambda _{\subrup }(T) 
\]
\end{corollary}
\proof From Lemma 5, $S_{\mathcal{H}}=\mathbf{s}S_{\mathcal{L}}\mathbf{s,}$ and
since
$\mathbf{ v}(T^{*})=\mathbf{v}(T),$%
\begin{eqnarray*}
S_{\mathcal{H}}(Y_{u}^{i}) &=&\mathbf{s}(S_{\mathcal{L}}(Y_{u^{*}}^{i})) \\
&=&\mathbf{s}\sum_{T\in \mathbf{RT}_{u^{*}}^{i}}(-1)^{\mathbf{v}(T)}\Lambda _{
\subrdown}(T)
\\ &=&\sum_{T\in \mathbf{RT}_{u}^{i}}(-1)^{\mathbf{v}(T)}\mathbf{s}\Lambda
_{\subldown}(T^{*}).
\end{eqnarray*}
If $y_{1}\rup y_{2}\ldots \rup y_{r}$ are the non-degenerate vertices in $T,$
then $y_{1}^{*}\lup \ldots \lup y_{r}^{*}$ in $T^{*},$ or changing notation,
$y_{r}^{*}\ldown \ldots \ldown y_{1}^{*}$. It follows that 
\[
\mathbf{s}\Lambda _{\subldown}(T^{*})=\mathbf{%
s}Y(y_{r}^{*})\cdots Y(y_{1}^{*})=Y(y_{1})\cdots Y(y_{r})=\Lambda
_{\subrup}(T), 
\]
and the desired formula follows.
\endproof

If $z_{1}\rdown z_{2}\rdown \ldots \rdown z_{r}$ are the non-degenerate
vertices of $T$ we define $\Lambda _{\subrdown }(T)=Y(z_{1})\cdots Y(z_{r})$.

\begin{corollary} The antipode of the right Lagrange algebra $\mathcal{R}$ is
given by
\[
S_{\mathcal{R}}(Y_{u}^{i})=\sum_{T\in \mathbf{RT}_{u}^{i}}(-1)^{\mathbf{v}%
(T)}\Lambda _{\subrdown }(T).
\]
\end{corollary}
\proof From the definition of the reflected algebra $\mathcal{R}$,
$S_{\mathcal{R}}=\mathbf{t}S_{\mathcal{H}}\mathbf{t}$, and thus
\begin{eqnarray*}
S_{\mathcal{R}}(Y_{u}^{i}) &=&\mathbf{t}(S_{\mathcal{H}}(Y_{u}^{i})) \\
&=&\mathbf{t}(\sum_{T\in \mathbf{RT}_{u}^{i}}(-1)^{\mathbf{v}(T)}\Lambda _{\subrup
}(T) )\\
&=&\sum_{T\in \mathbf{RT}_{u}^{i}}(-1)^{\mathbf{v}(T)}\mathbf{t}(\Lambda _{\subrup
}(T))
\end{eqnarray*}
If $z_{1}\rdown z_{2}\rdown \ldots \rdown z_{r}$ are the non-degenerate
vertices of $T$ then $z_{r}\rup \cdots \rup z_{1}$ implies that 
\begin{eqnarray*}
\mathbf{t(}\Lambda _{\subldown }(T)) &=&\mathbf{t(}Y(z_{r})\cdots Y(z_{1})) \\
&=&Y(z_{1})\cdots Y(z_{r})\\
&=&\Lambda _{\subrdown }(T)
\end{eqnarray*}
and we are done.
\endproof

\section{Cancellations in the the breadth first antipodal formula}

Given $i\in \langle N \rangle$ and $v\in \langle N \rangle^{*}$, we indexed
the summands of $S_{\mathcal{H}}(Y^{i}_{u})$ by trees in $\mathbf{LT}_{u}^{i}$,
the layered trees with root and leaf colorings $i$ and $u$ (see
(\ref{omegaantipode})). We wish to show that it suffices to use ``order
reduced'' simple trees.

It would be tempting to attempt to use (\ref{depth}) to obtain a formula
with reduced trees by simply contracting the edges in layered trees, keeping track
of the resulting ``multiplicities''. In the commutative
situation considered by Haiman and Schmitt, one has that $\Omega (\rho (T))=\Omega (T)$ and thus in
the formula for the antipode one may collect all the terms with the same reduced tree into a
multiple of $\Omega (\rho (T))$. The coefficient is a sum of positive and
negative 1's, and using a combinatorial argument they showed that all of the
non-reduced trees cancel.

In our situation, arbitrary contractions can disturb the $x\ll y$ ordering
on the non-degenerate leaves, and thus one need not have that $\Omega (\rho
(T))=\Omega (T)$. This can be seen in the third tree of the diagram below,
which was disordered by an ``improper'' contraction. One must therefore use
only contractions which are $\ll $ \emph{order preserving}. In our reduction
we will also modify the contraction so that the tree remains layered.

Let us suppose that $T$ is a layered tree and that $x$ is a non-degenerate
vertex in $T$. We say that $T$ is \emph{order contractible at} $x$ if 
\begin{list}{\alph{bean})}{\usecounter{bean}}
\item its parent $x^{\prime}$ is unary, i.e., $x$ has no siblings,
\item there does not exist a non-degenerate vertex to the \emph{right} of $x$
\item there does not exist a non-degenerate vertex to the \emph{left} of $x^{\prime}$.
\end{list}
If $x$ is a vertex in the $k$-th row which satisfies these conditons, the 
\emph{order contraction} ${\kappa}(T)={\kappa}_{x}(T)$ is the layered tree
obtained in the following manner: 
\begin{list}{\arabic{bean})}{\usecounter{bean}}
\item move $x$ to the position of its parent in the $(k-1)$-st row, 
\item attach each offspring $y$ of $x$
by a single line to a unary vertex $x^{\prime}$ in the $k$-th row,
\item leave all other vertices and edges alone,
\item if there are no other non-degenerate vertices in the $k$-th row, delete it.
\end{list}
Conditions a)-c) guarantee that the $\ll $ ordering on the non-degenerate
vertices is preserved. Thus the contraction on the non-degenerate vertex $%
x_{2}$ in the first tree below is allowed. On the other hand contracting on
the vertex $x_{1}$ would transpose the $\ll $ ordering for the two
non-degenerate vertices $x_{1}$ and $x_{2}$. \Treestyle{\vdist{16pt} %
\minsep{12pt}} 
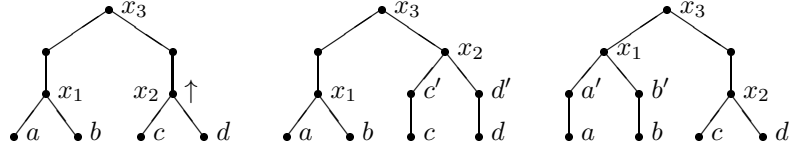
\[
\begin{Tree} 
\node{\external\type{dot}\rght{$a$}}                                   
\node{\external\type{dot}\rght{$b$}}                                             
\node{\type{dot}\rght{$x_{1}$}}
\node{\unary\type{dot}}
\node{\external\type{dot}\rght{$c$}}
\node{\external\type{dot}\rght{$d$}} 
\node{\type{dot}\lft{$x_{2}$}\rght{$\mathbf{\uparrow} $}}
\node{\unary\type{dot}}                      
\node{\type{dot}\rght{$x_{3}$}}                                                                                                                          
\end{Tree}
\hskip\leftdist\box\TeXTree\hskip\rightdist\qquad 
\begin{Tree}
\node{\external\type{dot}\rght{$a$}}                                   
\node{\external\type{dot}\rght{$b$}}                                             
\node{\type{dot}\rght{$x_{1}$}}
\node{\unary\type{dot}}
\node{\external\type{dot}\rght{$c$}}
\node{\unary\type{dot}\rght{$c^{\prime}$}}
\node{\external\type{dot}\rght{$d$}}
\node{\unary\type{dot}\rght{$d^{\prime}$}} 
\node{\type{dot}\rght{$x_{2}$}}                    
\node{\type{dot}\rght{$x_{3}$}}                                                                                                   
\end{Tree}
\hskip\leftdist\box\TeXTree\hskip\rightdist\qquad 
\begin{Tree}
\node{\external\type{dot}\rght{$a$}}
\node{\unary\type{dot}\rght{$a^{\prime}$}}
\node{\external\type{dot}\rght{$b$}}
\node{\unary\type{dot}\rght{$b^{\prime}$}}
\node{\type{dot}\rght{$x_{1}$}}
\node{\external\type{dot}\rght{$c$}}
\node{\external\type{dot}\rght{$d$}}
\node{\type{dot}\rght{$x_{2}$}}
\node{\unary\type{dot}}
\node{\type{dot}\bnth{
test}\rght{$x_{3}$}} 
\end{Tree}                                                                                                                        
\hskip\leftdist\box\TeXTree\hskip\rightdist\qquad 
\]

Since we will not consider general contractions in this section, we will
simply use the terms \emph{contractible} and \emph{contractions} for the
corresponding order preserving notions.

Given a tree $T\in \mathbf{LT}_{v}^{j}$ we define the \emph{\ canonical
expansion} $\Phi (T)\in \mathbf{LT}_{v}^{j}$ as follows. If $T$ is simple we
let $\Phi (T)=T.$ If $T$ is not simple, let $x_{r}$ be the first non-simple
non-degenerate vertex in the $\ll $ ordering, and let us suppose that it is
on the $k$-th level. We introduce a new level $L$ between the $k$-th and $%
(k+1)$-st levels in the following manner. 
\begin{list}{\alph{bean})}{\usecounter{bean}}
\item We move $x_{r}$ down to the level $L$ and we connect it to a new unary vertex
$x_{r}^{\prime}$  on the $k$-th level, and to the offspring of $x_{r}$ on the $(k+1)$-st
level,
\item If $y$ is a vertex on the $(k+1)$-st that is not an offspring of $x_{r}$, we
connect it by a single edge to a new unary vertex $y^{\prime}$ on the level $L$,
which we then connect to the parent of $y$ on the $k$-th level.  
\end{list}
We define $\Phi (T)$ to be the new tree. \Treestyle{\vdist{20pt}\minsep{8pt}}
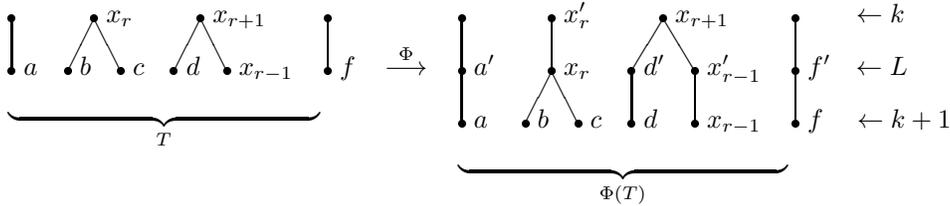
\[
\underbrace{%
\begin{Tree} 
\node{\external\type{dot}\rght{$a$}}
\node{\unary\type{dot}}                                                                                                                                                                        
\end{Tree}                                                                                  
\hskip\leftdist\box\TeXTree\hskip\rightdist\quad 
\begin{Tree}
\node{\external\type{dot}\rght{$b$}}
\node{\external\type{dot}\rght{$c$}}
\node{\type{dot}\rght{$x_{r}$}}                                                                                                                         
\end{Tree}
\hskip\leftdist\box\TeXTree\hskip\rightdist\quad \extended
\begin{Tree}
\node{\external\type{dot}\rght{$d$}}
\node{\external\type{dot}\rght{$x_{r-1}$}}
\node{\type{dot}\rght{$x_{r+1}$}}
\end{Tree}
\hskip\leftdist\box\TeXTree\hskip\rightdist\quad 
\begin{Tree} 
\node{\external\type{dot}\rght{$f\,\,\,\,\,\stackrel{\Phi }{\longrightarrow }
\,\,\,\,$}}
\node{\unary\type{dot}}                                                                                                                                                
\end{Tree}
}_{T} \hskip\leftdist\box\TeXTree\hskip\rightdist\ \underbrace{%
\begin{Tree} 
\node{\external\type{dot}\rght{$a$}}
\node{\unary\type{dot}\rght{$a^{\prime}$}}
\node{\unary\type{dot}}                                                                                                                                                                         
\end{Tree}                                                                                  
\hskip\leftdist\box\TeXTree\hskip\rightdist\quad 
\begin{Tree}
\node{\external\type{dot}\rght{$b$}}                                   
\node{\external\type{dot}\rght{$c$}}                                             
\node{\type{dot}\rght{$x_{r}$}}
\node{\unary\type{dot}\rght{$x_{r}^{\prime}$}}                                                                                                 
\end{Tree}
\hskip\leftdist\box\TeXTree\hskip\rightdist\quad \extended
\begin{Tree}
\node{\external\type{dot}\rght{$d$}}
\node{\unary\type{dot}\rght{$d^{\prime}$}}
\node{\external\type{dot}\rght{$x_{r-1}$}}
\node{\unary\type{dot}\rght{$x_{r-1}^{\prime}$}}
\node{\type{dot}\rght{$x_{r+1}$}}
\end{Tree}
\hskip\leftdist\box\TeXTree\hskip\rightdist\quad 
\begin{Tree} 
\node{\external\type{dot}\rght{$f\,\,\,\,\,\,\leftarrow k+1$}}
\node{\unary\type{dot}\rght{$f^{\prime}\,\,\,\,\leftarrow L$}}
\node{\unary\type{dot}\rght{$\,\,\,\,\,\,\,\,\,\,\,\leftarrow k$}}                     
\end{Tree}}_{\Phi (T)} \hskip\leftdist\box\TeXTree\hskip\rightdist\ 
\]
It should be noted that since there are no non-degenerate vertices to the
left of $x_{r}$, this operation will not affect the $\ll $ ordering on the
non-degenerate vertices. It also preserves the orderings $\!\lup\!$ and $%
\!\rup\!$ on the non-degenerate vertices, as is evident from the above
diagram.

We have that $x_{r}$ is a contractible vertex in $\Phi(T)$ because all the
other vertices on the new level are unary, and there are no non-degenerate
vertices to the left of $x_{r}^{\prime}$. If one contracts on this vertex,
the new level will contain only unary vertices, and thus will itself be
deleted (see the primed row in the right tree below). In this manner we see
that if we contract $\Phi(T)$ at the vertex $x_{r}$, we recover $T$. This is
illustrated in the following diagram, in which $e=x_{r-1}$ and $%
e^{\prime}=x_{r-1}^{\prime}$. 
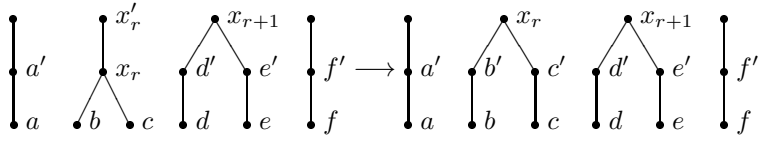
\[
\begin{Tree} 
\node{\external\type{dot}\rght{$a$}}
\node{\unary\type{dot}\rght{$a^{\prime}$}}
\node{\unary\type{dot}}                                                                                                                                                                         
\end{Tree}                                                                                  
\hskip\leftdist\box\TeXTree\hskip\rightdist\quad 
\begin{Tree}
 \node{\external\type{dot}\rght{$b$}}                                   
\node{\external\type{dot}\rght{$c$}}                                             
\node{\type{dot}\rght{$x_{r}$}}
\node{\unary\type{dot}\rght{$x_{r}^{\prime}$}}                                                                                                  
\end{Tree}
\hskip\leftdist\box\TeXTree\hskip\rightdist\quad \extended
\begin{Tree}
\node{\external\type{dot}\rght{$d$}}
\node{\unary\type{dot}\rght{$d^{\prime}$}}
\node{\external\type{dot}\rght{$e$}}
\node{\unary\type{dot}\rght{$e^{\prime}$}}
\node{\type{dot}\rght{$x_{r+1}$}}
\end{Tree}
\hskip\leftdist\box\TeXTree\hskip\rightdist\quad 
\begin{Tree} 
\node{\external\type{dot}\rght{$f$}}
\node{\unary\type{dot}\rght{$f^{\prime}\longrightarrow$}}
\node{\unary\type{dot}}                                                                                                                                                
\end{Tree}                                                                                  
\hskip\leftdist\box\TeXTree\hskip\rightdist\ 
\begin{Tree} 
\node{\external\type{dot}\rght{$a$}}
\node{\unary\type{dot}\rght{$a^{\prime}$}}
\node{\unary\type{dot}}                                                                                                                                                                         
\end{Tree}                                                                                  
\hskip\leftdist\box\TeXTree\hskip\rightdist\quad 
\begin{Tree}
 \node{\external\type{dot}\rght{$b$}}
\node{\unary\type{dot}\rght{$b^{\prime}$}}
\node{\external\type{dot}\rght{$c$}}
\node{\unary\type{dot}\rght{$c^{\prime}$}}
\node{\type{dot}\rght{$x_{r}$}}                                                                                                  
\end{Tree}
\hskip\leftdist\box\TeXTree\hskip\rightdist\quad \extended
\begin{Tree}
\node{\external\type{dot}\rght{$d$}}
\node{\unary\type{dot}\rght{$d^{\prime}$}}
\node{\external\type{dot}\rght{$e$}}
\node{\unary\type{dot}\rght{$e^{\prime}$}}
\node{\type{dot}\rght{$x_{r+1}$}}
\end{Tree}
\hskip\leftdist\box\TeXTree\hskip\rightdist\quad 
\begin{Tree} 
\node{\external\type{dot}\rght{$f$}}
\node{\unary\type{dot}\rght{$f^{\prime}$}}
\node{\unary\type{dot}}                                                                                                                                                
\end{Tree}                                                                                  
\hskip\leftdist\box\TeXTree\hskip\rightdist\ 
\]

Let us suppose that $x_{1}\ll \ldots \ll x_{p}$ are the non-degenerate
vertices of a tree $T.$ Turning to the breadth first ordering, there is a
unique sequence of indices $n_{1}<\ldots <n_{q}$ with 
\begin{equation}
\ldots x_{n_{1}}\Rsh \!\!\!\!\Rsh x_{n_{1}+1}\rup x_{n_{1}+2}\rup \ldots
\rup x_{n_{2}}\Rsh \!\!\!\!\Rsh x_{n_{2}+1}\ldots  \label{string}
\end{equation}
where we let $x \rup y$ if $x$ is a descendant of $y$. We call a maximal
sequence of the form $x_{n_{h-1}+1}\rup x_{n_{h-1}+2}\rup \ldots \rup
x_{n_{h}}$ an \emph{irreducible string} and we say that $x_{n_{h}}$ is its 
\emph{right end}.

\begin{lemma}
Suppose that $T$ is an arbitrary tree in $\mathbf{LT}_{u}^{i}$ with its
non-degenerate vertices $x_{1}\ll \ldots \ll x_{p}$ satisfying (\ref
{string}).

(i) Any contractible vertex in $T$ is a right end of an irreducible
string. 

(ii) If $T=E$ is simple then all of its right ends without siblings
are contractible. 

(iii) If one has $y\ll x$ in $T$ and both $y$ and $x$ are
contractible, then after a contraction at $x$, $y$ will still be
contractible. 
\end{lemma}

\proof (i) Let us suppose that $x=x_{j}$ is a contractible non-degenerate
vertex in $T$ on level $k.$ Then its parent $x_{j}^{\prime }$ is unary and
there are no non-degenerate vertices to the left of $x_{j}^{\prime }$. Since
every level is assumed to have a non-degenerate vertex, $x_{j+1}$ must lie
on the $(k-1) $-st row of $E$ to the right of $x_{j}^{\prime }$. It follows
that $x_{j}\Rsh \!\!\!\!\Rsh x_{j+1},$ and thus $x_{j}=x_{n_{h}}$ for some $%
h.$

(ii) Let us suppose that $x_{n_{h}}$ is a non-degenerate vertex on the $k$%
-th level. There are no non-degenerate vertices to the right of $x_{n_{h}}$
on its level since the level is simple. We have that $x_{n_{h}+1}$ must lie
on the previous level. Since $x_{n_{h}}\Rsh \!\!\!\!\Rsh x_{n_{h}+1}$, the
parent $x_{n_{h}}^{\prime }$ of $x_{n_{h}}$ lies to the left of $%
x_{n_{h}+1}, $ and $x_{n_{h}}$ is (order) contractible.

(iii) If $x$ is simple, then the contraction at $x$ will simply raise the
level of each vertex $y$ with $y\ll x$. If $x$ is not simple, then $x$ is
the only non-degenerate vertex that is affected. Since $y$ is assumed
contractible, there will not be any vertices to the left of it on its level.
On the other hand if $y$ is on level $k+1$, then by the same assumption, $x$
must lie to the right of the parent $y^{\prime}$. This will still be be the
case when one contracts $x$ to a higher level.\endproof

For each simple layered tree $E$ we let $\mathbf{T}_{E}$ be all the trees $%
T\in \mathbf{LT}_{v}^{j}$ with $E=\Phi ^{n}(T)$ for some $n.$ It is evident
that if $T$ has $n$ vertices and $k$ levels, then $E=\Phi ^{n-k}(T)$ is a
simple tree, and reversing the expansions as above, $T$ can be obtained by a
particular sequence of contractions of $E.$ More precisely, let $y_{1}\ll
\ldots \ll y_{q}$ be the right vertices of $T$ (or equivalently of $E$). We
have that there is a subsequence $y_{m_{1}}\ll \ldots \ll y_{m_{p}}$ with 
\[
T=\kappa _{y_{m_{1}}}\ldots \kappa _{y_{m_{p}}}(E). 
\]
Conversely given any such sequence, the subsequent right vertices remain
contractible as one proceeds, and we get a corresponding tree $T$. The tree $%
T$ uniquely determines the sequence $y_{m_{1}}\ll \ldots \ll y_{m_{p}}$
since the latter are by definition the non-simple non-degenerate vertices of 
$T$ in their given $\ll $ order.

We say that a layered tree is \emph{order reduced} if it does not have any
non-trivial ordered contractions, and we let $\mathbf{OST}_{v}^{j}$ 
be the set of all ordered reduced simple trees.

\begin{theorem}\label{ordantipodeth}
The antipode in $\mathcal{H}$ is given by 
\[
S_\mathcal{H}(Y_{v}^{j})=\,\,\sum_{E\in
\mathbf{OST}_{v}^{j}}(-1)^{\ell (E)}\Omega (E) 
\]
where $\ell (E)$ is the number of layers in $E.$
\end{theorem}

\proof It is evident that 
\[
\mathbf{LT}_{v}^{j}=\sqcup \{\mathbf{T}_{E}:E\in \mathbf{ST}_{v}^{j}\} 
\]
and that for any $T\in \mathbf{T}_{E}$ we have that $\Omega (T)=\Omega (E).$
Thus it suffices to show that if $E$ has contractions, then 
\[
\sum_{T\in \mathbf{T}_{E}}(-1)^{\ell (T)}\Omega (T)=\Omega (E)\sum_{T\in \mathbf{%
T}_{E}}(-1)^{\ell (T)}=0. 
\]
From our earlier discussion, $\mathbf{T}_{E}$ is in one-to-one correspondence
with the sequences $y_{m_{1}}\ll \ldots \ll y_{m_{p}}$ drawn from the $q$
right vertices in $E,$ or equivalently subsets drawn from ${1,\ldots ,q}$.
If the simple tree $E$ has $n$ non-degenerate vertices and thus $n$ levels,
the tree $T(m_{1},\ldots ,m_{p})$ has $n-p$ levels. There will be $\binom{q}{%
p}$ such sequence and thus 
\[
\sum_{T\in \mathcal{\mathcal{T}}_{E}}(-1)^{\ell (T)}=(-1)^{n-q}\sum_{p}%
\binom{q}{p}(-1)^{q-p}=(-1)^{n-q}(1-1)^{q}=0, 
\]
and we have proved the desired result.\endproof

\section{Breadth first and depth first duality}

There is a natural one-to-one correspondence between the reduced trees $%
\mathbf{RT}_{u}^{j}$ and the order-reduced simple layered trees
$\mathbf{OST}_{u}^{j}.$ Since the colorings as well as the set of leaves is
unaffected by the operations, we will use the notations
$\mathbf{T}$, $\mathbf{RT}$ and $\mathbf{OST}$ for the proper trees, the
reduced trees, and the order-reduced simple layered $N$-colored trees, respectively,
with a given
colored root and a given set of colored leaves. We let $\rho:\mathbf{OST}\to
\mathbf{RT}$ be the contraction mapping described in \S3.

We say that a vertex $x$ in a proper tree $T$ is \emph{weakly contractible} if it has no siblings,
i.e., its parent $x^{\prime}$ is unary. We say that $T$ is \emph{reduced} let
$\rho:\mathbf{T}\to
\mathbf{RT}$ be the contraction mapping. In this situation we say that any vertex without siblings
(this applies to leaves as well) is said to  be contractible. If
$\rho(T^{\prime})=T$, we may identify the non-leaf vertices of $T$ with the
non-degenerate vertices in $T^{\prime}$. Although we have seen that $\rho$
can disrupt the ordering $\ll$, it is evident that the depth first ordering 
$\rup$ is unaffected. 

\begin{theorem} Let $T$ be a layered tree.
\begin{itemize}
\item[(a)] The breadth first ordering $\ll $
coincides with the depth first ordering $\rup $ on the non-degenerate
vertices of $T$ if and only if $ T\in \mathbf{OST}$.
\item[(b)] The contraction mapping $\rho$ is a bijection of $\mathbf{OST}$ onto
$\mathbf{RT}$. 
\item[(c)] If $T^{\prime}$ is the unique tree in $\mathbf{OST}$ with 
$\rho(T^{\prime})=T$, then $\Omega (T^{\prime})=\Lambda_{\subrup} (T)$ and $\ell(T^{\prime})=\mathbf{v}(T) $.
\end{itemize}
\end{theorem}
\proof (a) Since these orderings are inverse to each other on the vertices in a
given row, there will be only one non-degenerate vertex in that row, i.e.,
$T$ is simple. Thus it suffices to prove the equivalence for simple trees.

Let us suppose that the non-degenerate vertices of $T$ are given by $x_{1}\ll
\ldots \ll x_{n}$. Since $T$ is simple, these vertices will lie on
the successively higher levels of $T$. A non-degenerate non-root vertex
$x_{j}$ is non-order contractible if and only if either its parent
$x_{j}^{\prime}$ is non-degenerate, and thus coincides with $x_{j+1}$, or
$x_{j+1}$ lies to the left of $x_{j}^{\prime}$. These are precisely the
conditions that $x_{j}\rup x_{j+1}$.  

(b)We define $\Psi :\mathbf{RT}\rightarrow \mathbf{OST}$ by sequentially inserting
singular edges. Let us suppose that $y_{1}\rup\ldots\rup y_{n}$ are the
non-degenerate vertices in $T$ (see the diagram below). Counting down from $n$, let us suppose that
$i$ is the last index for which
$y_{i-1}\gg y_{i}$ (this would be $y_{3}\gg y_{4}$ in the diagram). Then we lower
$y_{i-1}$ and the subtree from which it is a root by inserting the minimal number of unary vertices
and singular edges so that the resulting vertex $y_{i-1}^{\prime}$ satisfes $y_{i-1}^{\prime}\ll
y_{i}$. If the resulting tree is not proper, we eliminate the redundant rows having only
unary vertices. Relabelling, we may
assume that
$y_{i-1}\ll y_{i}$. By using singular edges to push down leaves to the bottom level, we
obtain a layered tree. Furthermore by
eliminating ``redundant'' levels, we may  assume that the layered tree is proper. If $y_{j-1}\gg
y_{j}$ is the next occurence of the relation
$\gg $ in the sequence, we again ``lower'' the vertex $y_{j-1}$, and we proceed as before. This will
not affect the fact that
$y_{i-1}\ll y_{i}$ and lowering does not change the $\rup$ relation. It
should be pointed out that ``lowering'' will in general disrupt the ordering in general, but we are
only concerned with successive terms in the given sequence. After ``correcting'' all of the reverse
orderings, we obtain 
$T^{\prime}=\Psi(T)$. It is evident that
$T^{\prime}$ satisfies the conditions in (a), and thus lies in $\mathbf{OST}$.
 
The mapping $\Psi$ is illustrated in the following diagram, in which one initially has
$y_{1}\rup y_{2}\rup y_{3}\rup y_{4}\rup y_{5}$ in both trees, $y_{1}\ll y_{2}
\gg y_{3}\gg y_{4}\ll y_{5}$ in $T$, and $y_{1}\ll y_{2}\ll y_{3} \ll y_{4} \ll
y_{5}$ in $T^{\prime}$,
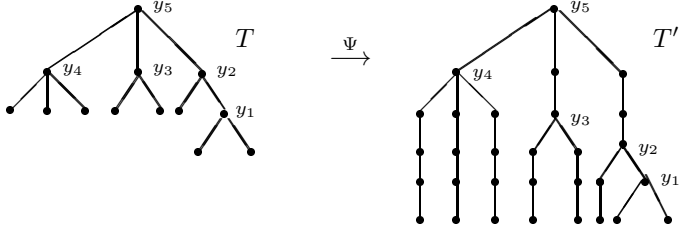
\vspace{-.9in}
\[\hspace{-1in}
\begin{picture}(4,4)
\linethickness{.2mm}
\put(2.8,1){$T$}
\put(.3,.67){\line(1,-1){.5}}
\put(.32,.67){\line(1,-1){.5}}
\put(1.2,.16){\line(2,3){.33}}
\put(1.22,.16){\line(2,3){.33}}
\put(1.8,.18){\line(-2,3){.33}}
\put(1.82,.18){\line(-2,3){.33}}
\put(.3,.15){\line(0,1){.5}}
\put(1.55,1.48){\line(1,-1){.8}}
\put(1.57,1.48){\line(1,-1){.8}}
\put(2.05,.17){\line(2,3){.3}}
\put(2.07,.17){\line(2,3){.3}}
\put(2.65,.17){\line(-2,3){.3}}
\put(2.665,.17){\line(-2,3){.3}}
\put(.3,.66){\circle*{.1}}
\put(.5,.66){$\scriptstyle{y_{4}}$}%
\put(1.5,.66){\circle*{.1}}
\put(1.7,.66){$\scriptstyle{y_{3}}$}%
\put(2.35,.64){\circle*{.1}}
\put(2.55,.64){$\scriptstyle{y_{2}}$}%
\put(1.2,.15){\circle*{.1}}
\put(1.8,.15){\circle*{.1}}
\put(2.05,.15){\circle*{.1}}
\put(2.3,-.4){\circle*{.1}}
\put(3,-.4){\circle*{.1}}
\put(2.3,-.4){\line(2,3){.3}}
\put(2.32,-.4){\line(2,3){.3}}
\put(3,-.4){\line(-2,3){.3}}
\put(3.02,-.4){\line(-2,3){.3}}
\put(2.65,.1){\circle*{.1}}
\put(-.2,.15){\circle*{.1}}
\put(2.8,.1){$\scriptstyle{y_{1}}$}
\put(.3,.15){\circle*{.1}}
\put(.8,.15){\circle*{.1}}
\put(-.2,.15){\line(1,1){.5}}
\put(.3,.7){\line(3,2){1.2}}
\put(.31,.7){\line(3,2){1.2}}
\put(1.5,.7){\line(0,1){.75}}
\put(1.5,1.5){\circle*{.1}}
\put(1.7,1.5){$\scriptstyle{y_{5}}$}%
\end{picture}
\hskip\leftdist\raisebox{.3in}{$\stackrel{\Psi}{\longrightarrow}$\!\!\!\!}\hskip\rightdist\qquad
\begin{picture}(3,4)
\put(2.8,1){$T^{\prime}$}
\linethickness{.2mm}
\put(.07,.3){\line(1,1){.1}}
\put(.7,.7){\line(1,-1){.1}}
\put(.18,.68){\line(1,-1){.53}}
\put(1.2,-.4){\line(2,3){.3}}
\put(1.21,-.4){\line(2,3){.3}}
\put(1.8,-.4){\line(-2,3){.3}}
\put(1.82,-.4){\line(-2,3){.3}}
\put(.2,-1.3){\line(0,1){2}}
\put(.7,-1.3){\line(0,1){1.4}}
\put(-.3,-1.3){\line(0,1){1.4}}
\put(1.2,-1.3){\line(0,1){.93}}
\put(1.8,-1.3){\line(0,1){.93}}
\put(2.6,-.4){$\scriptstyle{y_{2}}$}%
\put(1.7,0){$\scriptstyle{y_{3}}$}%
\put(.4,.6){$\scriptstyle{y_{4}}$}%
\put(2.1,-1.3){\line(0,1){.5}}
\put(1.55,1.48){\line(1,-1){.85}}
\put(1.57,1.48){\line(1,-1){.85}}
\put(2.1,-.75){\line(2,3){.3}}
\put(2.11,-.75){\line(2,3){.3}}
\put(2.3,-1.3){\line(2,3){.35}}
\put(2.3,-1.3){\line(2,3){.35}}
\put(2.32,-1.3){\circle*{.1}}
\put(3,-1.3){\circle*{.1}}
\put(2.68,-.75){\line(-2,3){.3}}
\put(2.7,-.75){\line(-2,3){.3}}
\put(2.69,-.75){\line(-2,3){.3}}
\put(2.4,-.35){\line(0,1){1}}
\put(.18,.66){\circle*{.1}}
\put(1.5,.66){\circle*{.1}}
\put(1.5,.1){\circle*{.1}}
\put(2.4,.64){\circle*{.1}}
\put(1.8,-.4){\circle*{.1}}
\put(1.2,-.4){\circle*{.1}}
\put(2.4,.1){\circle*{.1}}
\put(2.4,-.3){\circle*{.1}}
%\put(2.4,-.3){\circle*{.1}}%
\put(2.1,-.8){\circle*{.1}}
\put(1.8,-.8){\circle*{.1}}
\put(1.8,-1.3){\circle*{.1}}
\put(1.2,-.8){\circle*{.1}}
\put(1.2,-1.3){\circle*{.1}}
\put(.7,-.8){\circle*{.1}}
\put(.7,-1.3){\circle*{.1}}
\put(.18,-.8){\circle*{.1}}
\put(.18,-1.3){\circle*{.1}}
\put(-.3,-.8){\circle*{.1}}
\put(-.3,-1.3){\circle*{.1}}
\put(2.1,-.8){\circle*{.1}}
\put(2.1,-1.3){\circle*{.1}}
\put(2.7,-.8){\circle*{.1}}
\put(2.9,-.8){$\scriptstyle{y_{1}}$}%
\put(2.7,-.7){\line(1,-2){.3}}
\put(2.71,-.7){\line(1,-2){.3}}
\put(.18,.1){\circle*{.1}}
\put(.18,-.4){\circle*{.1}}
\put(.7,.1){\circle*{.1}}
\put(.7,-.4){\circle*{.1}}
\put(-.3,.1){\circle*{.1}}
\put(-.3,-.4){\circle*{.1}}
\put(-.3,.15){\line(1,1){.5}}
\put(.19,.68){\line(3,2){1.25}}
\put(.21,.68){\line(3,2){1.25}}
\put(1.5,1.45){\line(0,-1){1.4}}
\put(1.49,1.5){\circle*{.1}}
\put(1.7,1.5){$\scriptstyle{y_{5}}$}
\end{picture} 
\]

\vspace{.5in}
\noindent It is evident that $\rho(\Psi (T))=T$, and thus $\rho$ is surjective.  

Let us suppose that $T^{\prime}$ is any layered tree with 
$\rho(T^{\prime})=T\in \mathbf{RT}$. As can be see from the above diagrams (this does
not depend upon $T^{\prime}$ being order reduced),
$T^{\prime}$ may be constructed by inserting $d(v)$ singular edges between each vertex $v$
and its parent $v^{\prime}$ i.e.,
$T^{\prime}$ is characterized by the function $d:\mathbf{v}(T)\to \Bbb{N}\cup {0}$. For our purposes
it is not necessary to characterize the functions $d$ that arise in this manner. 

Let us suppose $T^{\prime}\in \mathbf{OST}$ and that $x_{1} \rup \ldots \rup
x_{n}$ are the non-degenerate (i.e., non-leaf) vertices of the reduced tree $T$. Then $x_{1} \rup
\ldots
\rup x_{n}$ are the non-degenerate vertices of
$T^{\prime}$, and from (a), $x_{1} \ll \ldots \ll x_{n}$. Since $T^{\prime}$ is
simple, each is on a different level, i.e., they reside on the successive levels
upwards. 

Let us suppose that
$x=x_{i}$ is a nonleaf vertex in $T$ with parent $x^{\prime}=x_{j}$ (recall that all the non-leaf
vertices in $T$ are non-degenerate). Since $x_{i}$ is an offspring of $x_{j}$, $x_{i}\rup x_{j}$ and
thus $i<j$. It follows that in $T^{\prime}$ there are $j-i$ levels between them. It
hence $d(x)=j-i-1$. On the other hand if $x$ is a leaf, then $d(x)$ is just the number of levels
between $x$ and the bottom level (there are precisely $n$ levels in $T^{\prime}$). Thus if $x$ is a
leaf on the $k$-th level of $T$, then $d(x)=n-k$. Thus the function
$d$ and the tree 
$T^{\prime}$ are uniquely determined by the
orderings $\rup$ and $\ll$ on the non-degenerate vertices in $T$. 

(c) Letting $x_{1} \rup \ldots \rup x_{n}$ be the non-leaf vertices in $T$, we have that 
$x_{1} \ll \ldots \ll x_{n}$ are the non-degenerate vertices in $T^{\prime}$, and thus 
\[
\Omega(T^{\prime})=Y(x_{1})\ldots Y(x_{n})=\Lambda_{\subrup}(T).
\] 
The number $\ell(T^{\prime})$ of layers in $T^{\prime}$ is equal to the number
of non-degenerate vertices in $T^{\prime}$, and thus the number $\mathbf{v}(T)$ of non-degenerate
vertices in $T$. 
\endproof 

If one uses the above result, Corollary \ref{redantipodeth} is an immediate consequence of Theorem 
\ref{ordantipodeth}. In this sense, the Haiman-Schmitt approach to the reduced formula for the antipode
in the Fa\`{a} di bruno algebra can be adapted to the Hopf algebra $\mathcal{H}$ of ordered 
colored partitions. 

\section{Formal Power Series}

Let us suppose that we are given a non-commutative unital algebra $A$ and
non-commuting variables $z_{1},\ldots ,z_{N}.$ Given a word $w=w(1)\cdots
w(p)\in [N]^{*},$ we let $z_{w}=z_{w(1)}\cdots z_{w(p)},$ and $z_{e}=1.$ A
multiple non-commutative power series with $N$ non-commuting variables $%
z_{1},\ldots ,z_{N}$ and non-commuting constants has the form $%
F(z)=(F^{1}(z_{1},\ldots z_{n}),\ldots ,F^{N}(z_{1},\ldots ,z_{N}))$ where 
\[
F^{j}(z)=F^{j}(z_{1},\ldots ,z_{N})=\sum f_{w}^{j}z_{w}, 
\]
and the ``constants'' $f_{w}^{j}$ lie in $A.$ We assume that variables
commute with constants. The latter enables us to multiply power series since
in particular,

\[
(az_{v})(bz_{w})=abz_{vw}. 
\]

Let us begin by computing the effect of \emph{substitution }on power series.
We do not use the term ``composition'' since there does not seem to be a
meaningful interpretation along those lines. Given a single power series of $%
N$ variables 
\[
F(z)=f_{e}+\sum f_{j}z_{j}+\sum f_{jk}z_{j}z_{k}+\cdots 
\]
and an $N$-tuple of power series without constant terms 
\[
G^{i}(z)=\sum g_{j}^{i}z_{j}+\sum f_{jk}z_{j}z_{k}+\cdots 
\]
we may substitute $G^{j}(z)$ for $z_{j}$ in the expression or $F.$ We will
denote the resulting power series by $H(z)=(F\circ G)(z).$

To see that substitution is non-associative, one need only consider the one
variable expressions $(F\circ G)\circ H$ and $F\circ (G\circ H),$ where $%
H(x)=az,G(x)=bz$, and $F(z)=z^{2},$ where $a$ and $b$ are not assumed to
commute.

We wish to compute the coefficients of substituted series. For this purpose
we consider the coefficient $h_{w}$ of $z_{w}.$ A typical summand of $h_{w}$
where $w=w(1)\ldots w(p)$ is obtained by taking a colored ordered interval
partition $((C_{1},\ldots ,C_{q}),j_{1}\cdots j_{q})$ of $((1,\ldots ,p),w).$
If $C_{1}=(1,\ldots ,p_{1}) $ then 
\[
g_{w(1)\ldots w(p_{1})}^{j_{1}}z_{w(1)}z_{w(2)}\ldots
z_{w(p_{1})}=g_{w|C_{1}}^{j_{1}}z_{w|C_{1}},
\]
and we have corresponding factors for $C_{2},\ldots ,C_{q}$. The relevant
summand of $h_{w}$ is given by 
\[
f_{j_{1}\ldots j_{q}}g_{w|C_{1}}^{j_{1}}\ldots g_{w|C_{q}.}^{j_{q}} 
\]
We conclude that 
\[
h_{w}=\sum_{q}\sum_{\pi =(C_{k})\in \mathcal{Y}_{q}([p])}f_{j_{1}\ldots
j_{q}}g_{w|C_{1}}^{j_{1}}\ldots g_{w|C_{q}}^{j_{q}}. 
\]
More generally we may substitute $G$ into an $N$-tuple $F(z)=(F^{1}(z),%
\ldots ,F^{n}(z)),$ obtaining $H=F\circ G$, where 
\[
h_{w}^{i}=\sum f_{j_{1}\ldots j_{q}}^{i}g_{w|C_{1}}^{j_{1}}\ldots
g_{w|C_{q}}^{j_{q}}, 
\]
where where we sum over all interval colored partitions 
\[
((C_{1},\ldots ,C_{q}),j_{1}\cdots j_{q})\,\,\,1\leq q\leq p 
\]
of the colored set $([p],u(1)\cdots u(p))$.

Given an algebra $A,$ we let $\mathcal{G}_{N}^{dif}(A)$ denote the set of
power series $F=F(z)$ with 
\[
F^{j}(z)=z_{j}+\sum_{\left| u\right| \geq
2}f_{u}^{j}z_{u},\,\,\,(f_{u}^{j}\in A) 
\]
i.e., without constant terms and with $f_{j}^{i}=\delta _{i}^{j}.$
Substitution of $G$ into $F$ provides us with a non-associative product $%
(F,G)\mapsto F\circ G$ on $\mathcal{G}_{N}^{dif}(A)$. From above, 
\begin{equation}
(F\circ G)_{u}^{i}=z_{i}+\sum f_{w}^{i}g_{u|C_{1}}^{w(1)}\cdots
g_{u|C_{q}}^{w(q)}  \label{substitution}
\end{equation}
where we sum over all interval colored partitions 
\[
((C_{1},\ldots ,C_{q}),w(1)\cdots w(q))\,\,\,1\leq q\leq p 
\]
of the colored set $([p],u(1)\cdots u(p))$.

Each generator $Y_{u}^{i}\in \mathcal{L}^{N}$ $(\left| u\right| >1)$ may be
used to select a corresponding coefficent $f_{u}^{i}$ in a power series $%
F(z).$ To be more precise, we define a linear mapping 
\[
\theta (Y_{u}^{i}):\mathcal{G}_{N}^{dif}(A)\rightarrow A 
\]
by letting $\theta (Y_{u}^{i})(F)=f_{u}^{i}.$ Since $\mathcal{L}^{N}$ is the
free algebra on these generators, we extend this to the basis elements $%
Y_{u_{1}}^{i_{1}}\cdots Y_{u_{q}}^{i_{q}}$ by letting 
\[
\theta (Y_{u_{1}}^{i_{1}}\cdots Y_{u_{q}}^{i_{q}}):\mathcal{G}
_{N}^{dif}(A)\rightarrow A:F\mapsto f_{u_{1}}^{i_{1}}\cdots
f_{u_{q}}^{i_{q}} 
\]
Extending linearly, we have a corresponding homomorhism 
\[
\theta :\mathcal{L}^{N}\rightarrow Lin(\mathcal{G}_{N}^{dif}(A),A) 
\]
and thus a bilinear mapping 
\[
\langle ,\rangle :\mathcal{L}^{N}\times \mathcal{G}_{N}^{dif}\rightarrow
A:(a,F)\mapsto \theta (a)(F). 
\]
From our definitions we have that 
\[
\langle ab,f\rangle =\langle a,f\rangle \langle b,f\rangle =m_{A}\langle
a\otimes b,f\otimes f\rangle . 
\]

Returning to the
substitution formula (\ref{substitution}), we have 
\begin{eqnarray*}
\langle Y_{u}^{i},F\circ G\rangle &=&m_{A}\sum_{w} \langle Y_{w}^{i}\otimes
Y_{u|C_{1}}^{w(1)}\cdots Y_{u|C_{s}}^{w(s)},F\otimes G\rangle \\
&=&m_{A}\langle \Delta ^{op }(Y_{u}^{i}),F\otimes G\rangle .
\end{eqnarray*}

\begin{theorem}
Given $F\in \mathcal{G}_{N}^{dif}(A),$ the left substitutional inverse of $F$
is given by the power series $G(z),$ where $g_{v}^{j}=\langle
S_{{\mathcal L}}(Y_{v}^{j}),f\rangle $ where $S_{{\mathcal
L}}=\mathbf{s}S_{\mathcal{H}}\mathbf{s}$ is the antipode of the left Lagrange Hopf
algebra $\mathcal{L}^{N}$. The right substitutional inverse of $F$ is given by
$H(z)$, where $h_{v}^{j}=\langle  S_{{\mathcal R}}(Y_{v}^{j}),f\rangle ,$ and
$S_{{\mathcal R}}=%
\mathbf{t}S_{\mathcal{H}}\mathbf{t}$ is
the antipode of the right Legendre Hopf algebra $\mathcal{R}^{N}$.
\end{theorem}

\proof Let $m_{A}$ denote the multiplication in $A.$ Defining $G$ as above,
we have that 
\begin{eqnarray*}
\langle Y_{u}^{j},G\circ F\rangle &=&m_{A}\langle \Delta ^{op
}(Y_{u}^{j}),G\otimes F\rangle \\
&=&m_{A}(\sum \langle Y_{w}^{i}\otimes Y_{u|C_{1}}^{w(1)}\cdots
Y_{u|C_{s}}^{w(s)},G\otimes F\rangle ) \\
&=&\sum \langle S_{\mathcal {H}}^{-1}(Y_{w}^{j}),F\rangle \langle
Y_{u|C_{1}}^{w(1)}\cdots Y_{u|C_{s}}^{w(s)},F\rangle \\
&=&\sum \langle S_{\mathcal
{H}}^{-1}((Y_{u}^{j})_{(2)})(Y_{u}^{j})_{(1)},F\rangle \\ &=&\langle
\varepsilon (Y_{u}^{j})1,F\rangle \\ &=&\delta _{u}^{j}1,
\end{eqnarray*}
and thus $F\ $is the left substitutional inverse of $G.$

On the other hand, if $H$ is defined as above, then using the fact that $%
\mathbf{t}S$ is an algebraic homomorphism and that $\mathbf{t}
(Y_{w}^{i})=Y_{w}^{i}$, 
\begin{eqnarray*}
\langle Y_{u}^{j},F\circ H\rangle &=&m_{A}\langle \Delta ^{op
}(Y_{u}^{j}),F\otimes H\rangle \\
&=&m_{A}(\sum \langle Y_{w}^{i}\otimes Y_{u|C_{1}}^{w(1)}\cdots
Y_{u|C_{s}}^{w(s)},F\otimes H\rangle ) \\
&=&m_{A}(\sum \langle Y_{w}^{i}\otimes (\mathbf{t}S)(Y_{u|C_{1}}^{w(1)})\cdots
(\mathbf{t}S)(Y_{u|C_{s}}^{w(s)})),F\otimes F\rangle ) \\
&=&m_{A}(\sum \langle Y_{w}^{i}\otimes (\mathbf{t}S)(Y_{u|C_{1}}^{w(1)})\cdots
Y_{u|C_{s}}^{w(s)}),F\otimes F\rangle \\
&=&\sum \langle Y_{w}^{j}(\mathbf{t}S)(Y_{u|C_{1}}^{w(1)}\cdots
Y_{u|C_{s}}^{w(s)})),F\rangle \\
&=&\langle \mathbf{t}(\sum S(Y_{u|C_{1}}^{w(1)}\cdots
Y_{u|C_{s}}^{w(s)})Y_{w}^{j}),F\rangle \\
&=&\langle \mathbf{t}(\sum S(Y_{u}^{j})_{(1)}(Y_{u}^{j})_{(2)}),F\rangle \\
&=&\langle \varepsilon (Y_{u}^{j})1,F\rangle \\
&=&\delta _{u}^{j}1
\end{eqnarray*}
and $H$ is the right substitutional inverse of $F.$
\endproof

\begin{corollary}
If the number of variables $N$ is greater than $1$, then the left and right
substitutional inverses of a power series are generally distinct.
\end{corollary}

\proof It suffices to show that 
\[
\mathbf{s}\circ S_{\mathcal{H}}\circ \mathbf{s}(Y_{1234}^{1})\neq \mathbf{t}\circ
S_{\mathcal{H}}\circ
\mathbf{t} (Y_{1234}^{1}) 
\]
In the following calculation we have used boldface subscripts to indicate
corresponding terms that equal. The bracketed terms cancel (these correspond
to the unique order preserving contraction). The sums are over the set of
colors, i.e., $k,\ell =1,2,3,4$ and they involve 13 varieties of layered
trees. 
\begin{eqnarray*}
S(Y_{1234}^{1}) &=&-Y_{1234}^{1}+\sum Y_{123}^{k}Y_{k4}^{1}+\sum
Y_{234}^{k}Y_{1k}^{1}+\sum Y_{12}^{k}Y_{k34}^{1} \\
&&\hspace{-.8in}+\sum Y_{23}^{k}Y_{1k4}^{1}+\sum Y_{34}^{k}Y_{12k}^{1}
+\left[ \sum Y_{12}^{k}Y_{34}^{\ell }Y_{k\ell }^{1}\right] -\sum
Y_{34}^{k}Y_{12}^{\ell }Y_{\ell k}^{1}-\left[ \sum Y_{12}^{k}Y_{34}^{\ell
}Y_{k\ell }^{1}\right] \\
&&-\sum Y_{12}^{k}Y_{k3}^{\ell }Y_{\ell 4}^{1}-\sum Y_{23}^{k}Y_{1k}^{\ell
}Y_{\ell 4}^{1}-\sum Y_{23}^{k}Y_{k4}^{\ell }Y_{1\ell }^{1}-\sum
Y_{34}^{k}Y_{2k}^{\ell }Y_{1\ell }^{1} \end{eqnarray*}
\begin{eqnarray*}
\mathbf{s}S\mathbf{s}(Y_{1234}^{1}) &=&-Y_{1234}^{1}+_{\mathbf{p}}\sum
Y_{k4}^{1}Y_{123}^{k}\,+\,\,_{\mathbf{q}}\sum Y_{1k}^{1}Y_{234}^{k}+\,_{%
\mathbf{a}}\sum Y_{k34}^{1}Y_{12}^{k}+ \\
&&\hspace{-.2in}+_{\mathbf{b}}\sum Y_{1k4}^{1}Y_{23}^{k} +_{\mathbf{c}}\sum
Y_{12k}^{1}Y_{34}^{k} -\sum Y_{\ell k}^{1}Y_{12}^{\ell }Y_{34}^{k} \\
&&\hspace{-.2in}-_{\mathbf{d}}\sum Y_{\ell 4}^{1}Y_{k3}^{\ell }Y_{12}^{k}-_{%
\mathbf{e}}\sum Y_{\ell 4}^{1}Y_{1k}^{\ell }Y_{23}^{k}-_{\mathbf{f}}\sum
Y_{1\ell }^{1}Y_{k4}^{\ell }Y_{23}^{k}-_{\mathbf{g}}\sum Y_{1\ell
}^{1}Y_{2k}^{\ell }Y_{34}^{k} \\
\mathbf{t}(Y_{1234}^{1}) &=&Y_{4321}^{1} \end{eqnarray*}
\begin{eqnarray*}
S\mathbf{t}(Y_{1234}^{1}) &=&S(Y_{4321}^{1})=-Y_{4321}^{1}+\sum
Y_{432}^{k}Y_{k1}^{1}+\sum Y_{321}^{k}Y_{4k}^{1}+\sum Y_{43}^{k}Y_{k21}^{1}
\\
&&+\sum Y_{32}^{k}Y_{4k1}^{1}+\sum Y_{21}^{k}Y_{43k}^{1} -\sum
Y_{21}^{k}Y_{43}^{\ell }Y_{\ell k}^{1} \\
&&-\sum Y_{43}^{k}Y_{k2}^{\ell }Y_{\ell 1}^{1}-\sum Y_{32}^{k}Y_{4k}^{\ell
}Y_{\ell 1}^{1}-\sum Y_{32}^{k}Y_{k1}^{\ell }Y_{4\ell }^{1}-\sum
Y_{21}^{k}Y_{3k}^{\ell }Y_{4\ell }^{1}\end{eqnarray*} 
\begin{eqnarray*}
\mathbf{t}S\mathbf{t}(Y_{1234}^{1}) &=&-Y_{1234}^{1}+_{\mathbf{q}}\sum
Y_{1k}^{1}Y_{234}^{k}+_{\mathbf{p}}\sum Y_{k4}^{1}Y_{123}^{k}+_{\mathbf{c}
}\sum Y_{12k}^{1}Y_{34}^{k} \\
&&\hspace{-.2in}+_{\mathbf{b}}\sum Y_{1k4}^{1}Y_{23}^{k}+_{\mathbf{a}}\sum
Y_{k34}^{1}Y_{12}^{k} -\sum Y_{k\ell }^{1}Y_{34}^{\ell }Y_{12}^{k} \\
&&\hspace{-.2in}-_{\mathbf{g}}\sum Y_{1\ell }^{1}Y_{2k}^{\ell }Y_{34}^{k}-_{%
\mathbf{f}}\sum Y_{1\ell }^{1}Y_{k4}^{\ell }Y_{23}^{k}-_{\mathbf{e}}\sum
Y_{\ell 4}^{1}Y_{1k}^{\ell }Y_{23}^{k}-_{\mathbf{d}}\sum Y_{\ell
4}^{1}Y_{k3}^{\ell }Y_{12}^{k}
\end{eqnarray*}

It follows that 
\[
\,\mathbf{s}\,S\mathbf{s}(Y_{1234}^{1})-\mathbf{t}S\mathbf{t}(Y_{1234}^{1})=-\sum
Y_{\ell k}^{1}Y_{12}^{\ell }Y_{34}^{k}+\sum Y_{k\ell }^{1}Y_{34}^{\ell
}Y_{12}^{k}. 
\]
\endproof

In particular, one can check that if $a$ and $b$ and $x$ and $y$ do not
commute, the substitional left and right inverses of the two-variable
polynomial function 
\begin{eqnarray*}
u &=&x+ax^{2}+by^{2} \\
v &=&y
\end{eqnarray*}
do not agree in the fourth order terms.

\end{document}